\documentclass[11pt]{amsart}

\usepackage{latexsym,graphicx}
\numberwithin{equation}{section}
\theoremstyle{plain}

\theoremstyle{remark}

\theoremstyle{definition}

\newcommand{\D}{{\mathcal D}}
\newcommand{\E}{\mathcal E}

\newcommand{\G}{{\mathcal G}}

\renewcommand{\H}{\mathbb H}

\newcommand{\K}{{\mathcal K}}
\renewcommand{\L}{{\mathcal L}}
\newcommand{\M}{{\mathcal M}}
\newcommand{\N}{\mathbb N}

\newcommand{\V}{{\mathcal V}}

\newcommand{\dist}{\operatorname{dist}}

\newcommand{\fp}{\operatorname{FP}}

\newcommand{\Int}{\operatorname{Int}}

\renewcommand{\span}{\operatorname{span}}

\newcommand{\supp}{\operatorname{Supp}}

\def\half{{1 \over 2}}

\newcommand{\oa}{\overrightarrow}

\newcommand{\ol}{\overline}

\def\XXint#1#2#3{{\setbox0=\hbox{$#1{#2#3}{\int}$}
      \vcenter{\hbox{$#2#3$}}\kern-.5\wd0}}

\begin{document}

\def\cal{\mathcal}

\font\tpt=cmr10 at 12 pt
\font\fpt=cmr10 at 14 pt

\font \fr = eufm10

%\font\AAA=Times.dfont  at 12pt
 %\font\BBB=Times.dfont at 8pt

%\font\AAA=cmr10 at 12pt
%\font\BBB=cmr10 at 8pt

%\def\AAA{\bf}
%\def\BBB{\bf}

\overfullrule=0in

\def\boxit#1{\hbox{\vrule
 \vtop{%
  \vbox{\hrule\kern 2pt %
     \hbox{\kern 2pt #1\kern 2pt}}%
   \kern 2pt \hrule }%
  \vrule}}

  \def\harr#1#2{\ \smash{\mathop{\hbox to .3in{\rightarrowfill}}\limits^{\scriptstyle#1}_{\scriptstyle#2}}\ }

\def\AA{1}
\def\BB{2}
\def\CC{3}
\def\DD{4}
\def\EE{5}
\def\FF{6}
\def\GGG{7}
\def\HH{8}
\def\II{9}
\def\JJ{10}
\def\KK{11}
\def\LL{12}
\def\MM{13}

 \def\GG{{{\bf G} \!\!\!\! {\rm l}}\ }

\def\GL{{\rm GL}}

\def\bll{I \!\! L}

\def\IFF{\qquad\iff\qquad}
\def\bra#1#2{\langle #1, #2\rangle}
\def\bbf{{\bf F}}
\def\bbj{{\bf J}}
\def\Jtn{{\bbj}^2_n}  \def\JtN{{\bbj}^2_N}  \def\JoN{{\bbj}^1_N}
\def\jt{j^2}
\def\jtx{\jt_x}
\def\Jt{J^2}
\def\Jtx{\Jt_x}
\def\bpp{{\bf P}^+}
\def\bpt{{\wt{\bf P}}}
\def\fsh{$F$-subharmonic }
\def\mo{monotonicity }
\def\jet{(r,p,A)}
\def\ss{\subset}
\def\sse{\subseteq}
\def\half{\hbox{${1\over 2}$}}
\def\smfrac#1#2{\hbox{${#1\over #2}$}}
\def\oa#1{\overrightarrow #1}
\def\dim{{\rm dim}}
\def\dist{{\rm dist}}
\def\codim{{\rm codim}}
\def\deg{{\rm deg}}
\def\rank{{\rm rank}}
\def\log{{\rm log}}
\def\Hess{{\rm Hess}}
\def\Hessyp{{\rm Hess}_{\rm SYP}}
\def\trace{{\rm trace}}
\def\tr{{\rm tr}}
\def\max{{\rm max}}
\def\min{{\rm min}}
\def\span{{\rm span\,}}
\def\Hom{{\rm Hom\,}}
\def\det{{\rm det}}
\def\End{{\rm End}}
\def\Sym{{\rm Sym}^2}
\def\diag{{\rm diag}}
\def\pt{{\rm pt}}
\def\Spec{{\rm Spec}}
\def\pr{{\rm pr}}
\def\Id{{\rm Id}}
\def\Grass{{\rm Grass}}
\def\Herm#1{{\rm Herm}_{#1}(V)}
\def\arr{\longrightarrow}
\def\supp{{\rm supp}}
\def\Link{{\rm Link}}
\def\Wind{{\rm Wind}}
\def\Div{{\rm Div}}
\def\vol{{\rm vol}}
\def\foral{\qquad {\rm for\ all\ \ }}
\def\fpsh{{\cal PSH}(X,\f)}
\def\Core{{\rm Core}}
\def\dis{f_M}
\def\Re{{\rm Re}}
\def\rn{\bbr^n}
\def\pp{\cp^+}
\def\plp{\cp_+}
\def\Int{{\rm Int}}
\def\cix{C^{\infty}(X)}
\def\Gr#1{G(#1,\rn)}
\def\Symn{{\Sym(\rn)}}
\def\SymN{{\Sym(\bbr^N)}}
\def\Gpn{G(p,\rn)}
\def\fd{{\rm free-dim}}
\def\SA{{\rm SA}}
 \def\cd{{\cal C}}
 \def\cdt{{\widetilde \cd}}
 \def\cm{{\cal M}}
 \def\cmt{{\widetilde \cm}}

\def\Theorem#1{\medskip\noindent {\bf THEOREM \bf #1.}}
\def\Prop#1{\medskip\noindent {\bf Proposition #1.}}
\def\Cor#1{\medskip\noindent {\bf Corollary #1.}}
\def\Lemma#1{\medskip\noindent {\bf Lemma #1.}}
\def\Remark#1{\medskip\noindent {\bf Remark #1.}}
\def\Note#1{\medskip\noindent {\bf Note #1.}}
\def\Def#1{\medskip\noindent {\bf Definition #1.}}
\def\Claim#1{\medskip\noindent {\bf Claim #1.}}
\def\Conj#1{\medskip\noindent {\bf Conjecture \bf    #1.}}
\def\Ex#1{\medskip\noindent {\bf Example \bf    #1.}}
\def\Qu#1{\medskip\noindent {\bf Question \bf    #1.}}
\def\Exercise#1{\medskip\noindent {\bf Exercise \bf    #1.}}

\def\HoQu#1{ {\AAA T\BBB HE\ \AAA H\BBB ODGE\ \AAA Q\BBB UESTION \bf    #1.}}

\def\pf{\medskip\noindent {\bf Proof.}\ }
\def\qed{\hfill  $\vrule width5pt height5pt depth0pt$}
\def\equdef{\buildrel {\rm def} \over  =}
\def\qedqed{\hfill  $\vrule width5pt height5pt depth0pt$ $\vrule width5pt height5pt depth0pt$}
\def\mathqed{  \vrule width5pt height5pt depth0pt}

\def\V{W}

\def\df{d^{\phi}}
\def\hk{\_{\rm l}\,}
\def\n{\nabla}
\def\w{\wedge}

\def\cu{{\cal U}}   \def\cc{{\cal C}}   \def\cb{{\cal B}}  \def\cz{{\cal Z}}
\def\cv{{\cal V}}   \def\cp{{\cal P}}   \def\ca{{\cal A}}
\def\cw{{\cal W}}   \def\co{{\cal O}}
\def\ce{{\cal E}}   \def\ck{{\cal K}}
\def\ch{{\cal H}}   \def\cm{{\cal M}}
\def\cs{{\cal S}}   \def\cn{{\cal N}}
\def\cd{{\cal D}}
\def\cl{{\cal L}}
\def\cp{{\cal P}}
\def\cf{{\cal F}}
\def\ccr{{\cal  R}}

\def\gerG{{\fr{\hbox{g}}}}
\def\gerB{{\fr{\hbox{B}}}}
\def\gerR{{\fr{\hbox{R}}}}
\def\p#1{{\bf P}^{#1}}
\def\vf{\varphi}

\def\wt{\widetilde}
\def\wh{\widehat}

\def\and{\qquad {\rm and} \qquad}
\def\arr{\longrightarrow}
\def\ol{\overline}
\def\bbr{{\mathbb R}}\def\bbh{{\mathbb H}}\def\bbo{{\mathbb O}}
\def\bbc{{\mathbb C}}
\def\bbq{{\mathbb Q}}
\def\bbz{{\mathbb Z}}
\def\bbp{{\mathbb P}}
\def\bbd{{\mathbb D}}

\def\a{\alpha}
\def\b{\beta}
\def\d{\delta}
\def\e{\epsilon}
\def\f{\phi}
\def\g{\gamma}
\def\k{\kappa}
\def\l{\lambda}
\def\o{\omega}

\def\s{\sigma}
\def\x{\xi}
\def\z{\zeta}

\def\D{\Delta}
\def\L{\Lambda}
\def\G{\Gamma}
\def\O{\Omega}

\def\bd{\partial}
\def\bdf{\partial_{\f}}
\def\lag{Lagrangian}
\def\psh{plurisubharmonic }
\def\ph{pluriharmonic }
\def\pph{partially pluriharmonic }
\def\omp{$\omega$-plurisubharmonic \ }
\def\ffl{$\f$-flat}
\def\PH#1{\widehat {#1}}
\def\lloc{L^1_{\rm loc}}
\def\dbar{\ol{\partial}}
\def\lp{\Lambda_+(\f)}
\def\lpp{\Lambda^+(\f)}
\def\bo{\partial \Omega}
\def\Ob{\overline{\O}}
\def\fc{$\phi$-convex }
\def\PSH{{ \rm PSH}}
\def\SH{{\rm SH}}
\def\totr{ $\phi$-free }
\def\BM{\lambda}
\def\Der{D}
\def\CH{{\cal H}}
\def\RH{\overline{\ch}^\f }
\def\pconv{$p$-convex}
\def\MA{MA}
\def\lagpsh{Lagrangian plurisubharmonic}
\def\hermsk{{\rm Herm}_{\rm skew}}
\def\PSHl{\PSH_{\rm Lag}}
 \def\ppsh{$\pp$-plurisubharmonic}
\def\fp{$\pp$-plurisubharmonic }
\def\fh{$\pp$-pluriharmonic }
\def\Symn{\Sym(\rn)}
 \def\ci{C^{\infty}}
\def\USC{{\rm USC}}
\def\LSC{{\rm LSC}}
\def\fa{{\rm\ \  for\ all\ }}
\def\ppc{$\pp$-convex}
\def\cpt{\wt{\cp}}
\def\ft{\wt F}
\def\ob{\overline{\O}}
\def\Be{B_\e}
\def\K{{\rm K}}

\def\M{{\bf M}}
\def\N#1{C_{#1}}
\def\ds{Dirichlet set }
\def\dir{Dirichlet }
\def\Fa{{\oa F}}
\def\TR{{\cal T}}
 \def\ISO{{\rm ISO_p}}
 \def\Span{{\rm Span}}

\def\ALL{1}
\def\AV{2}
\def\BTA{3}
\def\BL{4}
\def\BRE{5}
\def\CNS{6}
\def\CP{7}
\def\CPW{8}
\def\CIL{9}
\def\CRA{10}
\def\DON{11}
\def\DDD{12}
\def\DDR{13}
\def\GEO{14}
\def\HYP{15}
\def\BEL{16}
\def\SURVEY{17}
\def\AC{18}
\def\NOTES{19}
\def\AET{20}
\def\TANG{21}
\def\TANGG{22}
\def\LAG{23}
\def\JTY{24}
\def\KRY{25}
\def\PLI{26}
\def\RT{27}
\def\SLO{28}
\def\SPR{29}
\def\TRUU{30}
\def\TRU{31}
\def\TWC{32}
\def\TWCC{33}
\def\TWCCC{34}
\def\WAL{35}

\def\CHLP{1}
\def\CIL{2}
\def\CRA{3}
\def\DHGL{4}
\def\Eb{5}
\def\PTCG{6}
\def\PTCGG{7}
\def\DD{8}
\def\DDR{9}
\def\Survey{10}
\def\Hyp{11}
\def\Geo{12}
\def\NOTES{13}
\def\Gar{14}
\def\Pot{15}
\def\Lag{16}
\def\SMP{17}
\def\IDP{18}
\def\Edge{19}
\def\HP{20}
\def\HU{21}
\def\Kry{22}
\def\GRU{23}
\def\ST{24}
\def\SW{25}
\def\RU{26}

\vskip .4in

\def\E{E}
\def\fpsi{{F_f(\psi)}}
\def\bL{{\bf \Lambda}}
\def\bdf{{\bf f}}
\def\UU{U}
\def\bbm{{\bf M}}

\def\bbf{{\mathbb F}}
\def\bbe{{\mathbb E}}
\def\bbg{{\mathbb G}}

\font\headfont=cmr10 at 14 pt
\font\aufont=cmr10 at 11 pt

\def\bg{{\bf G}}

\def\bbem{\bbe_{\rm min}}
\def\bbgm{\bbg_{\rm max}}
\def\bbhc{\bbh_{\rm can}}
\def\bbgmt{\wt{\bbg}_{\rm max}}

\def\bbed{\bbe^{\diamondsuit}}
\def\bbgd{\bbg^{\diamondsuit}}
\def\bbhd{\bbh^{\diamondsuit}}

\title[A GENERALIZATION OF PDE'S FROM  A KRYLOV POINT OF VIEW]
{\headfont A GENERALIZATION OF PDE'S  FROM \\  A KRYLOV POINT OF VIEW}

\date{\today}
\author{ F. Reese Harvey and H. Blaine Lawson, Jr.}
\thanks
{Partially supported by the NSF}

\maketitle
\begin{abstract}
We introduce and investigate  the notion of a ``generalized equation'',
which  extends  nonlinear elliptic equations, and which 
is  based on the notions of  subequations and Dirichlet duality.  Precisely,
a subset $\bbh\ss \Symn$ is a {\bf generalized equation} if it is an intersection 
$\bbh = \bbe \cap (-\wt \bbg)$ where $\bbe$ and $\bbg$ are subequations and 
$\wt \bbg$ is the subequation dual to $\bbg$.  We utilize a viscosity  definition of ``solution'' 
to $\bbh$.  The $\bf mirror$ of $\bbh$ is defined by $\bbh^* \equiv \bbg\cap (-\wt \bbe)$.
One of the main results (Theorem 2.6) concerns the Dirichlet problem on arbitrary bounded domains
$\O\ss \rn$ for solutions to $\bbh$ with prescribed boundary function $\vf \in C(\bo)$.  We prove  that:

(A) Uniqueness holds  \qquad $\iff$ \qquad $\bbh$ has no interior, and

(B) Existence holds \ \    \qquad $\iff$ \qquad $\bbh^*$ has no interior.

\noindent
For (B) the appropriate boundary convexity of $\bo$ must be assumed.  Many examples of generalized 
equations are discussed, including the constrained Laplacian, the twisted Monge-Amp\`ere equation, and
the $C^{1,1}$-equation.

The closed sets $\bbh\ss\Symn$ which can be written as  generalized equations are intrinsically characterized.
For such an $\bbh$ the set of subequation pairs $(\bbe, \bbg)$ with $\bbh = \bbe \cap (-\wt \bbg)$ 
is partially ordered (see (1.10)). If $(\bbe,\bbg)\prec (\bbe', \bbg')$, then any solution for the first is
also a solution for the second.   Furthermore, in this ordered set there is a canonical least element, contained in all others.

A general form of the main theorem, which holds on any manifold, is also established.

\end{abstract}

\vfill\eject

\ \ 
\vskip.5in

\centerline{\headfont Table of Contents}

1.  \ \ Introduction 

2. \ \ Main Definitions

3. \ \ The Canonical Pair Defining a Given $\bbh$,
 
\qquad \ \ \ \  and an Intrinsic  Characterization of Generalized Equations

4. \ \ Examples

5. \ \ Unsettling Questions

6. \ \  A General Case of the Main Theorem

  Appendix  A.  The Quasi-Convexity Characterization of C$^{1,1}$.

%\vfill\eject
\vskip.3in

\noindent
{\headfont  1. Introduction.}

For some time  we have been studying fully nonlinear pde's from a perspective 
of generalized potential theory.  This was initiated by the discovery [\PTCG,\PTCGG] of a robust potential
theory attached to every calibrated manifold  --   a fact which generalized the classical
pluripotential theory in   K\"ahler geometry.  Our point of view had some reflections in
early work of Krylov [\Kry], and eventually went far outside calibrated geometries.
Good references for these techniques and results are  [\DD], [\DDR] and [\Survey]

The purpose of this paper is to examine, to the fullest extent, when these viscosity and
Dirichlet duality techniques can be employed to study  nonlinear differential relations.
Our fundamental Definition 2.2 is completely natural in the context of duality, but may seem
 to be too general or abstract.  However,  it turns out that there are lots of interesting 
examples: the constrained Laplace equation, the twisted Monge-Amp\`ere equation, the relation
$\det(D^2 u)\leq 0$ in dimension 2, the $C^{1,1}$-equation, and many, many more.
Furthermore,  our fundamental Theorem 2.6 gives a simple and somewhat surprising 
relationship with the two questions of existence and uniqueness.

For clarity and simplicity we restrict attention to the constant coefficient case until the last Section 6.

Before making more detailed comments about examples and results, we recall
the fundamental background ideas.

\vskip .3in 

\centerline{\bf Short Preliminaries}

We adopt the subequation point of view from  [\DD], [\DDR], where
a differential operator $f$ is replaced by the constraint set 
$\bbf\equiv \{A\in \Symn ; f(A)\geq 0\}$. (Here $\Symn$ denotes the space of quadratic forms
on $\rn$.)  The equation  $f(D^2 u)=0$ is replaced by the constraint condition 
$D^2 u\in \partial \bbf$. For that reason we will refer to $\partial\bbf$ as the 
{\bf equation associated to the subequation $\bbf$.}
 The ellipticity hypothesis  is assumed in the weakest possible form:  
$$
\bbf + \cp \ \ss\ \bbf.
\eqno{(1.1)}
$$
Here $\cp = \{P\in\Symn : P\geq 0\}$.  Any closed subset $\bbf\ss\Symn$ satisfying this {\sl positivity},
or $\cp$-{\sl monotonicty, condition} (1.1), is called a {\bf subequation}.\footnote{It is convenient occasionally in this paper
to allow $\bbf=\emptyset$ or $\bbf = \Symn$.}

The viscosity definition of a subsolution takes the following form.  Consider  an upper semi-continuous
function $u$ defined on an open set $X\ss \rn$ and taking values in $[-\infty, \infty)$.
An {\sl upper test function} for $u$ at a point $x\in X$ is a $C^2$-function $\vf$ defined  
near $x$ with $u\leq \vf$ and $u(x)=\vf(x)$.  The function $u$ is said to be $\bbf$-{\bf subharmonic}
or an  $\bbf$-{\bf subsolution} on $X$  if for every upper test function $\vf$ at any point $x\in X$ we have
$D^2_x \vf \in \bbf$.  For $C^2$-functions $u$, the {\bf consistency} of this definition with the classical
definition that $D^2_x u \in \bbf$ for all $x$, follows from the positivity condition (1.1).
(In fact, consistency mandates positivity.)
We will denote the space of  $\bbf$-subharmonics functions on $X$ by $\bbf(X)$.
(For a complete introduction to viscosity theory  see [\CIL, \CRA].)

The  {\bf Dirichlet dual} $ \wt \bbf$ of a subequation $\bbf$ is defined to be
$$
 \wt \bbf \ =\ \sim(-\Int\, \bbf) \ =\ -(\sim \Int \bbf)
\eqno(1.2)
$$
It is also a subequation and provides a true duality $\wt{\wt \bbf} = \bbf$.
Moreover, one has the key relationship
$$
\partial \bbf \ =\ \bbf \cap (- \wt \bbf)
\eqno(1.3)
$$
which enables one to replace $f(D^2 u)=0$ by $\partial \bbf$ via the  viscosity definitions
(see Def. 2.1).
It is easy to see that 
$$
\Int  \wt \bbf \ =\ -(\sim \bbf)\ =\ \sim(-\bbf),
\eqno(1.4)
$$
and to see that (1.1), together with $\bbf$ being closed, implies the 
{\bf topological property}
$$
\bbf \ =\\ \overline{\Int\,  \bbf}.
\eqno(1.5)
$$
(Also $\Int \bbf, \partial \bbf$ and $\bbf$ are all path-connected.)

By (1.3) applied to $\wt\bbf$, instead of $\bbf$,  we have
$$
\partial \wt\bbf \ =\ \wt\bbf \cap (-\bbf) \ =\ -\partial \bbf.
\eqno(1.6)
$$
It is easy to see that 
$$
\bbf \ =\ \partial \bbf +\cp
\eqno(1.7)
$$
and hence, by (1.6) we have
$$
\wt\bbf\ =\ -\partial \bbf + \cp
\eqno(1.8)
$$
This formula for the dual subequation  $\wt\bbf$ could have replaced (1.2) as the definition of the dual subequation.
It is frequently the easiest way to compute $\wt\bbf$ for 
examples, as in Section 4.

\vskip.2in

\centerline{\bf Generalized Equations}

A  {\bf   generalized equation} is simply a pair $(\bbe, \bbg)$ of subequations on $\Symn$.
We associate to this pair the constraint set
$$
\bbh \ =\ \bbh_{\bbe, \bbg} \ \equdef \ \bbe \cap(-\wt \bbg)
$$
and define an {\bf $\bbh$-harmonic} on an open $X\ss\rn$ to be a  function $h$  on $X$ with
$$
h \in \bbe(X) \and  -h\in \wt\bbg(X),
$$
i.e., $h$ is $\bbe$-subharmonic and $-h$ is $\wt\bbg$-subharmonic on $X$.  
An  $\bbh$-harmonic  on $X$ which is $C^2$ satisfies the differential relation
$$
D^2_x h \in \bbh \qquad\text{for all } \ x\in X
\eqno{(1.9)}
$$
(by the $C^2$ consistency referred to above).
One example is the {\bf constrained Laplacian} (see Example 2.3 below)  
 $$\bbh = \{\tr A=0 \  {\rm and}\  -r\Id \leq A\leq r \Id\} \qquad {\rm for }\ \ r\geq 0,$$  
where $\bbe=\bbh+\cp$ and  $-\wt \bbg = \bbh-\cp$, i.e.,  $\bbg=\wt \bbe$.
 Here the $\bbh$-harmonics are classical harmonic functions $h$ with 
$-r\Id \leq D^2 h\leq r\Id$.  This example can be generalized with $\bbh$ being any closed subset
of $\{\tr A=0\}$ (see Example 4.3).

Another Example is the twisted Monge-Amp\`ere Equation (Example 4.8) which has been studied
by Streets, Tian and Warren, [\ST,\SW].  This equation requires a splitting of space, and uses a mixture of a convex and concave Monge Ampere operator on the two pieces.  The authors proved an Evans-Krylov type estimate despite the nonconvexity of the operator.  We are grateful for Jeff Streets communicating their work on this equation and asking us if any of our methods could apply.

A very nice example $\bbh_\l$ is given by taking  $\bbe =\cp - \l\Id$ and $-\wt \bbg= -\bbe = -\cp+\l\Id$ 
for $\l\geq 0$. This is a special case of the 
quasi-convex/quasi-concave equation in Example 4.1. 
By a result of   Hiriart-Urruty and Plazanet  (which we present in Appendix A) we have 
$$
h\ \ \text{is $\bbh_\l$-harmonic \ \ $\iff$\ \ $h$ is $C^{1,1}$ with Lipschitz coefficient $\l$.}
$$

Related to this is the quasi-subaffine/quasi-superaffine equation in Example 4.2.  This equation
is the mirror $\bbh_\l^*$ of the $\bbh_\l$ above.  Let's let $\l=1$ and set $\bbh=\bbh_1$.  Then the intersection
$\bbh\cap \bbh^*$ is another generalized equation (see Remark before (2.4)).  The $\bbh\cap \bbh^*$-harmonics
$h$ of class $C^2$ satisfy
$$
D^2h +\Id \geq 0, \ \ \Id - D^2h \geq 0 \quad{\rm and}\quad \det(D^2h +\Id) =- \det(\Id -D^2h).
$$

A basic example is given by taking  $\bbe =\cpt$ and $\bbg=\cp$, i.e., $-\wt\bbg = -\cpt$ 
so that $\bbh = \cpt\cap (-\cpt) = \Symn -\Int \{\cp \cup (-\cp)\}$.  This is a special case of the 
quasi-subaffine/quasi-superaffine equation in Example 4.2 (see also Example 4.16(d)).  
When $n=2$, we see that
$\bbh = \{A : \l_1(A) \l_2(A)  =\det(A) \leq 0\}$.

When $\bbg = \bbe$, we have  $\bbh = \partial \bbe$ and we are in the case discussed
in the preliminaries above.  Here there is no other pair giving the same $\bbh$.  
%However for other $\bbh$'s there can be many pairs.

Note however that a generalized equation $\bbh$ can possibly be written as $\bbe\cap (-\wt \bbg)$ for many
pairs of subequations $(\bbe, \bbg)$.  These pairs are partially ordered by saying 
$$
(\bbe, \bbg) \prec (\bbe', \bbg') 
\quad  {\rm if}\quad
\bbe \ss\bbe'  \ \ {\rm and}\ \ \bbg' \ss\bbg.
\eqno{(1.10)}
  $$ 
   If this holds, then any
$\bbh_{\bbe, \bbg}$-harmonic on $X$ is automatically $\bbh_{\bbe', \bbg'}$-harmonic on $X$.

We shall see in Chapter 3 that  for any generalized equation $\bbh$, there exits a unique {\bf canonical
pair}  $(\bbe_{\rm min}, \bbg_{\rm max})$ defining $\bbh$, which is $\prec$ all other   pairs defining the set $\bbh$.  Thus an $\bbh_{\bbe_{\rm min}, \bbg_{\rm max}}$-harmonic on $X$, is harmonic 
for all pairs $(\bbe, \bbg)$ defining $\bbh$.  It will be called a {\bf canonical $\bbh$-harmonic} on $X$.
There remain some very interesting questions concerning this story (see  Chapter 5).
  
  One may wonder whether the closed sets $\bbh\ss\Symn$ which are generalized equations can be intrinsically characterized. They can be. They are exactly the sets satisfying:
  $$
  \bbh \ =\ (\overline{\bbh+\cp}) \cap ( \overline{\bbh-\cp}).
 \eqno{(1.11)}
  $$ 
  (See Theorem 3.5.)
  
  Given this, it is natural to consider the canonical $\bbh$-harmonics on an open set $X\ss\rn$.
  In Theorem 3.6 we show that   if $\bbh\ss\bbh'$ are generalized equations, 
  $$
  \text{
then any $\bbhc$-harmonic on $X$ is also
  a $\bbhc'$-harmonic on $X$.
  }
  $$

Now whether or not a closed subset $\bbh\ss\Symn$ satisfies (1.11), the set  
$$
\bbh^{\diamondsuit}\equiv  (\overline{\bbh+\cp}) \cap ( \overline{\bbh-\cp}) 
\eqno{(1.12)}
  $$ 
  is a generalized equation, and it is
  the smallest    generalized equation containing $\bbh$. That is, it is  contained in every other generalized equation containing $\bbh$.  (See Theorem 3.7 for all of these facts.)

Recall that if $\bbh= \partial \bbf \Rightarrow -\bbh = \partial \wt \bbf$.
More generally,  the negative of a generalized equation is also a generalized equation
$$
-\bbh_{\bbe, \bbg} \ =\    \bbh_{\wt\bbe,\wt \bbg}
\eqno{(1.13)}
  $$ 
This might be viewed as the dual of a generalized equation.

More importantly,  every generalized equation has a {\bf mirror}, 
which is obtained by interchanging $\bbe$ and $\bbg$.
A shortened form of our main result is the following.

\Theorem{2.6}  {\sl
Suppose $\bbh \equiv \bbh_{\bbe, \bbg} = \bbe\cap (-\wt \bbg)$ is a generalized equation with mirror $\bbh^*= \bbg \cap(-\wt\bbe)$,
and that $\O\ss \rn$ is a bounded domain.  Then
$$
\text{(A) \ \ Uniqueness for the (DP) for $\bbh$ holds on $\O \ \ \ \iff \ \ \Int \, \bbh=\emptyset$}
$$
Suppose that $\bo$ is smooth and strictly $\bbg$ and $\wt \bbg$-convex.  Then
$$
\text{(B) \ \ Existence  for the (DP) for $\bbh$ holds on $\O \ \ \ \iff \ \ \Int \, \bbh^*=\emptyset$}
$$
}

This allows us to divide generalized equations into four types depending on
whether the interiors of $\bbh$ and $\bbh^*$ are empty or not.  These types 
are strictly tied to the uniqueness and existence questions.

In proving Theorem 2.6 we established  some results of
independent interest.

\noindent
{\bf Proposition 2.16.}  {\sl Fix boundary values $\vf \in C(\bo)$.  If there exist solutions $h$ to the $\bbh$-(DP)
and $h^*$ to the $\bbh^*$-(DP) on $\O$, then $h=h^*$.
That is, $h=h^*$ is the common solution to the $\bbh$ and the $\bbh^*$ Dirichlet problems with boundary values
$\vf$.}

\noindent
{\bf Proposition 2.20.}  {\sl  Assume uniqueness holds for $\bbh$, i.e., $\Int \,  \bbh=\emptyset$, for the generalized equation
$\bbh = \bbe \cap (-\wt\bbg)$.  Given a domain as above and $\vf\in C(\bo)$,  then:
$$
\begin{aligned}
&\text{There exists $h\in C(\ob)$ with $h\bigr|_{\bo} = \vf$ and $h\bigr|_{\O}$ $\bbh$-harmonic}  \\
& \qquad\qquad\qquad \iff\ \qquad h_\bbe= h_\bbg, 
\quad \text{in which case $h=h_\bbe= h_\bbg$.}
\end{aligned}
$$
}

We would like to thank the referee for some very helpful suggestions.

%%%%%%%%%%%%%%%%%%%%%%%%%%%%%%%%%%%%%%%%%%%%%%%%%%%%%%
%%%%%%%%%%%%%%%%%%%%%%%%%%%%%%%%%%%%%%%%%%%%%%%%%%%%%%
%%%%%%%%%%%%%%%%%%%%%%%%%%%%%%%%%%%%%%%%%%%%%%%%%%%%%%
%%%%%%%%%%%%%%%%%%%%%%%%%%%%%%%%%%%%%%%%%%%%%%%%%%%%%%
%%%%%%%%%%%%%%%%%%%%%%%%%%%%%%%%%%%%%%%%%%%%%%%%%%%%%%
%%%%%%%%%%%%%%%%%%%%%%%%%%%%%%%%%%%%%%%%%%%%%%%%%%%%%%

%\vfill\eject
\vskip .3in

\noindent
{\headfont  2.  Main New Definitions and the Main Theorem.}

To begin we discuss more completely the notion  of a (determined) equation in the sense of [\Kry], [\DD] and [\DDR].

\Def {2.1. (Equation)}  A subset $\bbh\ss\Symn$ is a {\bf determined equation} or just 
 an {\bf equation}  if $\bbh = \partial \bbf$ for some subequation
$\bbf$.  In this case a {\bf solution to the equation} $\bbh$ (or an $\bbh$-{\bf harmonic}) on $X$,  is a function
$u$ such that $u \in \bbf(X)$  and $-u \in \wt\bbf(X)$.

Such functions are automatically continuous by definition.  For $C^2$-functions $u$ the consistency of
this definition  with the classical definition
that $D^2_x u\in \bbh = \partial\bbf$ for all $x$, follows from (1.3) above that 
$
\partial \bbf \ =\  \bbf\cap(- \wt \bbf),
$
and the consistency property for  subequations, mentioned in Section 1..

The main new concept in this paper is the following generalization.

\Def {2.2. (Generalized Equation)}  A subset $\bbh \equiv \bbh_{\bbe, \bbg}\ss\Symn$ is a {\bf generalized equation}
  if 
  $$
  \bbh =   \bbe\cap (-\wt \bbg) = \bbe\cap (\sim \Int \bbg)
  \eqno{(2.1)}
  $$ 
  for some pair of subequations $\bbe, \bbg$.  
  Just as in Definition 2.1 above, where $\bbe=\bbg =\bbf$, we define
  {\bf a solution to the (generalized) equation $\bbh_{\bbe, \bbg}$} (or an  $\bbh_{\bbe, \bbg}$-{\bf harmonic}) 
 to be a function $u$ with  
  $$
  u \in \bbe(X)  \and -u \in \wt\bbg(X),
 \eqno{(2.2)}
  $$ 
and let $\bbh(X)$ or $\bbh_{\bbe,\bbg}(X)$ 
  denote the space of $\bbh$-solutions on $X$.
  
As noted above, a generalized equation $\bbh$ can possibly be written as $\bbe\cap (-\wt \bbg)$ for many
pairs of subequations $(\bbe, \bbg)$.  These pairs are partially ordered by saying 
$(\bbe, \bbg) \prec (\bbe', \bbg')$ if $\bbe \ss\bbe'$ and $\bbg' \ss\bbg$.  If this holds, then any
$\bbh_{\bbe, \bbg}$-harmonic on $X$ is automatically $\bbh_{\bbe', \bbg'}$-harmonic on $X$.

  We present a guiding elementary example at this point to help the reader assimilate 
  the many definitions presented here.   Examples are discussed more fully in Section 4.
  
   \Ex{2.3.  (The Constrained Laplacian)} 
  Fix $r\geq 0$  and let 
   $$
   \bbh\equiv \{A : \tr A=0 \ \ {\rm and} \ \ -rI\leq A\leq rI\}.    
  \eqno{(2.3)}
   $$
   
%\vskip .4in
\centerline{   %\hskip 1.3in
           \includegraphics[width=.3\textwidth, angle=270,origin=c]{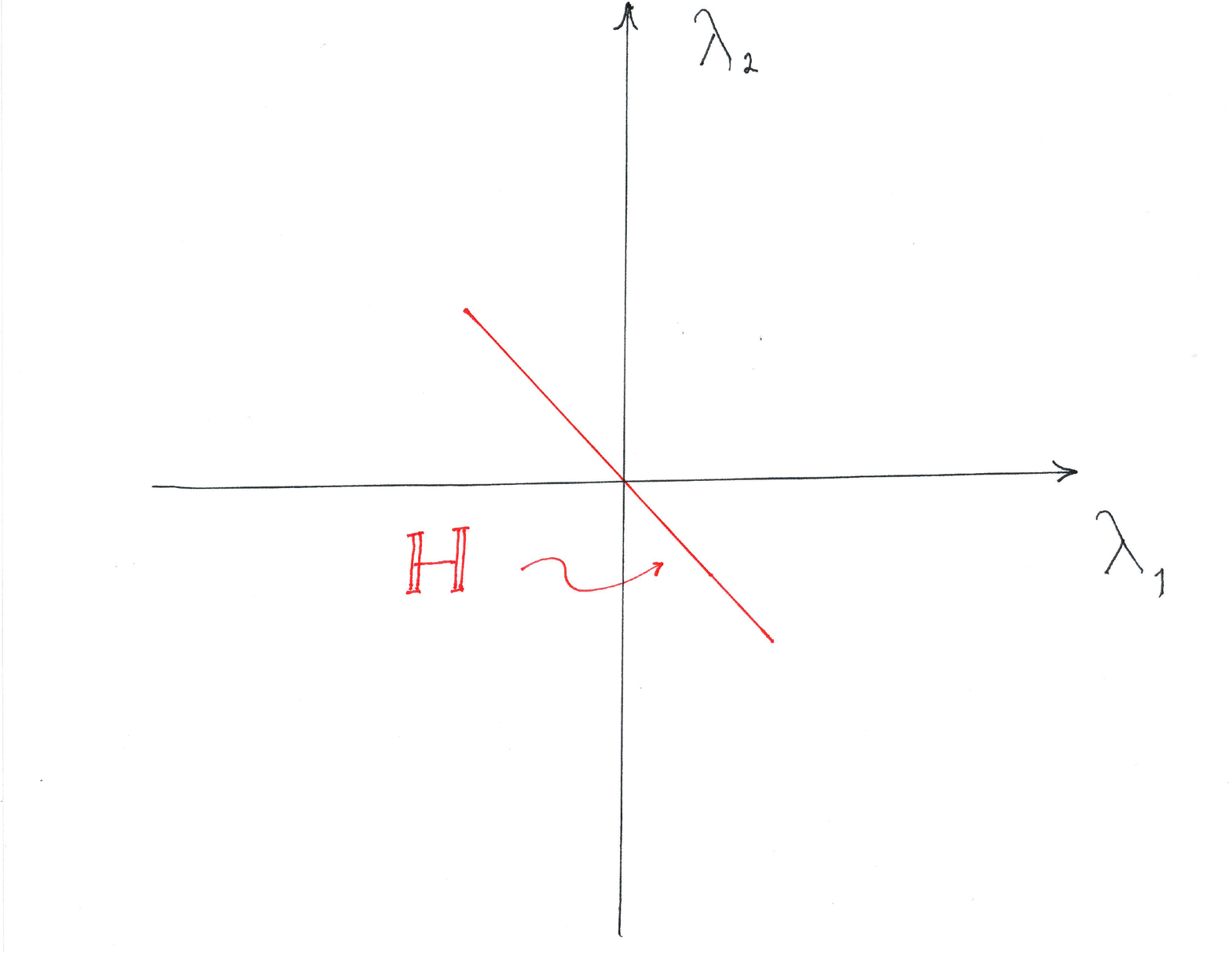}
          }

   \noindent
Define 
 $$
 \bbe_{\rm min} \ \equdef \  \bbh+\cp \ =\  \Delta\cap(\cp-rI)
 $$
 
\centerline{   %\hskip 1.3in
           \includegraphics[width=.3\textwidth, angle=270,origin=c]{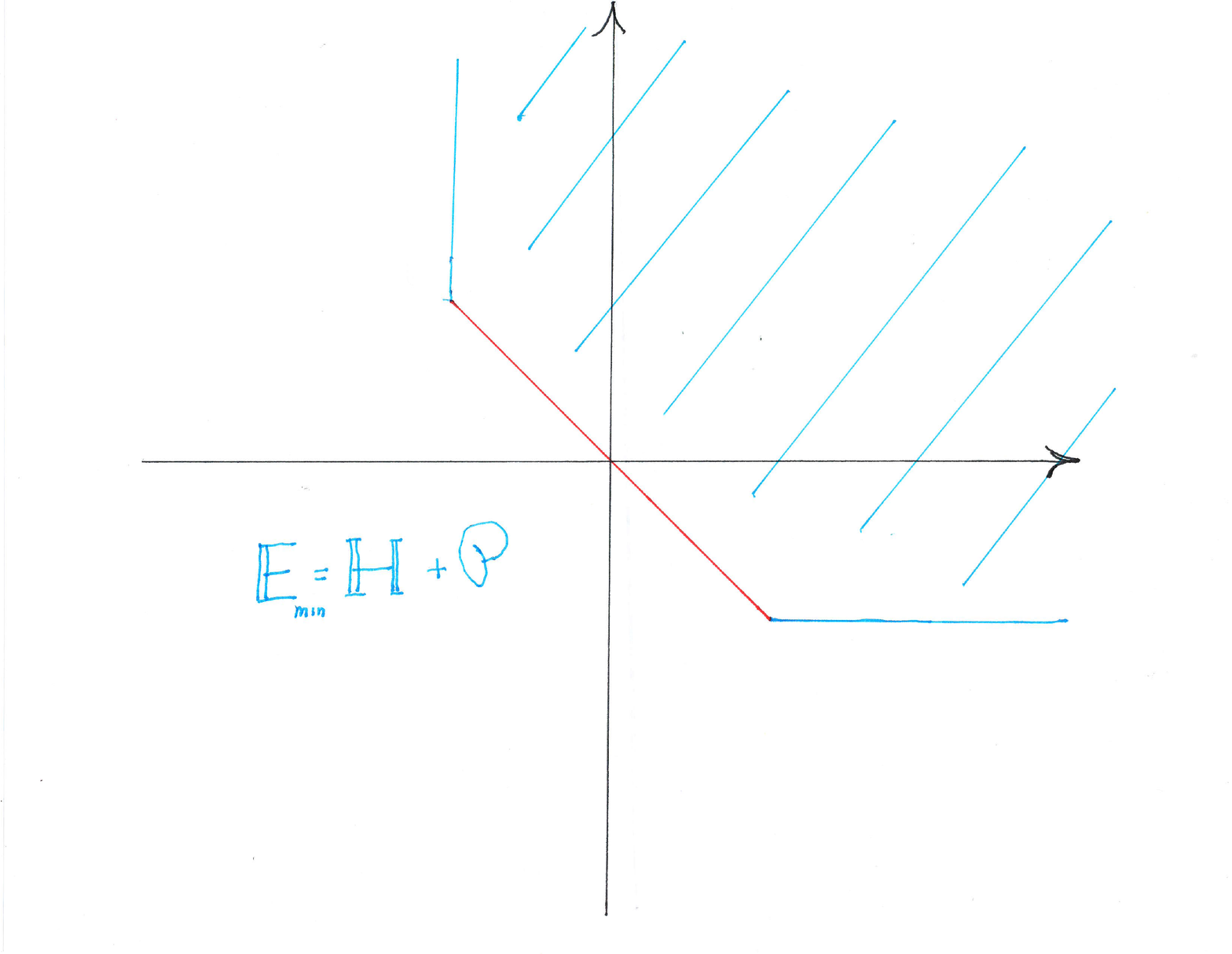}
          }

 \noindent
 and note that $-\bbem=\bbh-\cp$ since $-\bbh=\bbh$.  We then define $\bbgm \equiv \wt{\bbe}_{\rm min}$
 
\centerline{   %\hskip 1.3in
           \includegraphics[width=.3\textwidth, angle=270,origin=c]{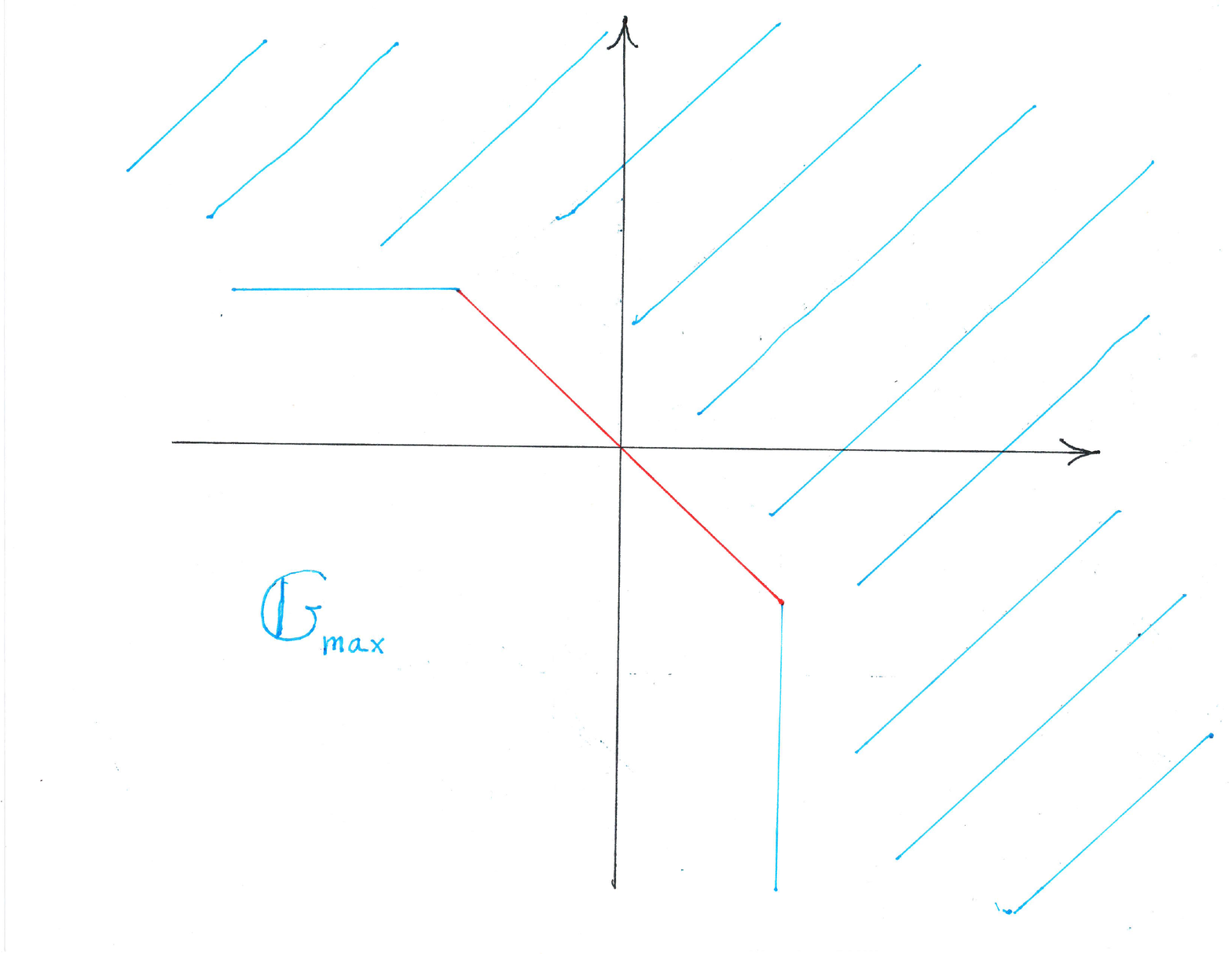}
          }

 \noindent
 Then we have
 $$
 \bbh\ =\ \bbem  \cap (- \wt{\bbg}_{\rm max})
 $$ 
 with $-\wt{\bbg}_{\rm max} = -\bbem = \bbh-\cp$
   
\centerline{   %\hskip 1.3in
           \includegraphics[width=.3\textwidth, angle=270,origin=c]{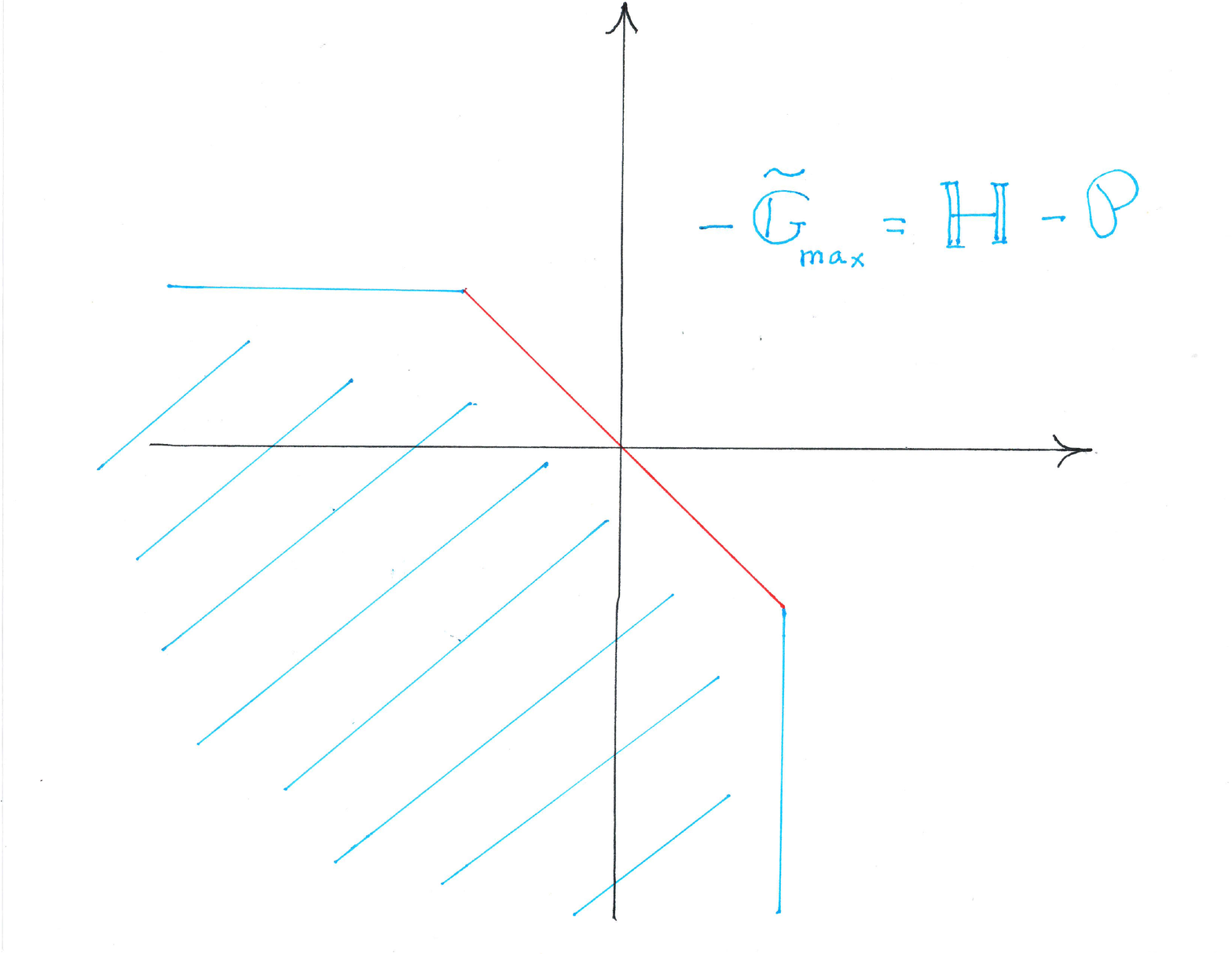}
          }

 We note that, with $\bbh$ defined by $\bbem, \bbgm$ as above,
 {\bf  $\bbh$-harmonic implies $\D$-harmonic}, which is of course obvious
for $C^2$-functions.  This follows from    Theorem 3.6 below.
 .

  \noindent
  {\bf Remark.\ (Intersections of GE's are GE's).} For   general intersections of subequations and their negatives
  $\bbh \equiv \bbe_1 \cap\cdots \cap \bbe_k \cap (-\wt\bbg_1)\cap\cdots\cap(-\wt\bbg_\ell)$,
we just get another generalized equation. The positivity condition (1.1) and the  closure  condition are both preserved under intersections. Hence $\bbe \equiv  \bbe_1 \cap\cdots \cap \bbe_k$ and
   $\wt\bbg \equiv \wt\bbg_1\cap\cdots\cap\wt\bbg_\ell$ are subequations, and $\bbh = \bbe\cap(-\wt\bbg)$.
     (Also since $\wt{\wt \bbf}= \bbf$ and $\wt{\bbf_1\cap \bbf_2} = \wt\bbf_1\cup\wt\bbf_2$ for
     subequations, one can show that $\bbg=\bbg_1\cup\cdots \cup \bbg_\ell$.)

  As before  by definition such functions are continuous with the coherence property  that if $u$ is $C^2$,  then
  $$
  \text{$u$ is $\bbh$-harmonic \qquad$\iff$\qquad $D^2_x u \in \bbh $ for all $x$}
\eqno{(2.4)}
$$

For any generalized equation $\bbh = \bbe\cap (-\wt \bbg)$ it is easy to see by (1.4) that the interior satisfies
$$
\Int \bbh\ =\ (\Int \bbe) \cap  (-\Int \wt \bbg)\ =\ (\Int \bbe) \cap  (\sim\bbg)
\eqno{(2.5)}
$$
   In particular,
$$\begin{aligned}
&\text{If $\bbh=\partial  \bbf = \bbf\cap (-\wt\bbf)$ is a determined equation,} \\
&\quad \text{then \ \  $\Int\,  \bbh = (\Int\, \bbf) \cap (\sim \bbf) = \emptyset.$}
\end{aligned}
\eqno{(2.6)}
$$

\vskip .3in

\centerline{\bf The Mirror}

\medskip

Each generalized equation  has a mirror, which we now define.

\Def {2.4. (The Mirror Equation)}  If  $$\bbh= \bbe \cap(-\wt\bbg)$$  is a generalized equation, its 
{\bf mirror} is defined  to be the generalized equation  $$\bbh^*= \bbg \cap(-\wt\bbe).$$

In Example 2.3 we have $\bbh^* = \bbgm \cap (-\wt{\bbe}_{\rm min}) = \bbgm\cap (-\bbgm)$
 since $\bbem$ and $\bbgm$ are dual to one another. 
 
 \centerline{   %\hskip 1.3in
           \includegraphics[width=.3\textwidth, angle=270,origin=c]{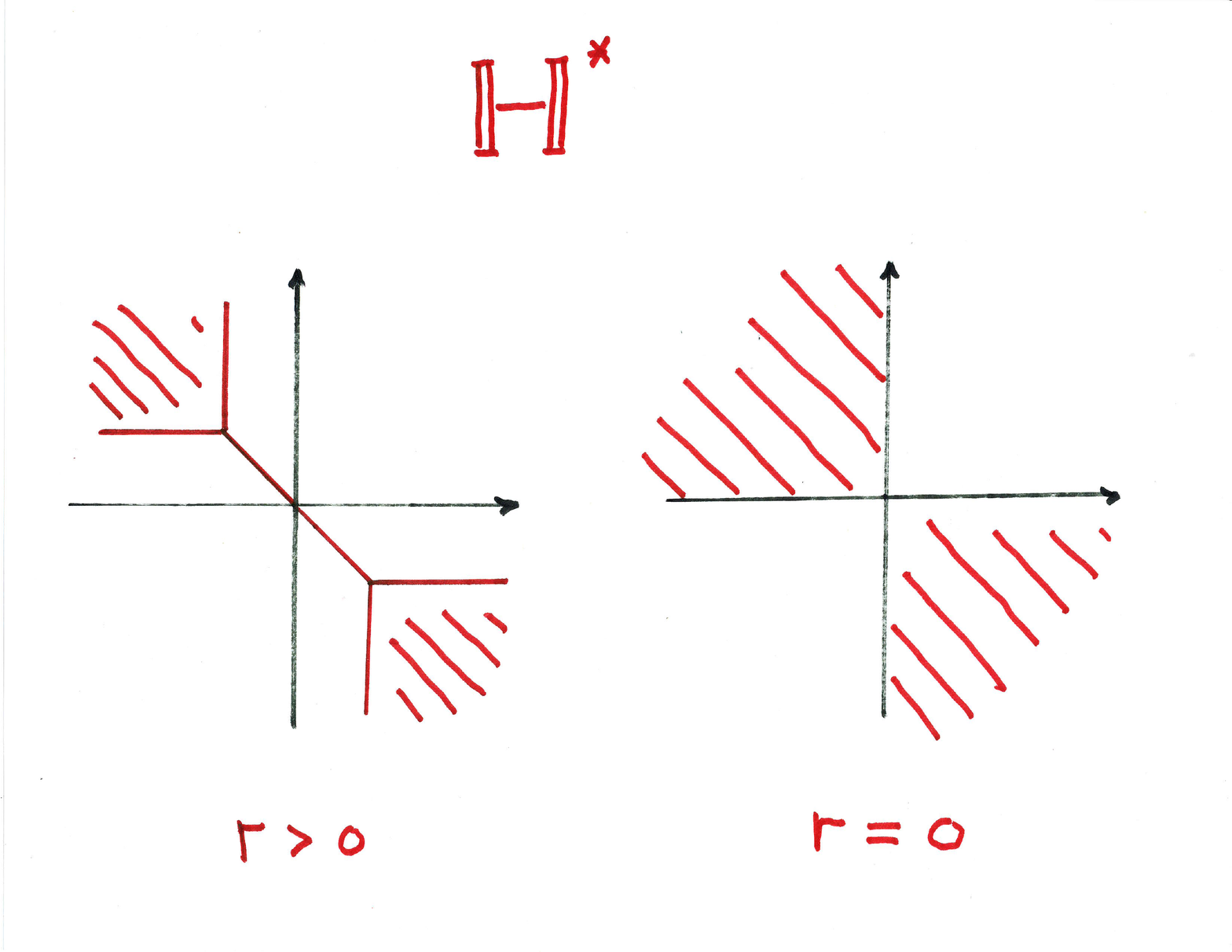}
          }
          
          This  $\bbh^*$ is a somewhat surprising 
 example of a generalized  equation from the following point of view. Consider $n=2$ and $r=0$. 
By using viscosity theory in applying Definition 2.2 to $\bbh^*$, we have   a notion
of when $D^2h$ has eigenvalues $\l_1$ and $\l_2$, which differ in sign, 
for  {\sl general continuous functions} $h$.  Moreover, this is consistent with the classical definition 
if $h\in C^2$.  On the other hand,  there is nothing ``elliptic'' about a mixed sign constraint on $D^2h$. 
(The reader could also consider $(\bbh^\diamondsuit)^*$ in Example 3.8 below where
the eigenvalues also satisfy $|\l_1|, |\l_2|\geq r$.)

\vskip.3in

\centerline{\bf Existence and Uniqueness}

%\medskip

Examination of  existence and uniqueness for the
 Dirichlet Problem for $\bbh$-harmonics leads to four distinct types of generalized equations,
as follows.

\Def {2.5}  Suppose $\O$ is a bounded domain in $\rn$. 
  We say that {\bf existence for the (DP) for} $\bbh$ holds on $\O$
if for all prescribed boundary functions $\vf \in C(\bo)$ there exists $h\in  C(\ob)$  satisfying

(a) \ \ $h\bigr|_{\O}$ is $\bbh$-harmonic, and

(b)\ \  $h\bigr|_{\bo}  =\vf$.

\noindent
 We say {\bf uniqueness for the (DP) for} $\bbh$ holds on $\O$
if for all $\vf \in C(\bo)$ there exists at most one $h\in  C(\ob)$
satisfying (a) and (b).

Now we can state our main result in this pure second-order constant coefficient case.

\Theorem{2.6}  {\sl
Suppose $\bbh \equiv \bbh_{\bbe, \bbg} = \bbe\cap (-\wt \bbg)$ is a generalized equation with mirror $\bbh^*= \bbg \cap(-\wt\bbe)$,
and that $\O\ss \rn$ is a bounded domain.  Then
$$
\text{(A) \ \ Uniqueness for the (DP) for $\bbh$ holds on $\O \ \ \ \iff \ \ \Int \, \bbh=\emptyset$}
$$
Suppose that $\bo$ is smooth and strictly $\bbg$ and $\wt \bbg$-convex.  Then
$$
\text{(B) \ \ Existence  for the (DP) for $\bbh$ holds on $\O \ \ \ \iff \ \ \Int \, \bbh^*=\emptyset$}
$$

\noindent
In fact, the following are equivalent:

(1) \ \  $\Int \,  \bbh = \emptyset$,

(2) \ \  Uniqueness for the (DP) for $\bbh$ holds on $\O$,

(3) \ \  $\bbe\ss \bbg$,

(4) \ \  $\bbh = \partial \bbe \cap \partial \bbg$, 

\noindent
and assuming that $\bo$ is smooth  and both strictly $\bbe$ and $\wt \bbe$ convex, these conditions 
are equivalent to

(5) \ \  Existence for the (DP) for $\bbh^*$ holds on $\O$.

\noindent
Interchanging $\bbe$ with $\bbg$  and $\bbh$ with $\bbh^*$, we have the mirror list of equivalences:

(1)$^*$ \ \  $\Int \,  \bbh^* = \emptyset$,

(2)$^*$ \ \  Uniqueness for the (DP) for $\bbh^*$ holds on $\O$,

(3)$^*$ \ \  $\bbg\ss \bbe$,

(4)$^*$ \ \  $\bbh^* = \partial \bbg \cap \partial \bbe$, 

\noindent
and assuming that $\bo$ is smooth  and both strictly $\bbg$ and $\wt \bbg$ convex,
these conditions  are equivalent to

(5)$^*$ \ \  Existence for the (DP) for $\bbh$ holds on $\O$.

}

\pf   It suffices to prove that (1) through (5) are equivalent because:
(A) is just the statement that  (2) $\iff$ (1), the mirror equivalences (1)$^*$ through (5)$^*$
are immediate from (1) through (5), and (B) is just the statement (1)$^*$ $\iff$  (5)$^*$.

Before proving the equivalence of (1) through (5) we list some trivial equivalences for any sets $\bbe$ and $\bbg$.

\Lemma{2.7}  {\sl
Suppose that $\bbe = \overline{\Int\, \bbe}$ and $\bbg = \overline{\Int\, \bbg}$, then
$$
(i) \ \ \bbe\ss \bbg  \quad
(ii)\ \ \Int\, \bbe \ss \bbg  \quad
(iii) \ \ \Int\, \bbe \ss \Int\, \bbg  \quad
(iv) \ \ \wt\bbg \ss \wt\bbe \quad 
(v) \ \ -\wt \bbg \ss -\wt\bbe
$$
are equivalent.   Interchanging $\bbe$ and $\bbg$ yields that 
$$
(i)^* \ \ \bbg\ss \bbe  \quad
(ii)^*\ \ \Int\, \bbg \ss \bbe  \quad
(iii)^* \ \ \Int\, \bbg \ss \Int\, \bbe  \quad
(iv)^* \ \ \wt\bbe \ss \wt\bbg \quad 
(v)^* \ \ -\wt \bbe \ss -\wt\bbg
$$
are equivalent.
}

\pf   Note that  $(i) \Rightarrow (ii)$ obviously,  $(ii) \Rightarrow (iii)$ since $\Int\, \bbe$ is an open subset contained in $\bbg$,
$(iii) \Rightarrow (i)$ by the hypothesis, $(i) \iff (iv)$ follows from the definitions of the duals, and 
$(iv) \iff (v)$  is trivial.\qed

\Cor{2.8}  {\sl
If $\bbh \equiv \bbe\cap (-\wt \bbg)$ is a generalized equation with mirror $\bbh^*= \bbg\cap (-\wt\bbe)$,
then 

(1) \ \ $\Int\, \bbh = \emptyset$ is equivalent to (i) through (v), and

(1)$^*$  \ \ $\Int\, \bbh^* = \emptyset$ is equivalent to (i)$^*$ through (v)$^*$.
}

\pf
Note that (2.5) says that 
$$
\Int\, \bbh \ =\ (\Int \,  \bbe) \sim \bbg.
%\eqno{(2.5)}
$$
Hence, $(1) \iff (ii)$.  Interchanging $\bbe$ and $\bbg$ yields $(1)^* \iff (ii)^*$.\qed

\noindent
{\bf Proof that (1) $\iff$ (3)}:  By  Corollary 2.8, (1) $\iff$ (i) which is (3). \qed

\noindent
{\bf Proof that (4) $\Rightarrow$ (1)}:  If $\bbh = \partial \bbe\cap\partial \bbg$, then in particular $\bbh\ss\partial \bbe$ which has
no interior.

\noindent
{\bf Proof that (1) $\Rightarrow$ (4)}: 
Note that  $\bbe = \partial \bbe\cup(\Int \,  \bbe)$ and $-\wt \bbg = \partial \bbg \cup(\sim \bbg)$ are
disjoint unions.  Hence, $\bbh \equiv \bbe\cap (- \wt\bbg)$ is the disjoint union of the four sets:
$\partial \bbe \cap  \partial \bbg$, $\partial \bbe \cap(\sim \bbg)$, 
$\Int \,   \bbe \cap  \partial \bbg$, and $\Int \,   \bbe \cap (\sim \bbg)$.
By (2.1) and (2.5), the last set $\Int \,   \bbe \cap (\sim \bbg) = \Int \,  \bbh = \emptyset$, so that
$\bbh$ is the disjoint union of the three remaining   sets.  However, 
$\Int \, \bbh = (\Int\, \bbe) \cap (\sim\bbg)  = \emptyset$ implies (3) and hence $\partial\bbe \ss\bbg$ or
$\partial\bbe \cap (\sim\bbg)  = \emptyset$.   By Lemma 2.7 (iii) $\Int\, \bbe\ss\Int\,\bbg$ so that 
$(\Int\, \bbe)\cap \partial \bbg = \emptyset$.
Thus three of these four sets are empty leaving 
$\bbh = \partial \bbe\cap \partial \bbg$.\qed

\noindent
{\bf Proof that (1) $\Rightarrow$ (2)}: 
Recall from [\DDR,  Def. 8.1] the following form
of comparison (C) for a subequation $\bbf$, which we will refer to as the {\bf zero maximum principle for sums}, and
abbreviate as either  (ZMP for sums) or (C).

\Def{2.9}  Given a relatively compact domain $\O$ we say that {\bf comparison  holds for $\bbf$ on $\O$}
if for all upper semi-continuous functions $u, v$ on $\ob$, with $u\bigr|_{\O}$ $\bbf$-subharmonic
and  $v\bigr|_{\O}$ $\wt\bbf$-subharmonic, one has
$$
u+v\ \leq\ 0\quad{\rm on}\ \ \bo
\qquad\Rightarrow\qquad
u+v\ \leq\ 0\quad{\rm on}\ \ \ob
\eqno{(ZMP\  for\  sums)}
$$

Comparison (C) always  holds for pure second-order subequations $\bbf \ss \Symn$
and domains $\O\ss\ss \rn$.  This was first established in [\DD, Rmk. 4.9 and Thm. 6.4]. 
(See [\CHLP] for many other constant coefficient situations where comparison always holds.
There are also some extensions in [\DDR]  to simply-connected, non-positively curved manifolds.)
More precisely we have:

\Theorem{2.10}  {\sl 
Suppose $\bbf\ss\Symn$ is a subequation and $\O\ss \rn$ is a bounded domain.  
Then comparison (C) holds for $\bbf$ on $\O$.
}

Now we can prove that (1) $\Rightarrow$ (2)

\Prop{2.11}  {\sl
Comparison (C) for both $\bbe$ and $\bbg$ on a domain $\O$ implies that:
$$
\Int\, \bbh  \ = \ \emptyset \qquad\Rightarrow \qquad  \text{uniqueness for the $\bbh$-(DP) on $\O$}
$$
}

\pf
By Corollary 2.8, $\Int \,  \bbh=\emptyset \Rightarrow \wt\bbg\ss\wt\bbe$.  
Therefore (C) for $\bbe$ implies  the (ZMP for sums)  if $u$ is $\bbe$-subharmonic and $v$ is $\wt\bbg$-subharmonic.
If $h_1, h_2$ are two solutions to the $\bbh$-(DP) on $\O$ with the same boundary values, 
then $u=h_1$ is $\bbe$--subharmonic and $v=-h_2$  is $\wt\bbg$--subharmonic on $\O$.
Since $u+v=0$ on $\bo$, the (ZMP) $\Rightarrow h_1\leq h_2$ on $\ob$. 
Interchanging $h_1$ and $h_2$ is possible since we are also assuming (C) for $\bbg$.  
This  proves $h_1=h_2$.\qed

\noindent
{\bf Proof that (2) $\Rightarrow$ (1)}: 
  
\Prop{2.12} {\sl
If there exists a function $h\in C^2(\O)\cap C(\ob)$ with 
$D^2_x h \in  \Int\, \bbh$ for all $x\in\O$, then uniqueness for the 
$\bbh$-(DP) on $\O$ fails.
}

\pf
Take $\vf = h\bigr|_{\bo}$.
For any function
 $\psi\in C^\infty_{\rm cpt}(\O)$, if $\e>0$ is sufficiently small, 
 we have $D^2_x (h+\e \psi)\in \bbh$ for all $x\in\O$.
 Thus the functions  $h+\e\psi$ give many solutions to the
  Dirichlet problem with the same boundary values $\vf$. \qed

The following trivial fact is peculiar to the pure second-order, constant coefficient case
(and the pure first-order case).

\Lemma{2.13}  {\sl
Given any non-empty subset $S\ss \Symn$,  there exists a function $h\in C^2(\rn)$ with $D^2_x h \in S$
for all $x\in\rn$.
}

\pf
Pick $A\in S$ and take $h(x) \equiv \half \bra {Ax}x$ so that $D^2_x h=A$ for all $x\in\rn$.\qed

Combining this Lemma with the previous Proposition proves the implication (2) $\Rightarrow$ (1) in the form:
$$
\begin{aligned}
&\Int\, \bbh \ \neq \ \emptyset
\qquad\Rightarrow \qquad   \\
\text{uniqueness for the }  &\text{$\bbh$-(DP) fails on all domains $\O\ss\rn$.}
\end{aligned}
\eqno{(2.7)}
$$
\qed

Next we treat the implication (1) $\Int \,  \bbh  =  \emptyset  \ \Rightarrow$\  (5) existence for the $\bbh^*$-(DP).

\Prop{2.14} {\sl
If existence for the $\partial \bbe$-(DP) holds on $\O$ (Definition 2.6), then 
$$
\Int \,  \bbh  =  \emptyset  \ \Rightarrow\ \text{existence for the $\bbh^*$-(DP) on $\O$.}
\eqno{(2.8)}
$$
}

\pf
By Corollary 2.8  
$\Int\, \bbh = \emptyset \ \Rightarrow \  \bbe\ss\bbg$.
Let $h$ denote the $\partial\bbe$-harmonic function solving the (DP) with boundary values $\vf$.
Since $h$ is $\bbe$-subharmonic and $\bbe\ss\bbg$, it is also $\bbg$-submarmonic.  
Since $-h$ is $\wt \bbe$-subharmonic, this proves that $h$ is $\bbh^*
=\bbg\cap(-\wt \bbe)$-harmonic.
\qed

 Recall the following from [\DD]. (As mentioned in \S 1, $\bbe(\O)$ denotes the space of $\bbe$-subharmonic functions on $\O$.)
 
 \Theorem{2.15. (Existence)}  {\sl Suppose $\O\ss\rn$ has a smooth boundary which is both $\bbe$ and $\wt \bbe$
strictly convex.  Given $\vf \in C(\bo)$, the Perron function 
$h(x) \equiv \sup \{ u\in \USC(\ob) : u\in \bbe(\O) \ \ {\rm and}\ \ u\bigr|_{\bo}\leq \vf\}$
 solves the $\partial \bbe$-(DP) on $\O$ for boundary values $\vf$.}

Combining Proposition 2.14 with Theorem 2.15 yields 
$$
(1)\ \ \Int\, \bbh  =  \emptyset  \ \Rightarrow\ (5)\ \  \text{existence for the $\bbh^*$-(DP)}
\eqno{(2.9)}
$$
on domains $\O$ with strictly $\bbe$ and $\wt\bbe$ convex smooth boundaries.

Before proving that  (5) $\Rightarrow$ (1), or that $\Int \,  \bbh  \neq  \emptyset$ 
implies non-existence for the $\bbh^*$-(DP), we need to establish
some preliminary facts, which are also of independent interest.

\Prop {2.16}  {\sl Fix boundary values $\vf \in C(\bo)$.  If there exist solutions $h$ to the $\bbh$-(DP)
and $h^*$ to the $\bbh^*$-(DP) on $\O$, then $h=h^*$.
That is, $h=h^*$ is the common solution to the $\bbh$ and the $\bbh^*$ Dirichlet problems with boundary values
$\vf$.}

\pf  By definition $h$ is $\bbe$-subharmonic and $-h$ is $\wt \bbg$-subharmonic on $\O$.
Also, $h^*$ is $\bbg$-subharmonic and $-h^*$ is $\wt \bbe$-subharmonic on $\O$.
Therefore,  
$$
h-h^* \ = \ 0 \ \ {\rm on}\ \ \bo \quad\Rightarrow \quad h-h^* \ \leq \ 0 \ \ {\rm on}\ \ \ob \ \ \  \text{by $\bbe$-comparison},
$$
$$
h^*-h \ = \ 0 \ \ {\rm on}\ \ \bo \quad\Rightarrow \quad h^*-h \ \leq \ 0 \ \ {\rm on}\ \ \ob \ \ \  \text{by $\bbg$-comparison}.
$$
Thus $h-h^*=0$ on $\ob$.\qed

\noindent
{\bf Note:}  Then  $h=h^*$ solves the generalized equation
$$
\bbh \cap \bbh^* \ =\ (\bbe\cap\bbg) \cap (-\wt{\bbe \cup \bbg}).
$$
(One can show that $\wt \bbe\cap \wt \bbg = \wt{\bbe\cup\bbg}$.
See [\DDR, Property (2) after Def.\ 3.1]   for arbitrary subsets of $J^2(X)$.)

\Prop{2.17} {\sl
Recall again that comparison holds for $\bbe$ and $\bbg$ on $\O$.  
From this we conclude the following.
If there exists a function $f\in C^2(\O)\cap C(\ob)$ with $D^2_x f \in \Int\, \bbh$
for all $x\in\O$, then there is no solution $h^*$ to the $\bbh^*$-(DP) on $\O$ 
with boundary values $\vf \equiv f\bigr|_{\bo}$. 
}

\pf
If $h^*$ exists, then since $f$ is an $\bbh$-solution, by  Proposition 2.16
we have  $h^*=f$, and hence $h^*$ is $C^2$. 
Thus $D^2f \in (\Int\, \bbh)\cap \bbh^* = (\Int \, \bbe\sim \bbg) \cap (\bbg \sim \Int\, \bbe) =\emptyset$. So this is impossible.\qed

\noindent
{\bf Proof that (5) $\Rightarrow$ (1) or that $\Int \, \bbh  \neq  \emptyset  \ \Rightarrow$\ Non-Existence for $\bbh^*$.}  
The fact that $\Int \, \bbh \neq \emptyset$ guarantees the existence of such a
function $f$ by Lemma 2.13,
 and hence the non-existence for the $\bbh^*$ Dirichlet problem.\qed

This completes the proof of Theorem 2.6.\qed

\vskip .3in
\centerline{\bf Four Types}

In light of Theorem 2.6, if one is given a generalized equation $\bbh=\bbe\cap(-\wt\bbg)$
with mirror $\bbh^*=\bbg\cap(-\wt\bbe)$, there are four distinct types possible which we label as follows.

\noindent
 {\bf Type I:} \ \ $\Int\, \bbh = \emptyset$ and  $\Int\, \bbh^* = \emptyset$

\noindent
{\bf Type II:} \ \  $\Int\, \bbh = \emptyset$ and  $\Int\, \bbh^* \neq \emptyset$

\noindent
{\bf Type III:} \ \  $\Int\, \bbh \neq \emptyset$ and  $\Int\, \bbh^* = \emptyset$

\noindent
{\bf Type IV:}  \ \  $\Int\, \bbh \neq \emptyset$ and  $\Int\, \bbh^* \neq \emptyset$

Note that Types I and II belong to part (A) of Theorem 2.6, and Types III and IV belong to part (B).
We shall now discuss each type.

\noindent
 {\bf Type I:}  $\Int\, \bbh =  \Int\, \bbh^* =\emptyset$.  This type is a ``determined equation''  $\partial \bbf$
 as defined in Definition 2.1, because by (1)  $\iff$ (3), and (1)$^*$  $\iff$ (3)$^*$,
 this is just the case where $\bbe$ and $\bbg$ are equal.  We will call this subequation 
 $\bbf$.  Thus $\bbh$ and $\bbh^*$ are $\bbf\cap(-\wt\bbf)= \partial \bbf$.  Theorems 2.10 and 2.15 apply directly.
 Comparison holds for all bounded domains, and existence holds if $\bo$ is smooth and strictly $\bbf$ and $\wt\bbf$
 convex using results from [\DD].

\noindent
{\bf Type II:} \ \  $\Int\, \bbh = \emptyset$ and  $\Int\, \bbh^* \neq \emptyset$. 
 Collecting together (1)--(5) and the negations of (1)$^*$--(5)$^*$ we have that \medskip
 
 \centerline{$\bbe$ is a proper subset of $\bbg$ and $\bbh=\partial \bbe\cap\partial \bbg \neq \bbh^*$}

\noindent
Uniqueness but not existence holds for $\bbh$ on any bounded domain $\O$.  The opposite
is true for $\bbh^*$, namely uniqueness for $\bbh^*$ fails on all bounded domains $\O$,
but if $\bo$ is smooth and both strictly $\bbe$ and $\wt\bbe$ convex, then existence holds for $\bbh^*$
on $\O$.  In addition $\bbh$ is a proper subset of both $\partial \bbe$ and $\partial \bbg$.
This is proven in (2.10) below.

For Type III we interchange $\bbe$ with $\bbg$ and $\bbh$ with $\bbh^*$.

\noindent
{\bf Type III:} \ \    $\Int\, \bbh \neq \emptyset$ and $\Int\, \bbh^* = \emptyset$. 
 Collecting together (1)$^*$--(5)$^*$ and the negations of  (1)--(5) we have that \medskip
 
 \centerline{$\bbg$ is a proper subset of $\bbe$ and $\bbh^*=\partial \bbg\cap\partial \bbe \neq \bbh$}

\noindent
Uniqueness but not existence holds for $\bbh^*$ on any bounded domain $\O$.  The opposite
is true for $\bbh$, namely uniqueness for $\bbh$ fails on all bounded domains $\O$,
but if $\bo$ is smooth and both strictly $\bbg$ and $\wt\bbg$ convex, then existence holds for $\bbh$
on $\O$. Also, $\bbh^*$ is a proper subset of both $\partial \bbe$ and $\partial \bbg$ by (2.10).

\noindent
{\bf Type IV:} \ \    $\Int\, \bbh \neq \emptyset$ and $\Int\, \bbh^* \neq \emptyset$. 
By (2.5) this is equivalent to 
$$
(\Int\, \bbe) \cap (\sim \bbg) \ \neq\ \ \emptyset
\qquad{\rm and} 
\qquad
(\Int\, \bbg) \cap (\sim \bbe) \ \neq\ \ \emptyset.
$$
Because of Lemma 2.7 (iii) and (iii)$^*$ this is equivalent to 
$$
\Int\, \bbe \ \not\ss \ \Int\, \bbg
\qquad{\rm and} 
\qquad
\Int\, \bbg \ \not\ss \ \Int\, \bbe.
$$
The main point here is that both existence and uniqueness for the (DP) for both $\bbh$ and $\bbh^*$ fail.

Next we consider the following.

\Prop{2.18.\ (Boundaries of Subequations)} {\sl If $\bbh$ is a determined  equation, then the subequation $\bbf$ with $\bbh = \partial \bbf$ is uniquely determined  by $\bbh$.  In fact, $\partial \bbe \ss\partial \bbg$ is enough to conclude that $\bbe=\bbg$ for any two subequations $\bbe$ and $\bbg$.  }

\pf
The first statement follows from (1.7) which holds for any subequation $\bbf$.

For the second statement note that one has  $\partial \bbe\ss\partial \bbg \Rightarrow \bbe = \partial\bbe+\cp \ss\partial \bbg +\cp = \bbg$.
However, $\partial \bbe\ss\partial \bbg  \iff -\partial \bbe\ss -\partial \bbg$, but $-\partial \bbe = \partial \wt \bbe$
and  $-\partial \bbg = \partial \wt \bbg$.  Hence, $\partial \wt \bbe\ss\partial \wt \bbg$, and this implies that 
$\wt \bbe\ss\wt\bbg$, which is equivalent to $\bbg\ss\bbe$.\qed

It follows that:
$$
\begin{aligned}
&\text{If $\bbh\equiv \bbe\cap(-\wt\bbg)$ is Type II, then} \\
&\text{$\bbh = \partial \bbe\cap \partial \bbg$ is a proper subset of both 
$\partial \bbe$ and $\partial \bbg$.}
\end{aligned}
\eqno{(2.10)}
$$

\pf  If $\bbh \equiv \partial \bbe = \partial \bbe\cap \partial \bbg$, then $\partial \bbe \ss\partial \bbg$,
so that by Proposition 2.18, $\bbe=\bbg$ and $\bbh$ is Type I.\qed

Next we begin to examine to what extent a generalized equation $\bbh=\bbe\cap(-\wt\bbg)$
determines the subequations $\bbe$ and $\bbg$. 
The answer in the determined case is given by next proposition.
However, open questions concerning more general cases can be found is Chapter 5.

\Prop{2.19.\ (Uniqueness of the Defining Pair $\bbe, \bbg$ in the Determined Case)}  {\sl 
Suppose that $\bbh=\partial \bbf$ where $\bbf$ is a subequation.
If $\bbh = \bbe \cap (-\wt\bbg)$ for subequations $\bbe$ and $\bbg$, then 
$\bbe=\bbg=\bbf$.}

\noindent
{\bf Proof.}
By (2.6) we have $\Int \bbh = \Int (\partial \bbf) =\emptyset$.
By (2.5) we have $\Int\bbh = (\Int \bbe) \cap (\sim \bbg)$.
Therefore, $\Int \bbe \ss\bbg$, and so $\bbe= \overline{\Int\bbe} \ss\bbg$.

Now  by (1.7) $\bbf = \partial\bbf +\cp =\bbh +\cp$, and by the hypothesis
 $\bbh=\bbe\cap(-\wt\bbg)$, we have $\bbh\ss\bbe$.
Thus, $\bbf = \bbh+\cp \ss\bbe+\cp=\bbe$.
By the same hypothesis $\bbh=\bbe\cap(-\wt\bbg)$, we have $-\bbh\ss\wt\bbg$.
Hence by (1.8), $\wt\bbf = -\bbh+\cp \ss \wt\bbg+\cp=\wt\bbg$,  so that $\bbg\ss\bbf$.  This proves that 
$\bbf\ss\bbe\ss\bbg\ss\bbf$.\qed

\noindent
\centerline{\bf The Question of Characterizing}
\centerline{\bf  the Boundary Functions for Existence}

Here we turn to a natural question which arises in the Non-Existence cases, Types II and IV.

\medskip

\centerline{
For which boundary functions $\vf\in C(\bo)$
}
\centerline{ does there exists a solution to the $\bbh$-Dirichlet problem?
}

First we discuss the Type II case: Non-Existence/Uniqueness for $\bbh$.
This is interesting, for example, for the constrained Laplacian above.

We make the assumption that $\O$ is a bounded domain with smooth boundary which
is both strictly $\bbe$ and $\wt\bbg$-convex.  Using the equivalent  versions (1)--(5) of the uniqueness
hypothesis for $\bbh$ in Theorem 2.6:
$$
\Int \, \bbh \ =\ \emptyset  \qquad\iff\qquad
\bbe\ \ss \ \bbg    \qquad\iff\qquad
\wt\bbg\ \ss\ \wt\bbe,
$$
this implies that $\bo$ is also strictly $\bbg$ and $\wt\bbe$-convex.
Let $h_\bbe \in C(\ob)$ denote the (unique) $\partial\bbe$-harmonic function on $\O$
with $h_\bbe\bigr|_{\bo} = \vf$,  and $h_\bbg \in C(\ob)$ denote the (unique) $\partial\bbg$-harmonic function on $\O$ with $h_\bbg\bigr|_{\bo} = \vf$.  One answer to the question is the following.

\noindent
{\bf Proposition 2.20.}  {\sl  Assume uniqueness holds for $\bbh$, i.e., $\Int \,  \bbh=\emptyset$, for the generalized equation
$\bbh = \bbe \cap (-\wt\bbg)$.  Given a domain as above and $\vf\in C(\bo)$,  then:
$$
\begin{aligned}
&\text{There exists $h\in C(\ob)$ with $h\bigr|_{\bo} = \vf$ and $h\bigr|_{\O}$ $\bbh$-harmonic}  \\
& \qquad\qquad\qquad \iff\ \qquad h_\bbe= h_\bbg, \quad \text{in which case $h=h_\bbe= h_\bbg$.}
\end{aligned}
$$
}

\noindent
{\bf Proof.}  Suppose that  $h_\bbe= h_\bbg$ and set $h=h_\bbe= h_\bbg$.  Then $h\bigr|_{\bo} = \vf$,
$h$ is  $\bbe$-subharmonic, and $-h$ is $\wt\bbg$-subharmonic on $\O$,  
which proves that $h$ is a solution to the $\bbh$ Dirichlet Problem (DP) on $\O$ with boundary values $\vf$.

Conversely, if these exists such an $h\in C(\ob)$, then $h$ also solves the $\partial\bbe$ (DP)
since $h$ is $\bbe$-subharmonic and $-h$ is $\wt\bbe \supset \wt\bbg$-subharmonic.   
Hence, $h_\bbe=h$.   
Similarly, $h_\bbg=h$ since $h$ is $\bbg\supset\bbe$-subharmonic and $-h$ is $\wt \bbg$-subharmonic. \qed

The   question  posed  above  can also be asked in the Type IV case.  
  It is intriguing for the ``$C^{1,1}$-equation'' in Example 4.2 below.
Here we take $r_1=r_2=\l >0$, and ask the question: How do we characterize the functions on the
boundary of a domain $\O$ which have a $C^{1,1}$-extension with Lipschitz constant $\l$ to all of $\O$?
This is related to $C^{1,1}$ Glaeser-Whitney extensions, for which there is a very large literature. 
The interested reader could consult [\GRU], [\DHGL] for results and some history.

We finish this section with a general comment.

\noindent
{\bf Remark 2.21. (Nice Properties of Generalized Equations).} Recall that 
generalized equations are preserved under:

(1)  Taking the mirror,

(2)  Taking Intersections (see Remark before  (2.4)),

(3)   Taking the negative (see (1.13) $\bbh_{\bbe, \bbg} = \bbh_{\wt\bbe, \wt\bbg}$).

\vfill\eject
%\vskip.2in
\noindent
{\headfont 3.   The Canonical Pair Defining a Given $\bbh$,
and an Intrinsic Characterization of Generalized Equations.}

In this section we look at the question of which closed subsets of $\Symn$
are generalized equations,  and we characterize them.  We start with the following.

\Lemma{3.1} {\sl
Suppose $\bbh\ss\Symn$ is any closed subset.  Then
$$
\bbe_{\rm min} \ \equdef\ \overline{\bbh +\cp}
\eqno{(3.1)}
$$
is a subequation containing $\bbh$, and it is minimal with respect
to these properties, i.e., if $\bbe$ is any other subequation containing $\bbh$,
then $\bbe_{\rm min} \ss\bbe$.
}

\noindent
{\bf Proof.}   The set $\bbh+\cp$ obviously satisfies positivity and so does its 
closure $\bbe_{\rm min}  = \overline{\bbh +\cp}$.  To see this, let $B\in \bbe_{\rm min}$
and choose  a sequence $B_k = A_k+P_k \to  B$ with $A_k \in\bbh$ and $P_k\geq 0$.
Now note that for all  $P\geq0$, $B_k+P    \to B+P$. \qed

\Cor{3.2} {\sl
The set 
$
-(\overline{\bbh-\cp}) \ =\ \overline{-\bbh+\cp}
$
is the minimal subequation containing $-\bbh$.  Let 
$$
\wt{\bbg}_{\rm max} \ \equdef \ \overline{-\bbh+\cp} \ \ \text{denote this subequation}.
\eqno{(3.2)}
$$
Then ${\bbg}_{\rm max}$ denotes its dual $\wt{\wt {\bbg}}_{\rm max}$.
 }

\noindent
{\bf Example.}  Let $\bbh=\{A\in\Sym(\bbr^2) : \l_1<0, \l_2>0 \ {\rm and}\ \l_1 \l_2 = -1\}$,
where $\l_1 \leq \l_2$ are the ordered eigenvalues of $A$.
Then $\bbh+\cp = \{A: \l_2>0\}$ is not closed.  
Nor is $\bbh-\cp = \{A:\l_1<0\}$ closed.  

At the moment we do not have an example
of this sort with $\bbh$  a generalized equation.

Although it is only in the determined case that $\bbh$ uniquely determines
the defining pair $\bbe, \bbg$ (namely,  $\bbe = \bbh+\cp = \bbg$),
we always have the following.

\Prop {3.3. (The Canonical Pair)}  {\sl
Suppose $\bbh \equiv \bbe_0\cap (-\wt\bbg_0)$ is any generalized equation. 
Then there exists a canonical choice for the 
subequation pair defining $\bbh$, namely $\bbe_{\rm min}$ and $\bbg_{\rm max}$:
$$
\bbhc \ \equiv\ \bbe_{\rm min}\cap(-\wt\bbg_{\rm max}) \ =\ (\overline{\bbh+\cp})  \cap  (\overline{\bbh-\cp}).
\eqno{(3.3)}
$$
This canonical  pair is characterized by the following property:  
 
If $\bbe, \bbg$ is any other   subequation pair 
yielding the same generalized equation $\bbh$, i.e., if  
{$\bbh=\bbe\cap(-\wt\bbg),$}
then
$$\bbe_{\rm min}\ss\bbe \qquad{\rm and} \qquad \bbg\ss\bbg_{\rm max}, $$
\noindent{ i.e., $\bbe_{\rm min}$ is minimal and $\bbg_{\rm max}$ is maximal.}

\noindent
In particular, if $h$ is $\bbhc= \bbe_{\rm min}\cap (-\wt \bbg_{\rm max})$ harmonic for the canonical
 min/max pair $\bbem, \bbg_{\rm max}$, then $h$ is   also
$\bbh = \bbe\cap (-\wt \bbg)$ harmonic for all other pairs $\bbe, \bbg$.
 }

\pf  
Since $\bbh\ss\bbem=\overline{\bbh+\cp}$ and $\bbh \ss -\bbgmt= \overline{\bbh-\cp}$, we have that
 $$
 \bbh \ \ss\ \bbem \cap (- \bbgmt) .
 $$

Now assume that $\bbh=\bbe\cap(-\wt\bbg)$ for a subequation pair $(\bbe, \bbg)$.  
Then $\bbh\ss \bbe$ and so 
$\bbem  \equiv \overline{\bbh+\cp} \ss \overline{\bbe+\cp} = \overline{\bbe} = \bbe$.
We also have $-\bbh \ss \wt\bbg$ which implies that 
$\bbgmt \equiv \overline{-\bbh+\cp} \ss\overline{\wt\bbg+\cp}=\overline{\wt \bbg} = \wt\bbg$.
Thus we have $\bbem \ss\bbe$ and $\bbg\ss\bbgm$.  
Therefore, $\bbem\cap (-\bbgmt) \ss \bbe\cap (-\wt\bbg)  = \bbh$.
With the display above, this implies that  $\bbh=  \bbem\cap (-\bbgmt)$.\qed

The next result concerns a generalized equation $\bbh$ belonging to part (A) of
Theorem 2.6 -- the {\bf uniqueness case}.  This is the case where $\Int\,\bbh=\emptyset$,
or equivalently by Theorem 2.6, this is the case of a generalized equation which is either of Type I or Type II.
Given a generalized equation $\bbh$, another characterization  is
$$
\Int\,\bbh=\emptyset \qquad\iff\qquad \bbh\ \ss\ \partial \bbf \ \ \text{for some subequation $\bbf$.}
$$
 (This follows easily from Theorem 2.6.)

\Prop {3.4}  {\sl
Suppose $\bbh$ is a generalized equation with $\Int \bbh=\emptyset$ (that is, $\bbh$ belongs to 
part (A) of Theorem 2.6 -- the uniqueness case).  Let $\bbe_{\rm min}, \bbg_{\rm max}$ 
denote the canonical min/max pair
with $\bbh = \bbem\cap(-\bbgmt)$.
Any subequation $\bbf$ with $\bbh\ss\partial \bbf$ must satisfy}
$$
\bbe_{\rm min} \ \ss\ \bbf\ \ss\  \bbg_{\rm max}.
$$

\noindent
{\bf Proof.}  Note that $\bbh\ss\partial \bbf \ \Rightarrow\ \bbe_{\rm min} \equiv \overline {\bbh+\cp}
 \ss\overline{\partial\bbf+\cp} =\overline\bbf=\bbf$.  
 Now $\bbh\ss\partial \bbf \ \iff\ -\bbh\ss\partial \wt\bbf = -\partial\bbf$.
Hence, $\bbh\ss\partial \bbf \ \Rightarrow\ -\bbh\ss\partial\wt \bbf
 \ \Rightarrow\ -\bbh +\cp \ss\partial\wt\bbf +\cp= \wt\bbf 
\ \Rightarrow\  \wt \bbg_{\rm max} \equiv  \overline{ -\bbh +\cp} \ss \wt\bbf \ \iff\ \bbf \ss  \bbg_{\rm max}$.\qed

One may wonder whether the closed sets $\bbh\ss\Symn$, which are generalized equations, can
be intrinsically characterized.  They can be.

\Theorem {3.5. (The Characterization of Generalized Equations)}  {\sl
A closed subset $\bbh\ss\Symn$ is a generalized equation if and only if }
$$
\bbh \ =\ ( \overline{\bbh+\cp}) \cap  (\overline{\bbh-\cp}).
\eqno{(3.4)}
$$

\noindent
{\bf Proof.}
First, by (3.3) this equation holds for any generalized equation $\bbh$.

For the converse,  recall from Lemma 3.1 that the sets 
$\bbe \equiv \overline{\bbh+\cp}$ and $\wt\bbg \equiv \overline{(-\bbh+\cp)}$ are always 
subequations, so that $\bbe\cap (-\wt \bbg)$ is a generalized equation.
\qed

As we have seen, if $\bbh$ is a generalized equation then the set of pairs $(\bbe, \bbg)$
with $\bbh=\bbe\cap(-\wt\bbg)$is partially ordered with a unique  minimal element
$(\bbem, \bbgm)$.  So associated to $\bbh$ we have this canonical subequation pair $(\bbem, \bbgm)$,
and it is natural to consider the associated {\bf canonical $\bbhc$-harmonics.}

\Theorem{3.6}  {\sl
Let $\bbh\ss\bbh'$ be two generalized equations as in Theorem 3.5.  Then any function $h$ which is 
$\bbhc$-harmonic on an open set $\O\ss\rn$ is also $\bbhc'$-harmonic on $\O$.
}

\noindent
{\bf Proof.}
We shall use the second form of (3.3).
Note that $\overline{\bbh+\cp} \ss \overline{\bbh'+\cp}$ and 
$\overline{\bbh-\cp} \ss \overline{\bbh'-\cp}$.  Therefore,
$\bbem\ss\bbem'$, and $- \wt \bbgm \ss -\wt{\bbgm'}$ which means that 
$\wt \bbgm \ss \wt{\bbgm'}$.  Thus if $h$ is $\bbem$-subharmonic and $-h$ is 
$\wt\bbgm$-subharmonic, this also holds for the primed subequations.\qed

This theorem says that the partial ordering by inclusion on the family of closed subsets
$\bbh\ss\Symn$ which are subequations, carries over to their canonical harmonics  on any open
$\O\ss\rn$. The reader will also recall that a canonical $\bbhc$-harmonic on $\O$ is
also an $\bbh_{\bbe, \bbg}$-harmonic for any other pair $(\bbe, \bbg)$ with $\bbh=\bbe\cap(-\wt\bbg)$.
One could certainly wonder whether every $\bbh_{\bbe, \bbg}$-harmonic is canonically $\bbhc$-harmonic.
(This is discussed as the   ``Broadened Equation Question'' in Section 5.)  If this were true, the theory
would be very tidy.

\Theorem {3.7. (The  GE Associated to a Closed Set)}  {\sl
Let $\bbh\ss\Symn$ be any closed subset.  Then the pair of subequations
$$
(\bbe^\diamondsuit, \bbg^\diamondsuit) \ \equdef \  (\overline{\bbh+\cp}, (\overline{-\bbh +\cp})^\sim)
\eqno{(3.5)}
$$
  gives   a generalized equation 
  $$
  \bbh^\diamondsuit  \ \equdef \   \bbe^\diamondsuit  \cap (-\wt\bbg^\diamondsuit ) = 
(\overline{\bbh+\cp})\cap (\overline{\bbh-\cp})
\eqno{(3.6)}
$$
containing $\bbh$, and it is the smallest such.  That is, if
$(\bbe, \bbg)$ is a subequation pair with $\bbh\ss \bbe\cap(-\wt \bbg)$, then
$\bbh^\diamondsuit   \ss \bbe\cap(-\wt\bbg)$.
Moreover, $(\bbe^\diamondsuit, \bbg^\diamondsuit)$ is the canonical min/max pair 
defining $\bbh^\diamondsuit$.

More succinctly, every closed subset  $\bbh\ss \Symn$ gives rise to a {\bf minimal generalized equation}
containing $\bbh$, namely}
$$
\bbh^\diamondsuit \ \equiv\  (\overline{\bbh+\cp})\cap (\overline{\bbh-\cp}).
$$

\noindent
{\bf Proof.}
From Lemma 3.1 and Corollary 3.2 we know that $(\overline{\bbh+\cp}, (\overline{-\bbh +\cp})^\sim)$
 is a pair of subequations defining the generalized equation $ \bbh^\diamondsuit$ given in (3.6),
which clearly contains $\bbh$.  

Suppose $\bbe, \bbg$ are subequations with $\bbh\ss \bbe\cap( -\wt\bbg)$. Then $\bbh\ss\bbe$ and so $\overline{\bbh +\cp} \ss \overline{ \bbe+\cp}=\overline{\bbe}=\bbe$. 
Also $\bbh \ss -\wt\bbg$, so $-\bbh \ss \wt\bbg$ and therefore $\overline{-\bbh+\cp}\ss
\overline{\wt \bbg +\cp} \ss \wt \bbg +\cp \ss\wt\bbg$.
Hence, $\bbh^\diamondsuit = (\overline{\bbh+\cp})\cap (\overline{\bbh-\cp}) \ss \bbe\cap(-\wt\bbg)$ as claimed.
Finally, $\bbem \equdef (\overline{\bbhd+\cp}) \supset (\overline{\bbh+\cp}) \equdef \bbed$
(since $\bbhd\supset \bbh$), and 
$\bbhd \ss \overline{\bbh+\cp} \equdef \bbed$  implies $\overline{\bbhd+\cp} \ss \bbed$,
which proves  $\bbem= \bbed$.
The proof that $-\wt{\bbg}^\diamondsuit = \overline{\bbh-\cp}$, so that $\bbgm=\bbgd$, is similar.
\qed

%On the other hand, if the equation holds, then $\bbh$ is a generalized equation,
%because $ {\bbh+\cp}$ is a subequation, and ${\bbh-\cp}= -({-\bbh+\cp})$
%is $-\wt\bbg$ where $\bbg=({-\bbh+\cp})^{\sim}$ is a subequation.\qed

\noindent
{\bf Example 3.8.}    Let 
$$
\bbh \equiv \{t\cdot  {\Id} : -1\leq t\leq 1\}.
$$
\centerline{   %\hskip 1.3in
           \includegraphics[width=.3\textwidth, angle=270,origin=c]{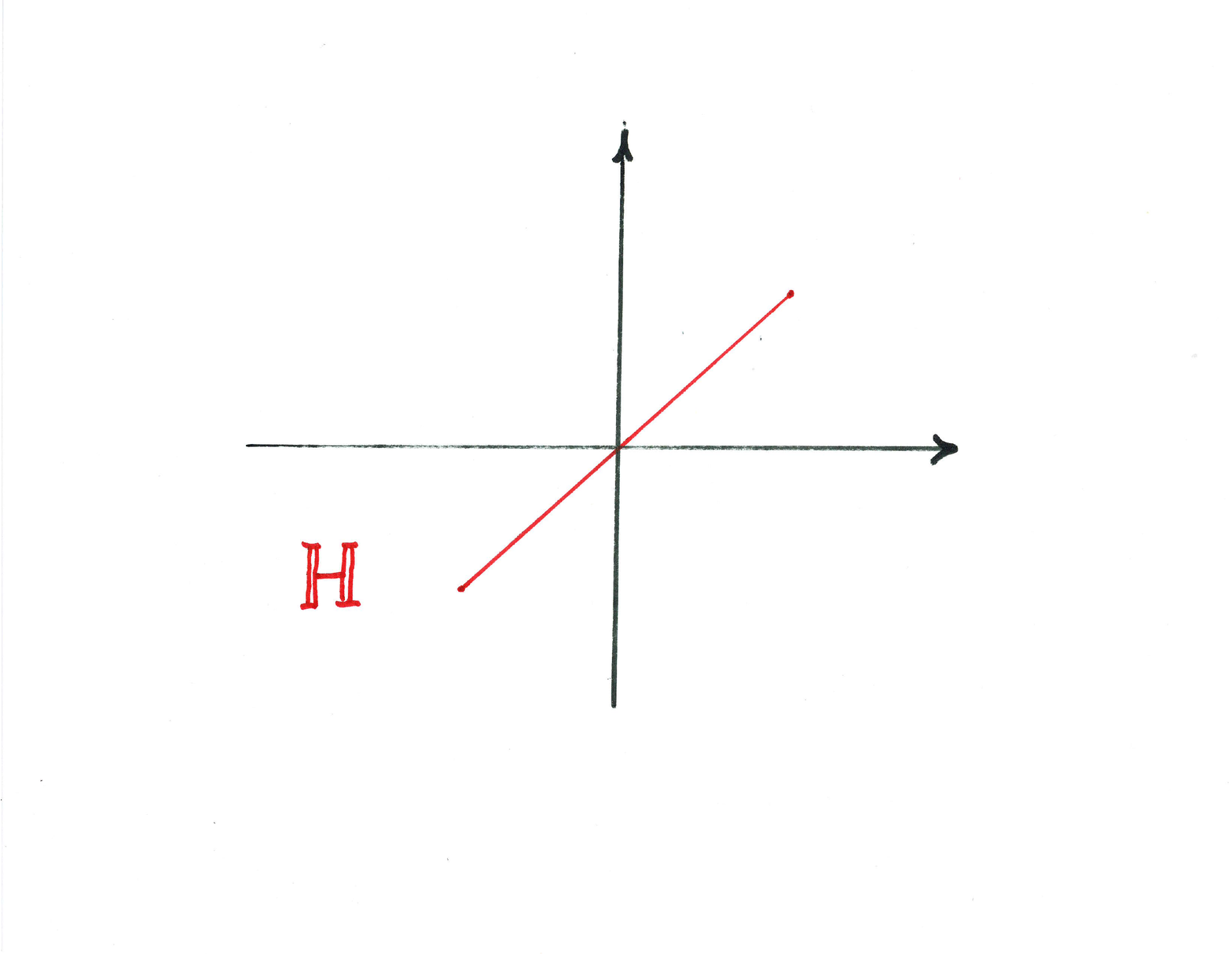}
          }
          \noindent Then
$$
\bbem \ =\   \bbh+\cp \ =\ \{A:A\geq -\Id\} \and -\wt{\bbg}_{\rm max}   \ =\   \bbh-\cp    \ =\ \{A:A \leq  \Id\}.
$$
 {   %\hskip 1.3in
           \includegraphics[width=.3\textwidth, angle=270,origin=c]{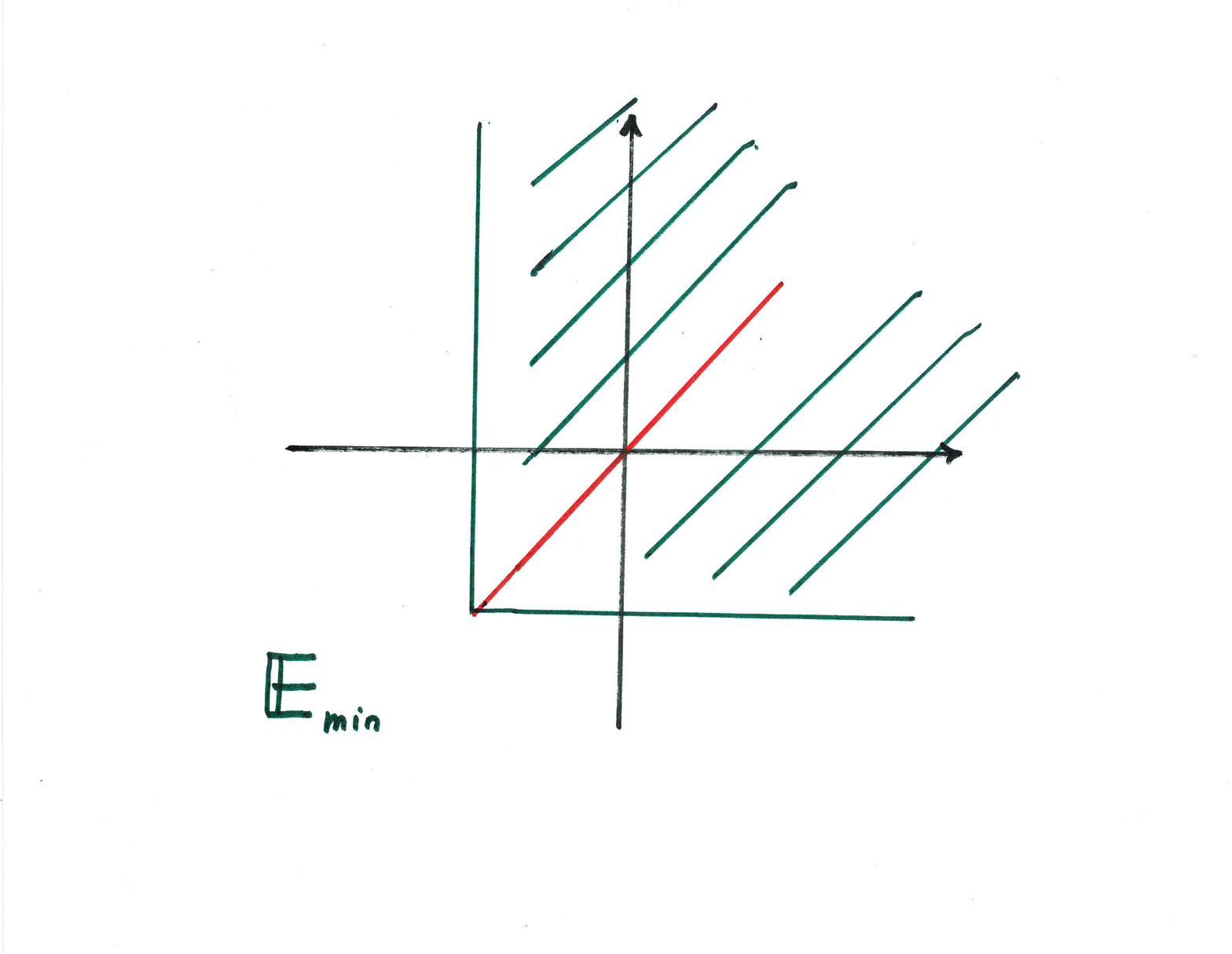}
          }  \hskip1.2in
          {   %\hskip 1.3in
           \includegraphics[width=.3\textwidth, angle=270,origin=c]{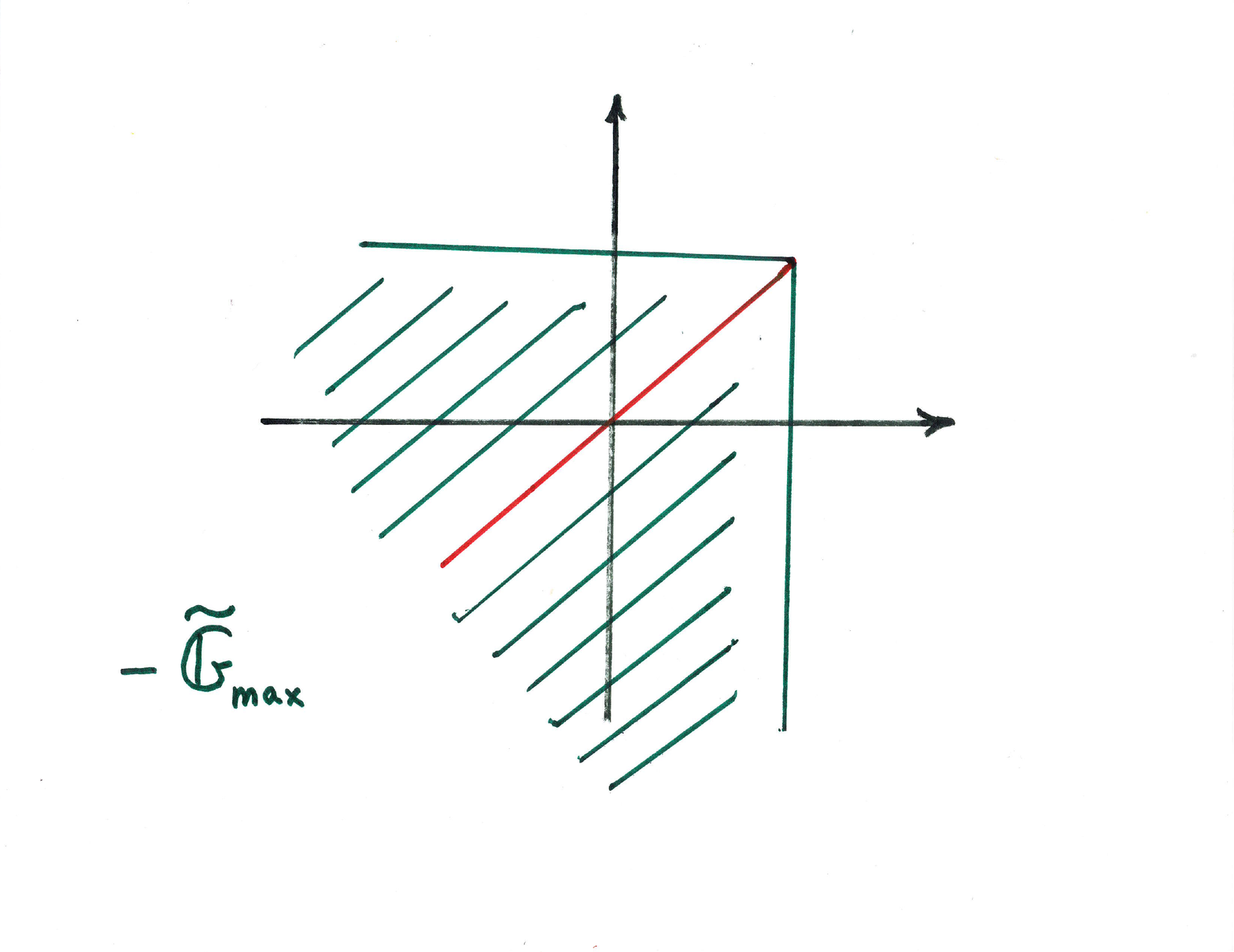}
          }

%\vfill\eject

   \noindent Therefore,  the minimal generalized equation containing $\bbh$ is 
$$
\bbh^\diamondsuit \ =\ (\bbh+\cp) \cap (\bbh-\cp) \ =\ \{A : -\Id \leq A\leq \Id\}.
$$
\centerline{   %\hskip 1.3in
           \includegraphics[width=.3\textwidth, angle=270,origin=c]{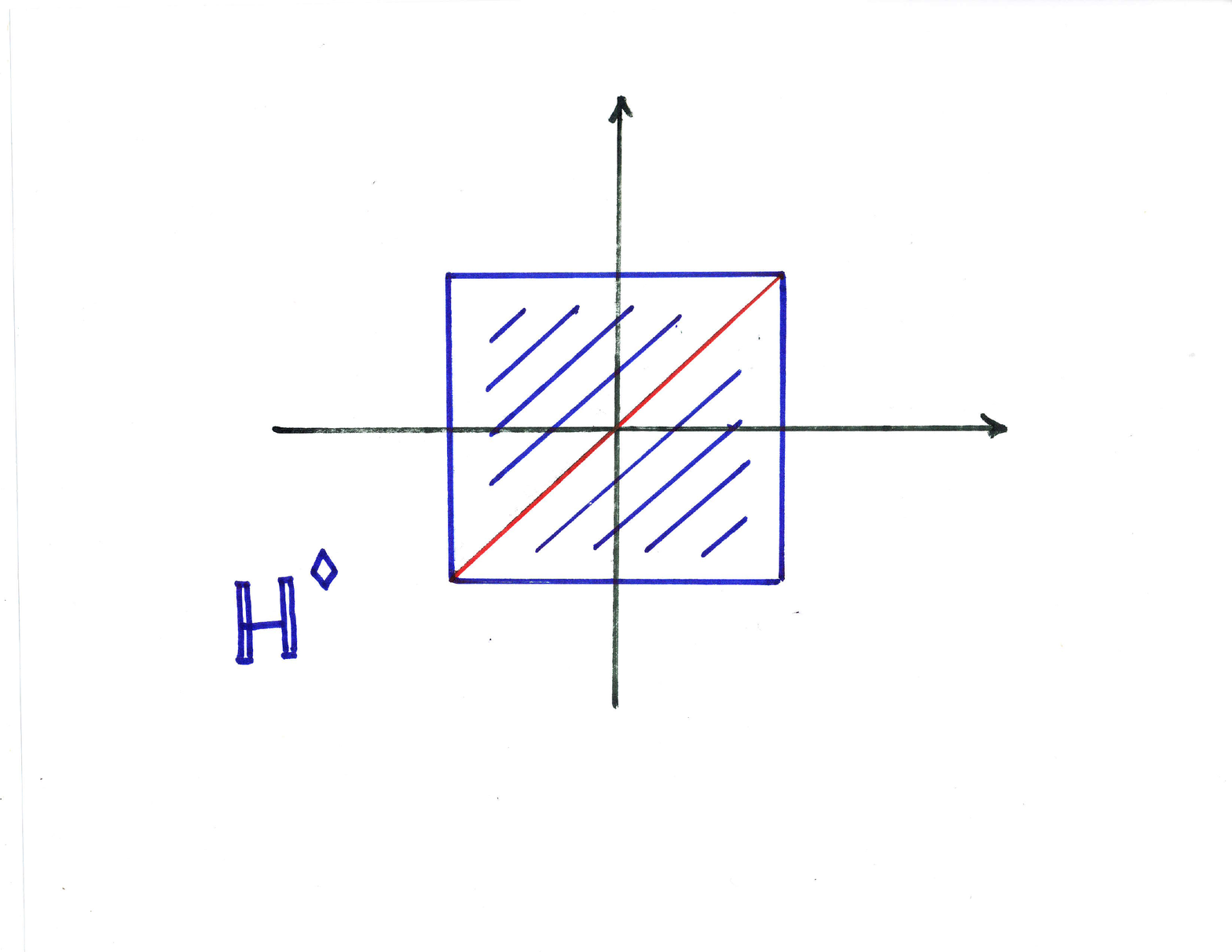}
          }

   \noindent Note that  the canonical   pair $\bbem,\bbgm$ for this generalized equation
   $\bbhd$ is:
$$
\bbem  \ =\ \{A: A\geq -\Id\} \and \ \ {\bbg}_{\rm max}   \ =\ \{A:  A - \Id \in\cpt\}.
$$
 {   %\hskip 1.3in
           \includegraphics[width=.3\textwidth, angle=270,origin=c]{Figure11.pdf}
          }  \hskip.9in
          {   %\hskip 1.3in
           \includegraphics[width=.3\textwidth, angle=270,origin=c]{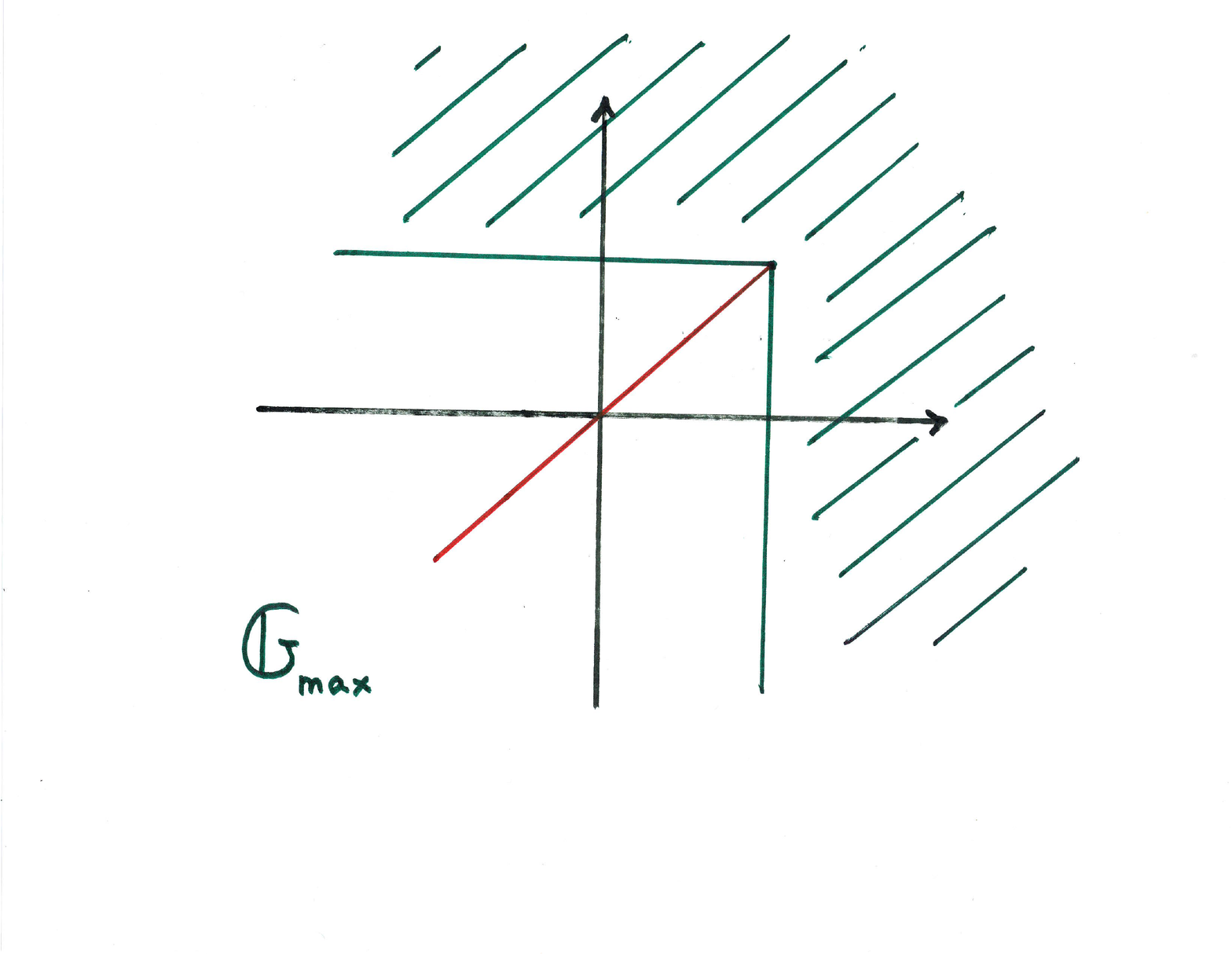}
          }
\noindent
{\bf Some Facts:}  Here $\bbgm$ and $\bbem$ are dual to one another, and so
$$
\bbh^\diamondsuit \ \equiv \ \bbem \cap (-\wt{\bbg}_{\rm max}) \ =\ \bbem \cap (-{\bbem})
$$
is a Class II generalized equation (see Section 4).

The {\bf mirror} of $\bbh^\diamondsuit$  is 
$$
(\bbh^\diamondsuit)^* \ =\ \bbgm \cap (-\bbgm) \ =\ \{A : A-\Id \in\cpt \ {\rm and}\ -A-\Id \in \cpt\}.
$$
\centerline{   %\hskip 1.3in
           \includegraphics[width=.3\textwidth, angle=270,origin=c]{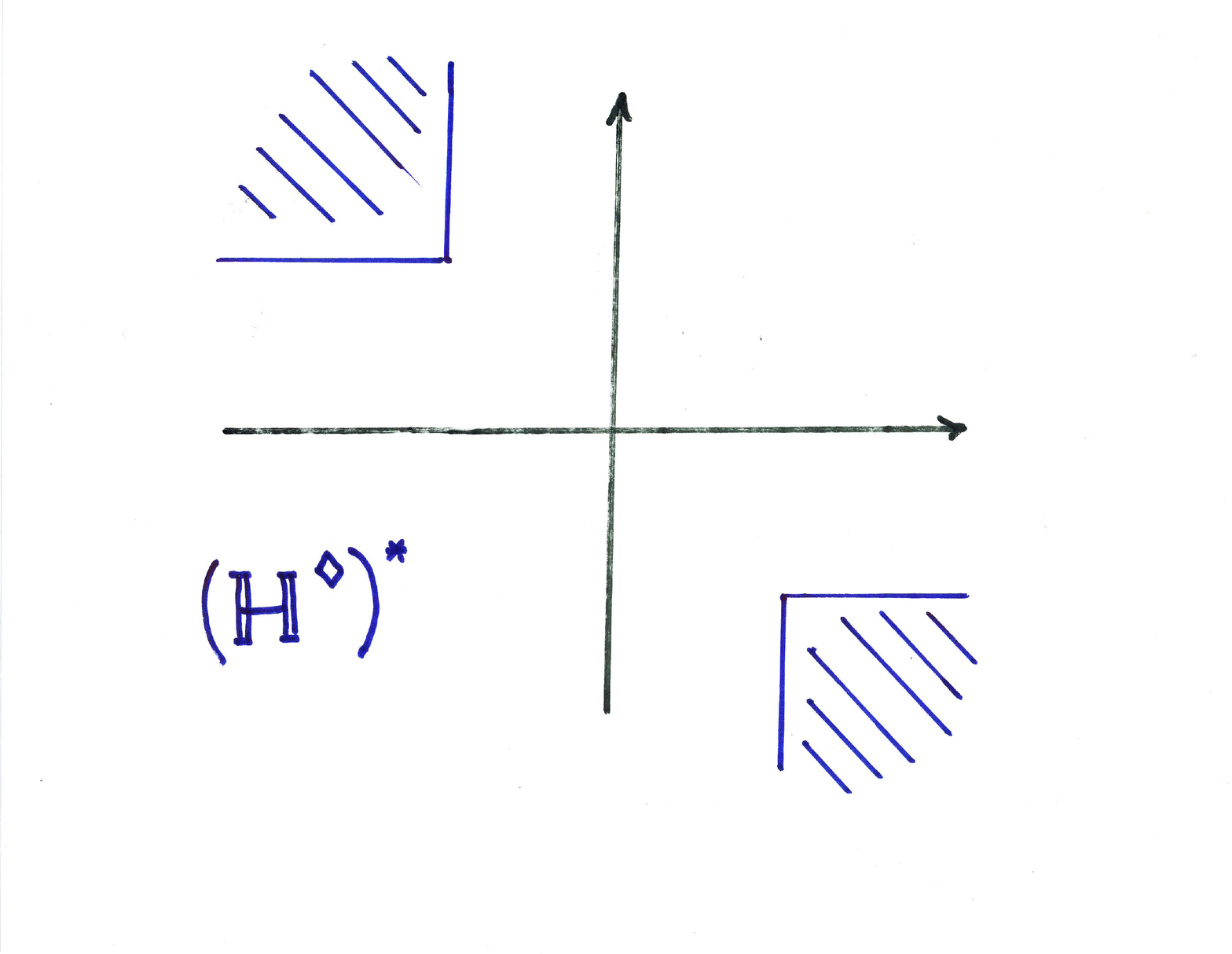}
          }
Both $\Int \bbh^\diamondsuit \neq \emptyset$ and $\Int (\bbh^\diamondsuit)^* \neq \emptyset$,
so that $\bbh^\diamondsuit$ is Type IV.

\vskip.3in

%\vfill\eject

\noindent
{\bf Example 3.9.}    Let 
$$
\bbh \equiv \{\tr(A) =0\}\cup  \{t\cdot  {\Id} : -1\leq t\leq 1\}.
$$
\centerline{   %\hskip 1.3in
           \includegraphics[width=.3\textwidth, angle=270,origin=c]{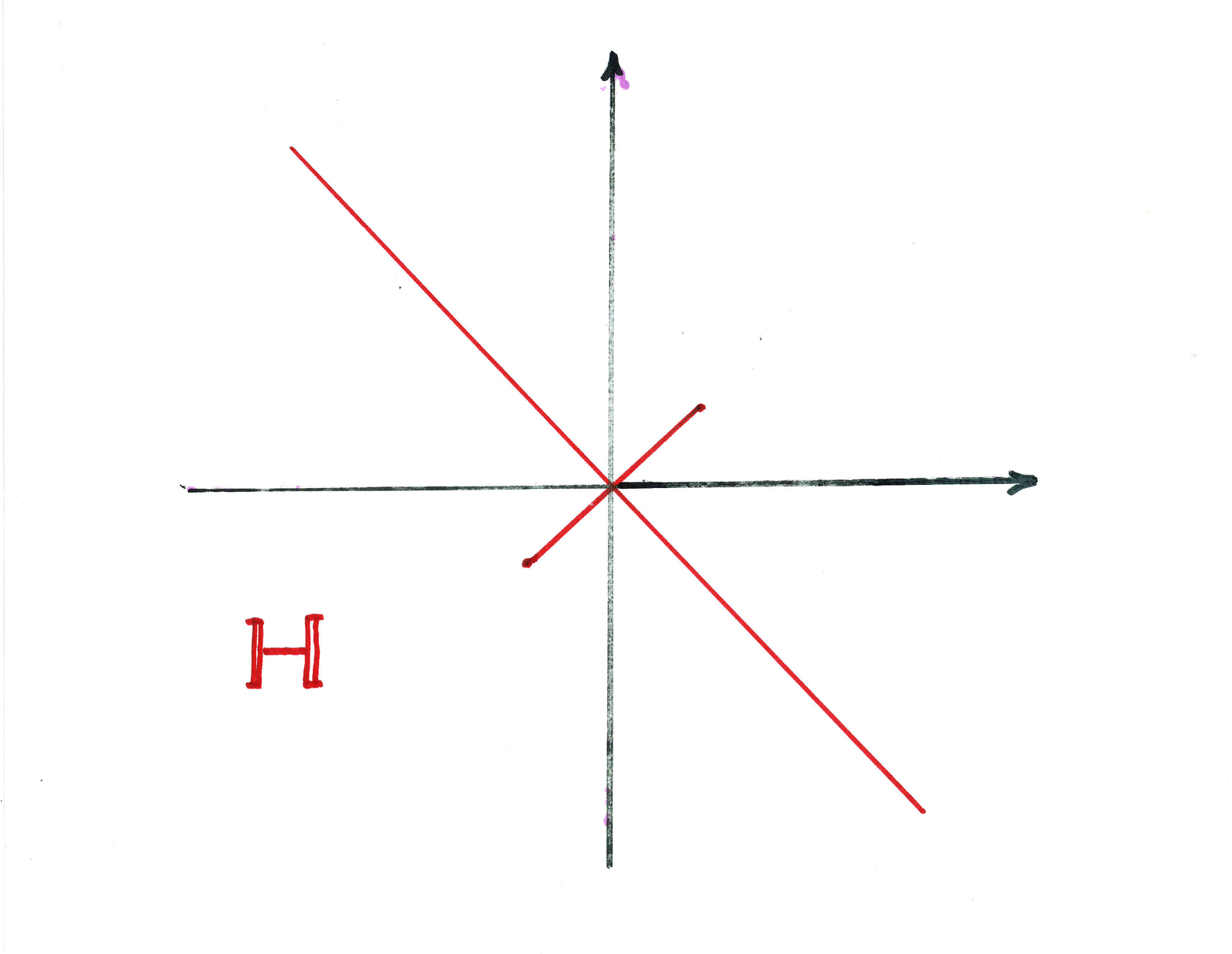}
          }
          \noindent Then we have 

$$
\bbem \ =\   \bbh+\cp \qquad \and \qquad -\wt{\bbg}_{\rm max}   \ =\ \bbh - \cp 
$$
 {   %\hskip 1.3in
           \includegraphics[width=.3\textwidth, angle=270,origin=c]{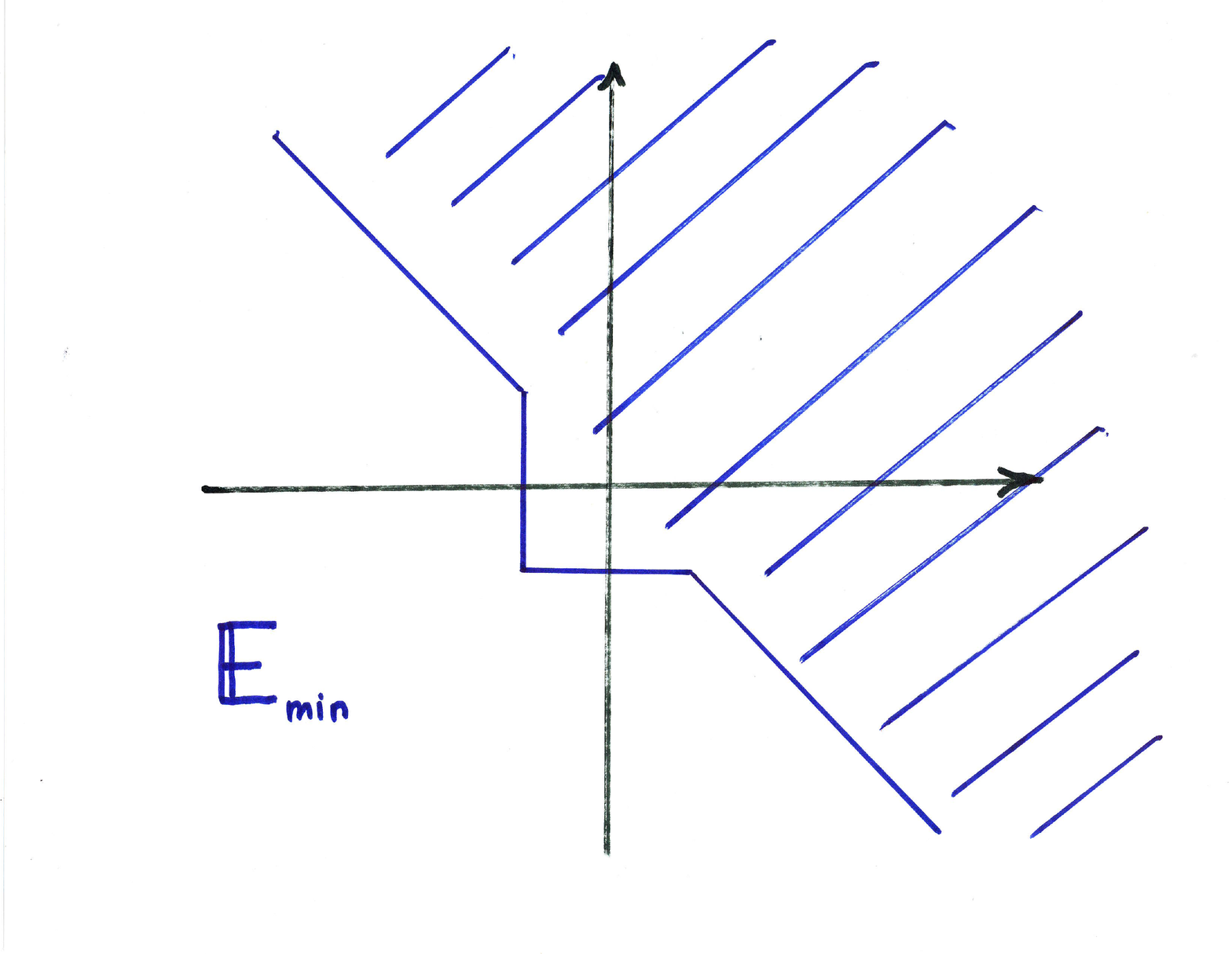}
          }  \hskip1.1in
          {   %\hskip 1.3in
           \includegraphics[width=.3\textwidth, angle=270,origin=c]{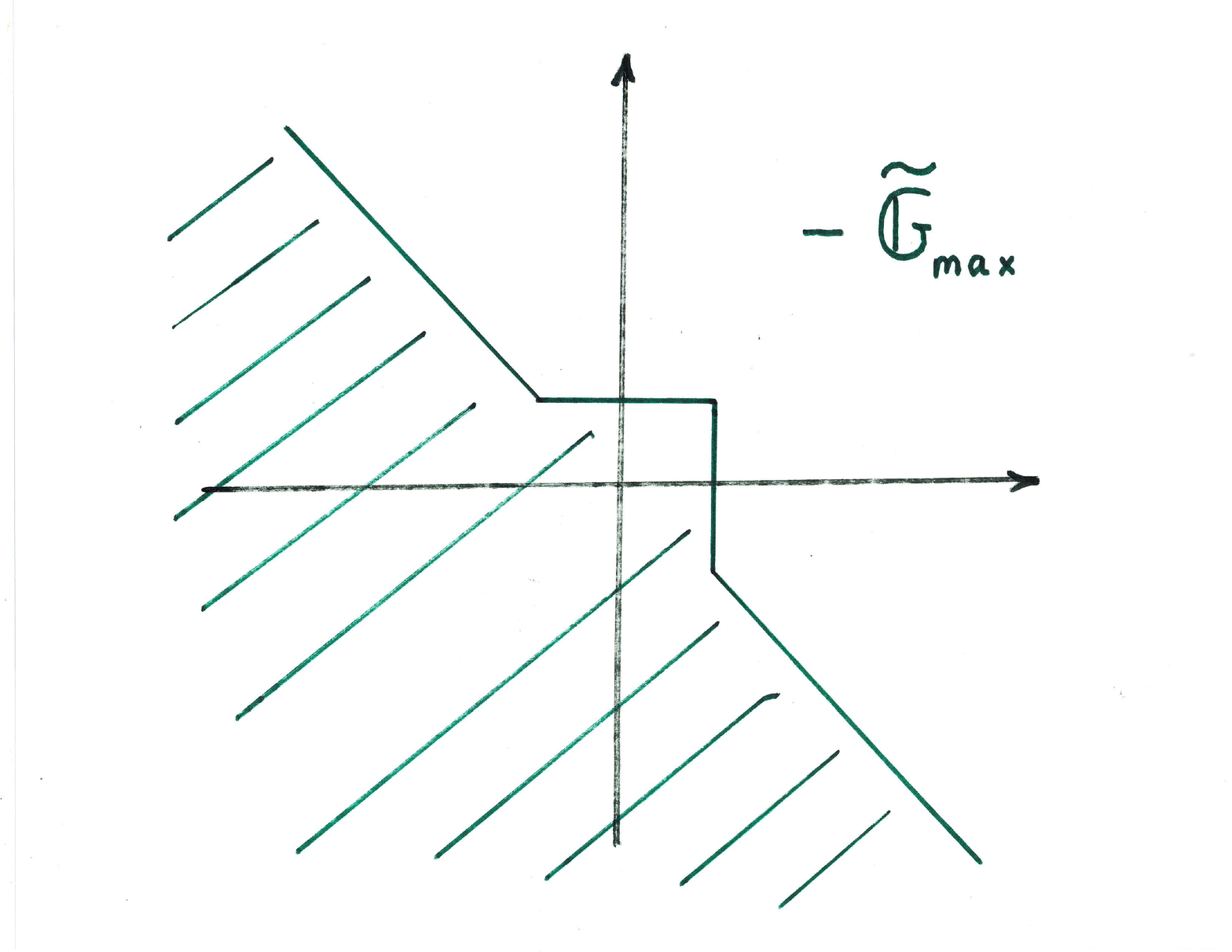}
          }

   \noindent Therefore,  the minimal generalized equation containing $\bbh$ is 
$$
\bbh^\diamondsuit \ =\ (\bbh+\cp) \cap (\bbh-\cp) \ =\ \{A : -\Id \leq A\leq \Id\} \cup \{A:\tr A=0\}
$$
\centerline{   %\hskip 1.3in
           \includegraphics[width=.3\textwidth, angle=270,origin=c]{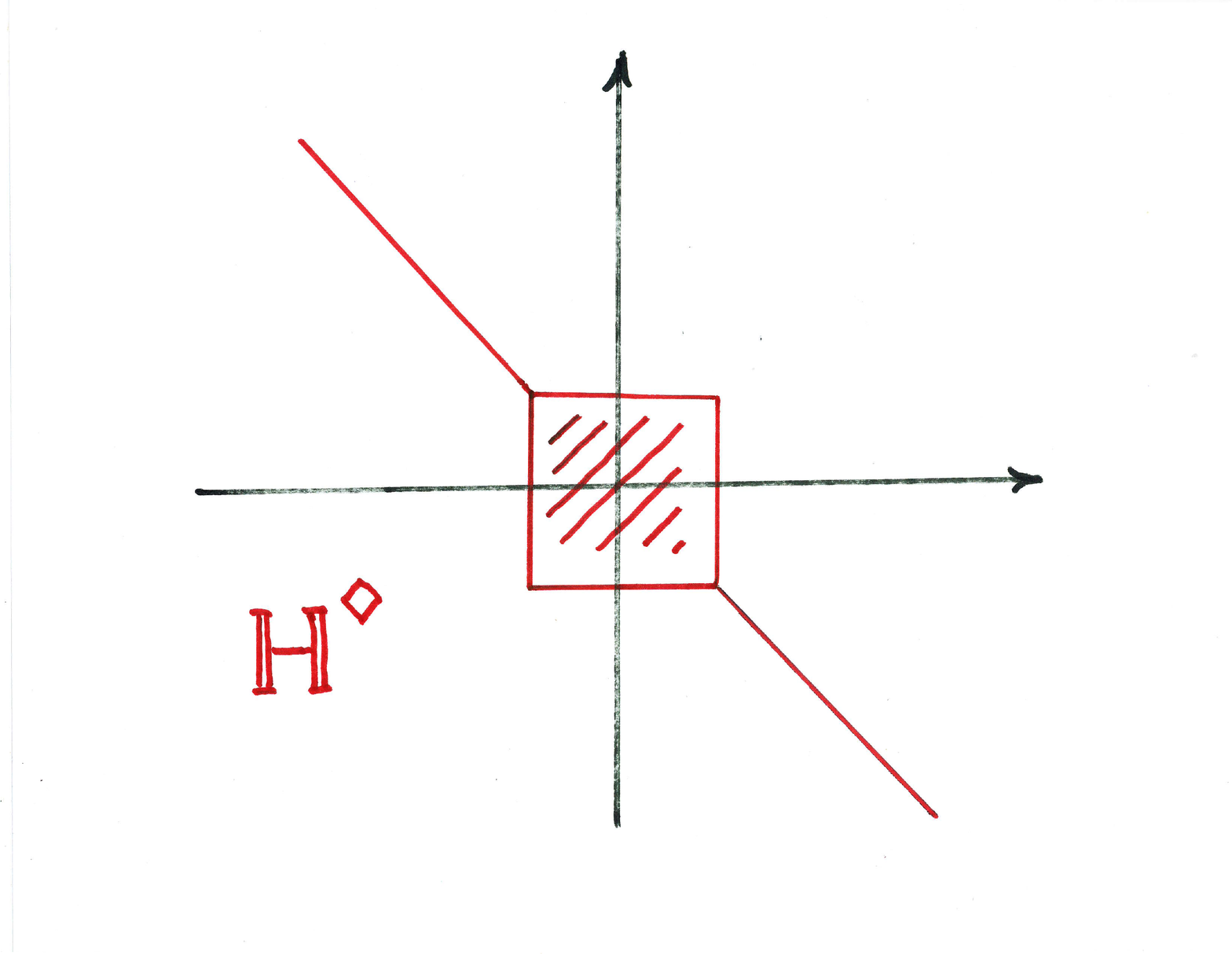}
          }

\vfill\eject

   \noindent Note that  the two subequations are:
$$
\bbem  \ =\ \{A\geq -\Id\} \, \cup \, \{\tr A\geq 0\} \and {\bbg}_{\rm max} \ =\   \{A - \Id \in\cpt\} \, \cap \, \{\tr A\geq 0\} .
$$
 {   %\hskip 1.3in
           \includegraphics[width=.3\textwidth, angle=270,origin=c]{Figure21.pdf}
          }  \hskip1.2in
          {   %\hskip 1.3in
           \includegraphics[width=.3\textwidth, angle=270,origin=c]{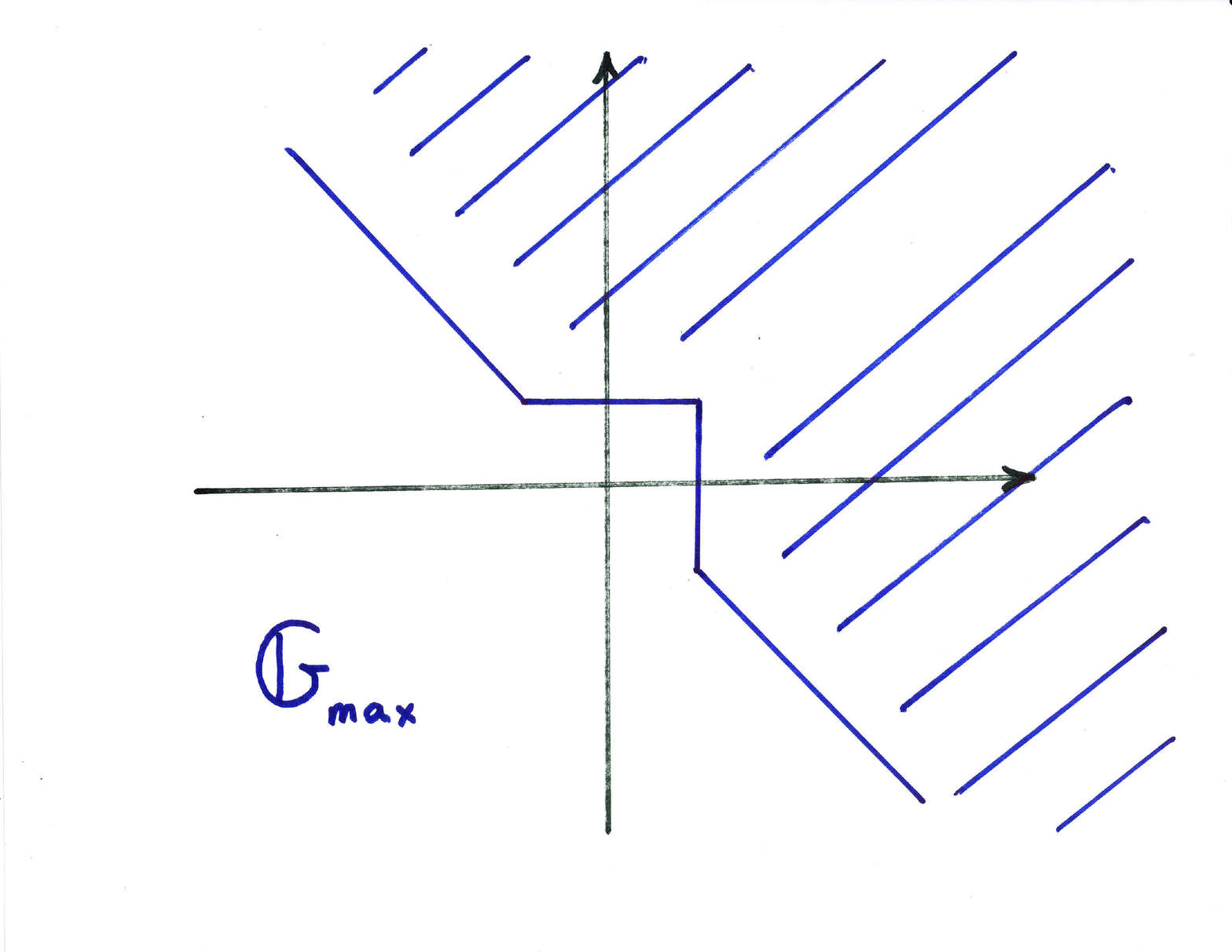}
          }
\noindent
{\bf Some Facts:}  Here $\bbgm$ and $\bbem$ are dual to one another, and so
$$
\bbh^\diamondsuit \ \equiv \ \bbem \cap (-\wt{\bbg}_{\rm max}) \ =\ \bbem \cap (-{\bbem})
$$
is a Class II generalized equation (see Section 4).

The {\bf mirror} of $\bbh^\diamondsuit$  is 
$$
(\bbh^\diamondsuit)^* \ =\ \bbgm \cap (-\bbgm) \ =\ \{  \tr A = 0\} \sim \{ -\Id < A < \Id\}.
$$
\centerline{   %\hskip 1.3in
           \includegraphics[width=.3\textwidth, angle=270,origin=c]{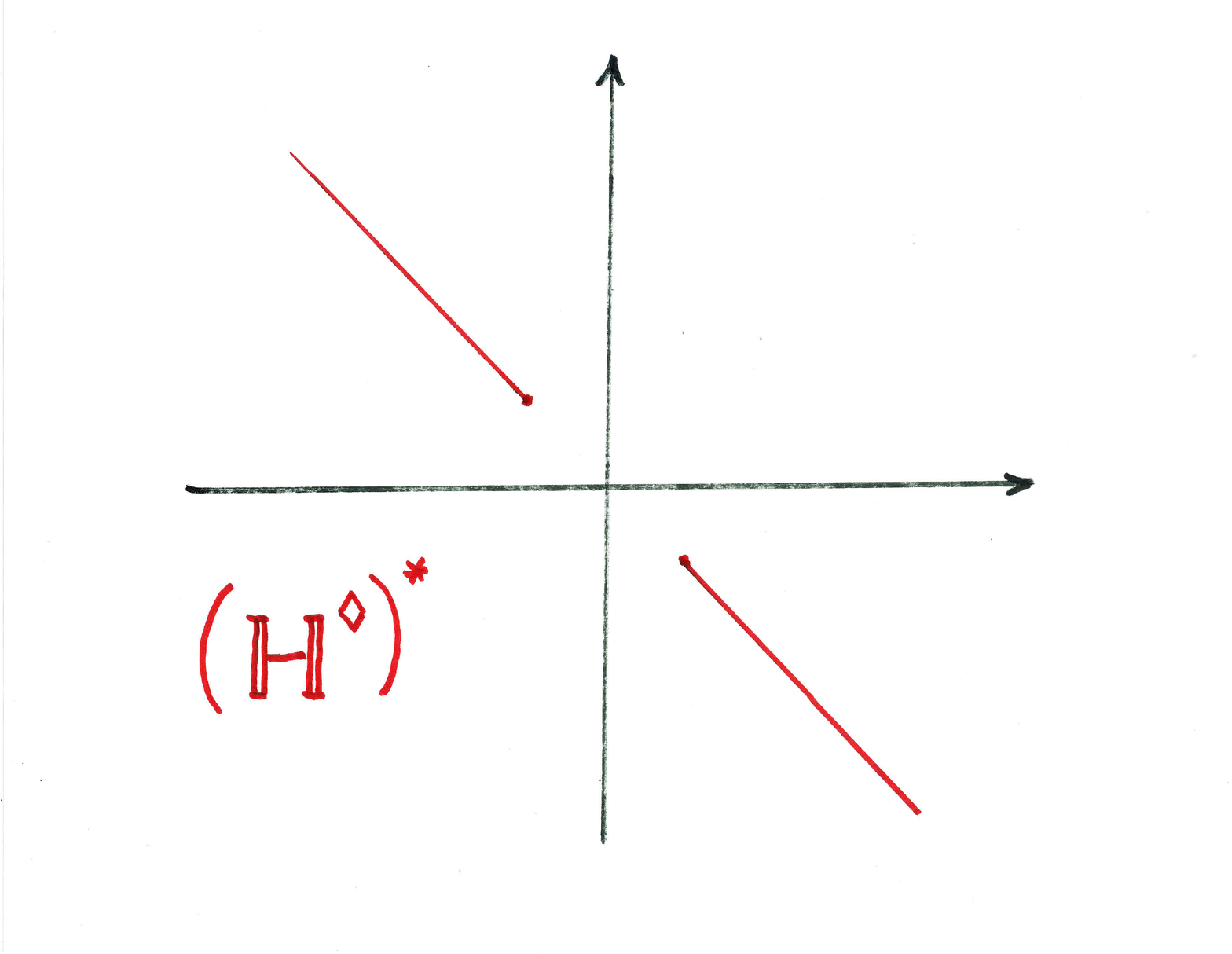}
          }
Here we have  $\Int \bbh^\diamondsuit \neq \emptyset$ and $\Int (\bbh^\diamondsuit)^*  =\emptyset$,
so that $\bbh^\diamondsuit$ is Type  III.

%%%%%%%%%%%%%%%%%%%%%%%%%%%%%%%%%%%%%%%%%%%%%%%%%%%
%%%%%%%%%%%%%%%%%%%%%%%%%%%%%%%%%%%%%%%%%%%%%%%%%%%
%%%%%%%%%%%%%%%%%%%%%%%%%%%%%%%%%%%%%%%%%%%%%%%%%%%
%%%%%%%%%%%%%%%%%%%%%%%%%%%%%%%%%%%%%%%%%%%%%%%%%%%
%%%%%%%%%%%%%%%%%%%%%%%%%%%%%%%%%%%%%%%%%%%%%%%%%%%
%%%%%%%%%%%%%%%%%%%%%%%%%%%%%%%%%%%%%%%%%%%%%%%%%%%
%%%%%%%%%%%%%%%%%%%%%%%%%%%%%%%%%%%%%%%%%%%%%%%%%%%

\vfill\eject
%\vskip.2in
\noindent
{\headfont 4.   Examples of Generalized Equations $\bbh= \bbe\cap(-\wt\bbg)$.}

 We start with some classes of examples.  First and foremost is the following.
 
 \noindent
 {\bf  Class I.  Type I Determined Equations (The $\bbg=\bbe$ Case).}  Here $\bbh = \partial \bbf = \bbf\cap (-\wt\bbf)$
 is the boundary of a subequation $\bbf$.  We refer to [\DD], [\DDR] and [\Survey] for an abundance of important
specific  examples.

Another important class of examples is

\noindent
 {\bf Class II.   (The $\bbg=\wt\bbe$ Case).}   Here $\bbh = \bbe\cap (-\bbe)$
 and $\bbh^* = \wt \bbe\cap (- \wt \bbe) =$ $  \sim [(\Int\, \bbe) \cup (-\Int\,\bbe)]$.
   Note that the overlap
 between Classes I and II consists of the boundaries $\bbh = \partial \bbf$ of self-dual subequations 
 (where the dual $\wt\bbf$ equals $\bbf$).

\noindent
 {\bf  Class IIa.    (Edges).}  In this Class II, the most basic examples are when $\bbe$ is a convex cone subequation.   Then $\bbh = \bbe\cap(-\bbe)$ is a vector subspace called the {\bf  edge} of  the cone $\bbe$.
Since $\bbe$ is a proper subspace, $\bbh$ is also a proper subspace, and hence $\Int \bbh=\emptyset$.
Now note that $\sim \bbe \ss -\Int \bbe$ is false for a proper convex cone $\bbe\ss \Symn$
because $-\Int \bbe$ is an open convex cone.   Equivalently, $\bbg \equiv \, \sim(-\Int \bbe)\ss\bbe$
is false.  By the equivalence of (1)* and (3)* in Theorem 2.6, this proves that $\Int \bbh^*\neq \emptyset$.
 Thus $\bbh$ is Type II.
 Such edge harmonics include: (i) Affine functions, where $\bbh= \{0\}$ and $\bbe=\cp$,  (ii) pluriharmonic functions in complex analysis, where $\bbe=\cp_\bbc$ and the edge $\bbh$  is SkewHerm$(\bbc^n)$, 
 (iii) $\D$-harmonic functions, and many others. The ``edge'' generalized equations are the subject of [\Edge].

Three specific non-edge Class II examples are as follows.

\vskip.2in

\centerline{ \bf   Some $\bbh$ Non-Uniqueness Examples}
 
 \noindent
 {\bf Example 4.1.  (The Quasi-convex/Quasi-concave Equation)}  Choose $r_1, r_2 \in \bbr$ with  
 $-r_1 \leq r_2$,  and let
 $$
 \bbh \ \equiv \ (\cp-r_1I) \cap (-\cp+r_2I).
 $$
Here $\bbe \equiv \cp-r_1I$ is the subequation for $r_1$-quasiconvex functions,
and $\wt \bbg\equiv \cp - r_2 I$   is the subequation for $r_2$-quasiconvex functions.
Thus $\bbh \equiv \bbe\cap (-\wt \bbg)$ is the generalized equation for functions that are
both $r_1$-quasiconvex and $r_2$-quasiconcave.
Note that $A\in \bbh \iff -r_1 I \leq A\leq r_2 I$.
{\bf A function $u$ is $\bbh$-harmonic $\iff$ $u+r_1 {|x|^2\over 2}$ is convex and 
$u-r_2 {|x|^2\over 2}$ is concave.}  

Observe now that if $u$ satisfies  this generalized  $(r_1, r_2)$ equation, then the function 
$u(x) + {\rho\over 2} |x|^2$ satisfies the  generalized  $(r_1-\rho, r_2-\rho)$ equation. Thus,
 by simply adding multiples of $|x|^2$ we translate the equation up and down the line
 $\{t\Id : t\in\bbr\}$.  So we can assume that $r_1,r_2\geq 0$, in fact let's assume 
 $r_1 = r_2 = \l\geq0$. In this case  $\wt\bbg =\bbe$ and 
$\bbh=\bbe\cap (-\bbe)$ is class II above, and Type IV.  Furthermore, there is the following
result  of Hiriart-Urruty and Plazanet in [\HP]. An alternate proof appeared in [\Eb].
 For the benefit of the reader we include a proof in Appendix A.

We say that a function is $\l-C^{1,1}$ if it is $C^1$ and
 the first derivative is Lipschitz with Lipschitz coefficient $\l$.
 
\Theorem {A.1}  {\sl For a function  $u$ on  a  convex domain $\O\ss\rn$},
$$
u \ \    \l-C^{1,1} 
\qquad\iff\qquad
{\rm both}\ \ 
\pm u \ \ \text{are}  \ \l-{\rm quasi-convex}
$$

 {   %\hskip 1.3in
           \includegraphics[width=.3\textwidth, angle=270,origin=c]{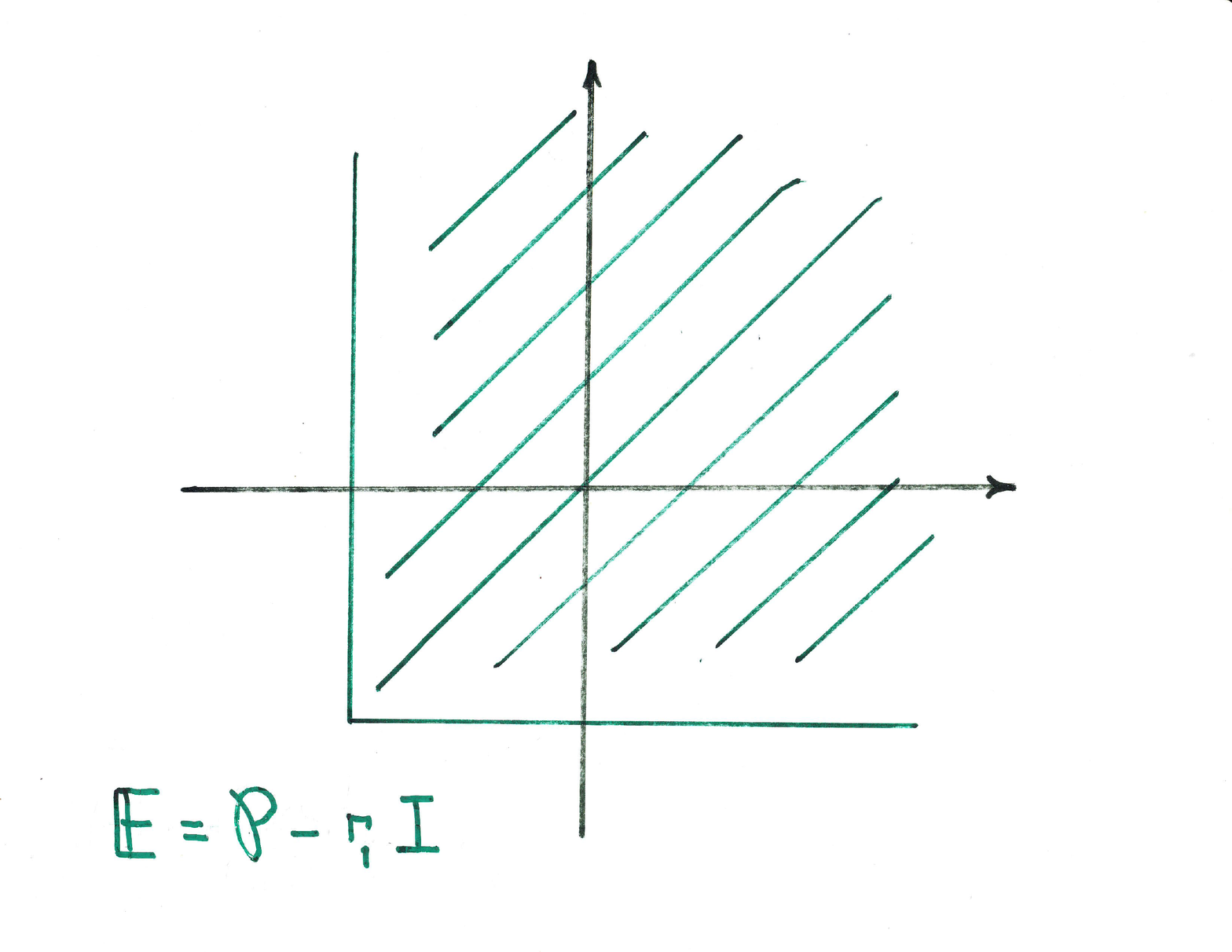}
          }  \hskip1.1in
          {   %\hskip 1.3in
           \includegraphics[width=.3\textwidth, angle=270,origin=c]{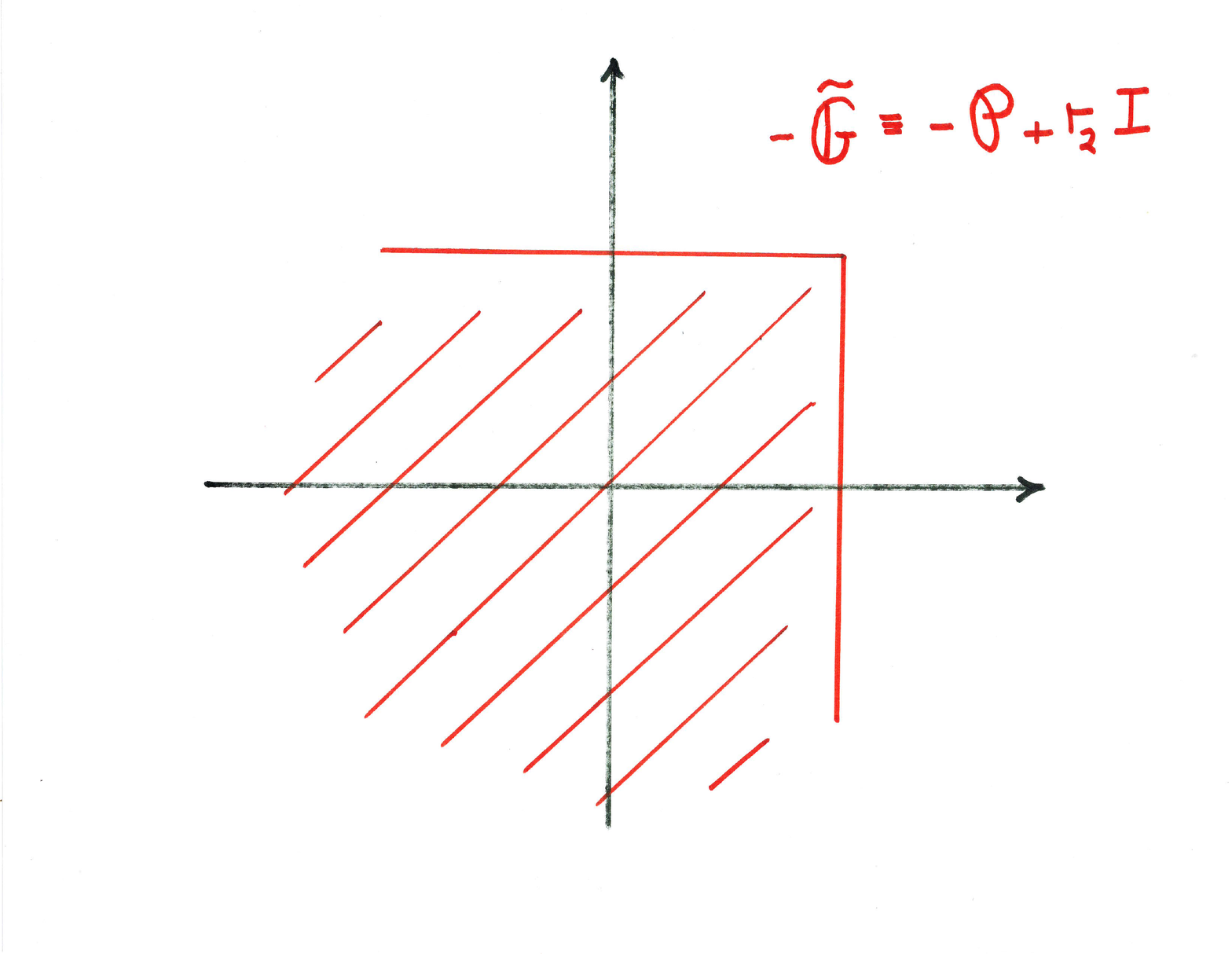}
          }

 \centerline{\includegraphics[width=.3\textwidth, angle=270,origin=c]{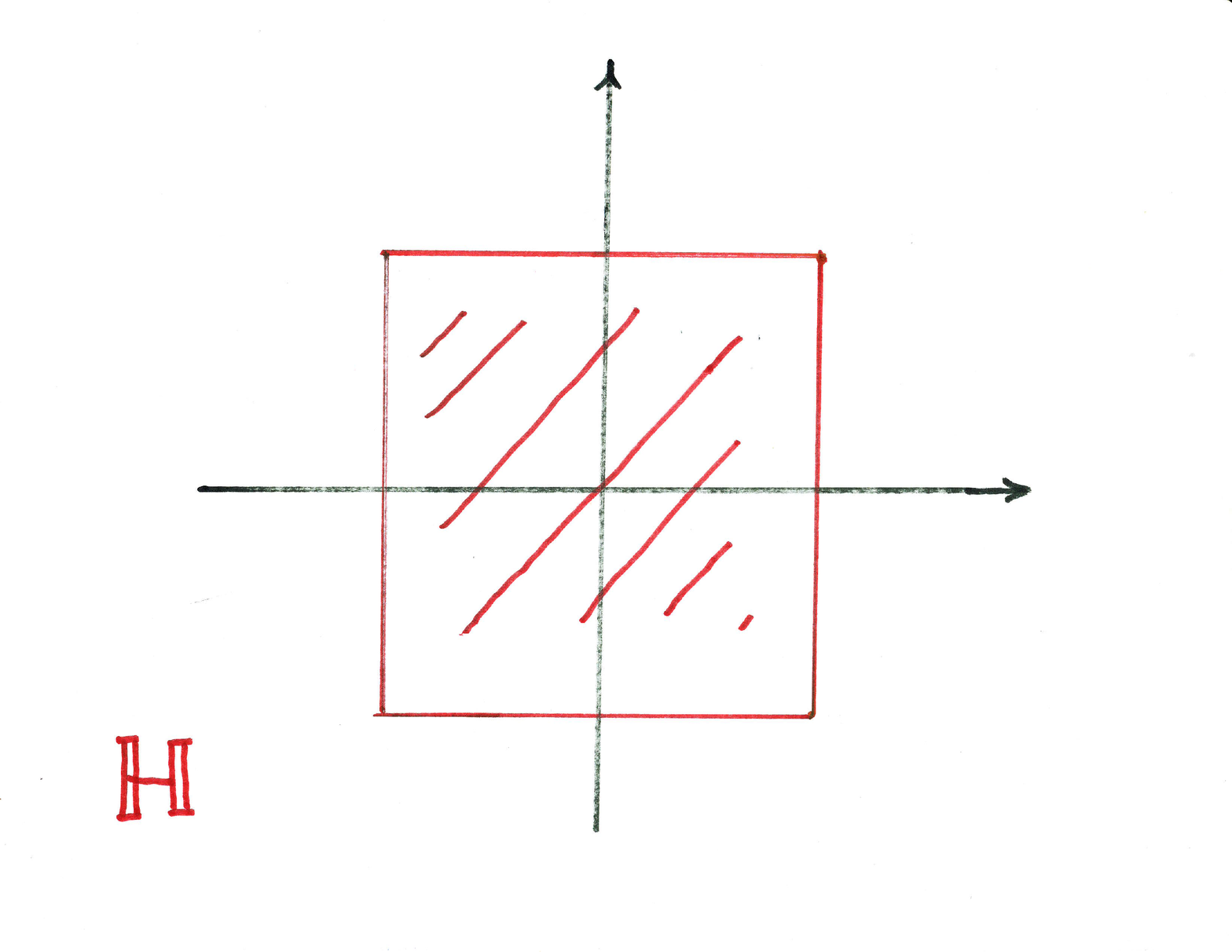}}

Here is a related example which produces the {\bf mirror equation} to the one above.

 \Ex{4.2.  (The Quasi-subaffine/Quasi-superaffine Equation)}  Choose $r_1, r_2\in\bbr$ and set
 $$
 \bbh \ \equiv \ (\cpt-r_1I) \cap (-\cpt+r_2I).
 $$
Here $\bbe \equiv \cpt - r_1 I$ is the subequation for $r_1$-quasi-subaffine functions, i.e., 
$u(x) + {r_1\over 2} |x|^2$ is subaffine, and $\wt \bbg = \cpt - r_2 I$ 
is again the subequation for $r_2$-quasi-subaffine functions.    Again  if $r_1=r_2$, then $\wt\bbg =\bbe$ and 
$\bbh=\bbe\cap (-\bbe)$ is a special case of Class II above.

If $r_1, r_2 \leq 0$, then consider $r_1'= - r_1, r_2' = - r_2$ and $\bbe' \equiv \cp- r_1'\Id$ in Example 4.1.
Then here in Example 4.2, the subequation $\bbe \equiv \cpt - r_1\Id$ is equal to $\wt{\bbe'} = \cpt+ r_1'\Id$.
Thus the mirror of $\bbh'$ in Example 4.1 is $\bbh$ in Example 4.2.
%\vfill\eject

 \centerline{\includegraphics[width=.3\textwidth, angle=270,origin=c]{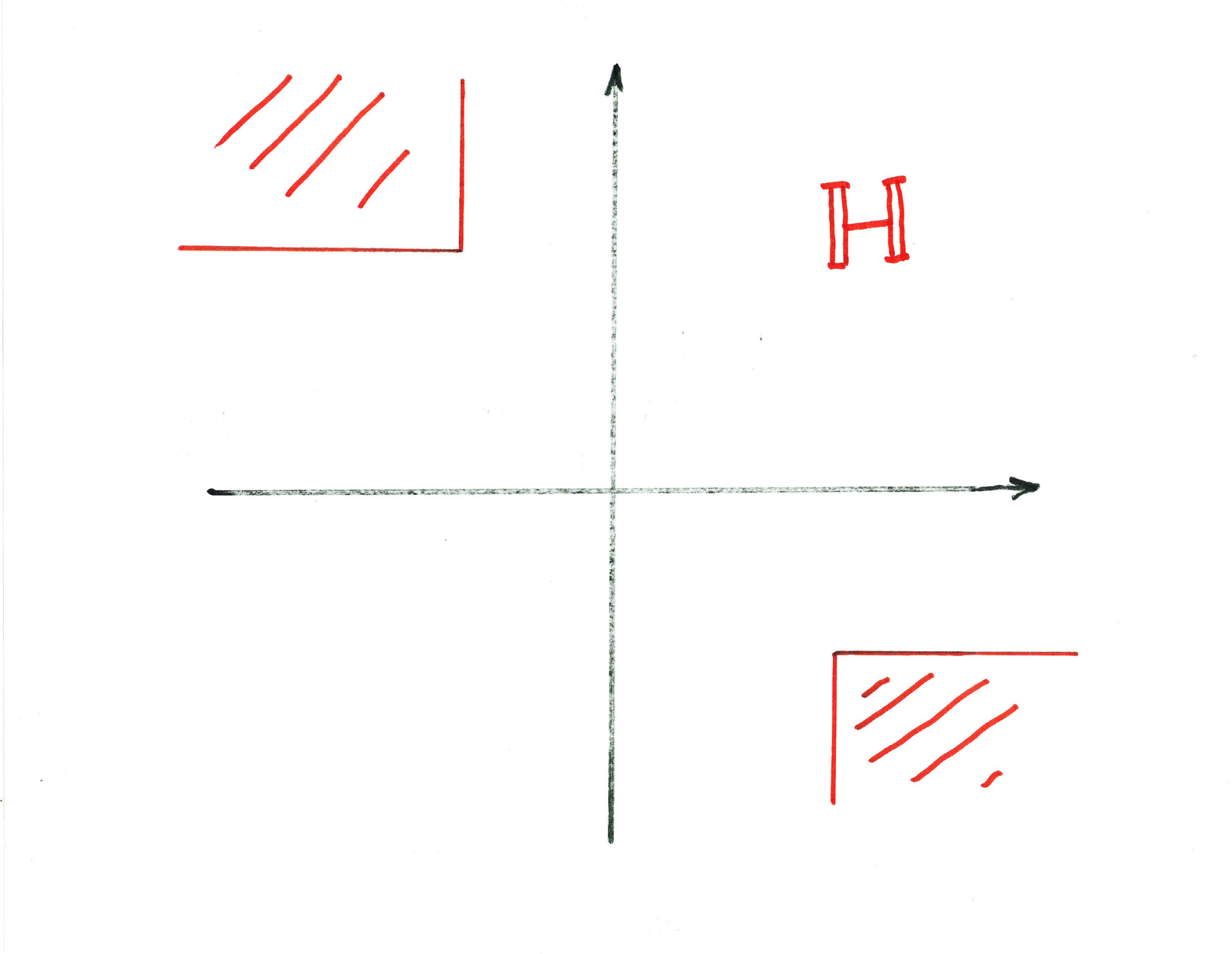}}

\noindent
{\bf The Intersection of Examples 4.1 and 4.2.}
From Example 4.1 we have
$$
\bbe = \{A \geq -\Id\} = \{ \l_{\rm min}(A) \geq -1  \}
\quad{\rm and}\quad
-\wt\bbg = \{A \leq \Id\} = \{ \l_{\rm max}(A) \leq 1  \},
$$
$$
\bbh \ =\ \{ -\Id \leq A\leq \Id\} \ =\ \{-1 \leq \l_{\rm min}(A) \ {\rm and}\  \l_{\rm max}(A) \leq 1\}.
$$
Recall that the GE in Example 4.2 is just the mirror $\bbh^* = \bbg\cap(-\wt\bbe)$ where
$$
\bbg  \ =  \ \sim \{A < \Id\} = \{ \l_{\rm max}(A) \geq 1  \}
\quad{\rm and}\quad
-\wt\bbe \  = \  \sim \{A  > -\Id\} = \{ \l_{\rm min}(A) \leq -1  \},
$$
Hence we have
$$
\bbh\cap \bbh^* \ =\ (\partial \bbe)\cap \partial (-\wt\bbg) \  = \   \{-1 = \l_{\rm min}(A) \ {\rm and}\  \l_{\rm max}(A) = 1\}
$$
{The $C^2$-harmonics for this GE are $h$'s with $ \l_{\rm min}(D^2 h)=-1$ and $  \l_{\rm max}(D^2 h) = 1$.}
This can be expressed invariantly as 
$$
D^2h +\Id \geq 0, \ \ \Id - D^2h \geq 0 \quad{\rm and}\quad \det(D^2h +\Id) + \det(\Id -D^2h) \ =\ 0,
$$
or in keeping with the Twisted Monge-Amp\`ere Example 4.8
$$
D^2h +\Id \geq 0, \ \ \Id - D^2h \geq 0 \quad{\rm and}\quad \det(D^2h +\Id) =- \det(\Id -D^2h).
$$

\noindent
{\bf Note:}   There are lots of such harmonics. For example,
$$
h(x_1, x_2, y) \ =\ \half x_1^2 - \half x_2^2 + u(y)
\quad {\rm where}\quad
-\Id \ \leq\ D^2u \ \leq\ \Id.
$$

\vskip.2in

\centerline{ \bf   Some Type II Examples of  Non-Existence and Uniqueness for $\bbh$}

 \Ex{4.3.  (Generalized Constrained Laplacians)} 
 
 \Def{4.4}  A closed subset $\bbh$ of $\partial \D \equiv \{A: \tr A=0\}$
 (with $\emptyset \neq \bbh\neq \partial \D$) will be called a {\bf (generalized) constrained Laplacian}.
 
 \Lemma{4.5}  {\sl
 If $\bbh$ is a (generalized) constrained Laplacian, then $\bbe\equiv \bbh+\cp$ 
 is closed and hence a subequation.  Furthermore, $\bbh-\cp$ is  closed, and
 $$
 \bbh\ =\ (\bbh+\cp) \cap (\bbh-\cp),
 $$
 so that $\bbh$ is a generalized equation with $\bbem=\bbh+\cp$ and $- \wt {\bbg}_{\rm max}  =\bbh-\cp$.
 }
 
 \noindent
 {\bf Proof.}
 Suppose $B\in\overline{\bbh+\cp}$, i.e.,  $B=\lim_j(A_j+P_j)$ with $A_j\in\bbh$
 and $P_j\geq0$.  Set $B_j=A_j+P_j$.  Since $\bbh\ss \{\tr A=0\}$, we have 
 $\tr B_j = \tr(A_j+P_j) = \tr P_j$.  Thus $\tr P_j \to \tr B$.  The set $\{P\in\cp: \tr P\leq c\}$ with $c>0$
 is compact.  Therefore we can extract a convergent subsequence of $\{P_j\}$, which we again call 
$\{P_j\}$, with $P_j\to P\geq0$.  Therefore, $A_j = B_j- P_j \to B-P\equdef A$.
Since $\bbh$ is closed and each $A_j\in\bbh$, we have $A\in\bbh$.  
Thus $B=A+P$ with $A\in\bbh$ and $P\geq 0$, and so we have proved that $\bbh+\cp$ is closed.
Replacing $\bbh$ by $-\bbh$ we have that $-\bbh+\cp = -(\bbh-\cp)$ is closed. Thus $\bbh-\cp$ is closed.

Obviously $\bbh\ss(\bbh+\cp) \cap (\bbh-\cp)$. Suppose now that $B\in (\bbh+\cp) \cap (\bbh-\cp)$,
i.e., $B=A+P =A'-P'$ for $A,A'\in \bbh$ and $P,P'\in\cp$.  Then $\tr B =\tr P = - \tr P'$.  Since
$P, P'\geq0$, this implies that $P=P'=0$, and hence $B=A=A'\in\bbh$.\qed

Note that $\Int \bbh=\emptyset$, so $\bbh$ must be type I or type II.  However, it cannot be
type I without $\bbh=\partial \D$ (cf. Propositions 2.18 and 2.19).  This proves
$$
\text{\sl A generalized constrained Laplacian is of type II.}
$$

\noindent
{\bf Question 4.6.} Suppose $\bbe$ and $\bbg$ define $\bbh\  (=\bbe\cap(-\wt\bbg))$
and $\bbh \ss\partial\D$ is a generalized constrained Laplacian.
If $h$ is $H_{\bbe,\bbg}$-harmonic, then is $h$ $\D$-harmonic?

\Ex{4.7} \ Let $\rn=\bbr^k\oplus\bbr^\ell$ and
$$
\bbh \equiv \left \{A\equiv \left( \begin{matrix} a& c \\ c^t & b \end{matrix} \right) : \tr A=0, \ a \geq0
\ {\rm and}\ b\leq 0\right\}
$$
with $a\in \Sym(\bbr^k),\  b\in \Sym(\bbr^\ell)$.  Then
$$
\begin{aligned}
\bbe_{\rm min} \ =\ \{A : \  &a\geq 0 \ {\rm and}\  \tr A\geq 0\}, \qquad
\bbg_{\rm max} \ =\ \{A : b\geq 0\} \cup \bbe_{\rm min}, \\
\ \ &{\rm and}\ \ 
  \wt \bbg_{\rm max}  = \{A : b\geq0\ \ {\rm and}\ \ \tr A\geq0\}.
  \end{aligned}
$$

One could also look at this from the universal eigenvalue point of view
 (see the subsection of Section 5 in [\Hyp] and Remark 4.12 below), 
 by taking the eigenvalues  of $a$ and $b$.  Let  $Q^+(\bbr^k)$ denote
  the positive orphant defined by $x_j\geq 0$ for all $j$, and let
$Q^-(\bbr^k)$ be  similar with coordinates $y_i\leq 0$.  Then $\bbh, \bbem,\bbgm$ can be defined by
$$
H \ =\ \{(x,y) : x\in Q^+(\bbr^k), y\in Q^-(\bbr^\ell), \tr(x,y)=0\}
$$
$$
E_{\rm min} \ =\ \{(x,y) :  x\in Q^+(\bbr^k),  \tr(x,y)\geq 0\},
$$
$$
G_{\rm max} \ =\ \{(x,y) :  y\in Q^+(\bbr^\ell)\} \cup E_{\rm min}, \qquad
\wt G_{\rm max} \ =\ \{(x,y) :  y\geq 0 \ \ {\rm and}\ \ \tr(x,y)\geq0)\}.
$$

A great example of a non-existence/uniqueness $\bbh$ equation (Type II)
 has been introduced and studied in  [\ST] and [\SW].
This is discussed next.

\Ex{4.8. (The Universal Version of the  Twisted Monge-Amp\`ere  Equation)}
The real twisted Monge-Amp\`ere equation is defined by $\bbh \ss\Sym(\bbr^k\times\bbr^\ell)$  
consisting of all 
$$
\left(\begin{matrix} A & C \\ C^t & B  \end{matrix}\right) \quad {\rm such\ that\ } \ A\geq 0, \ \ B\leq 0
\ \ {\rm and} \ \  \log\, \det A - \log\, \det (-B) =0
$$
i.e.,  $\det A = \det (-B)$, or $\det A = |\det B|$.

As in Example 4.7, the universal version  of this equation is defined  on $\rn = \bbr^k\times \bbr^\ell$   by
$$
H \ \equiv\  \{(x,y) : x\in Q^+(\bbr^k), y\in Q^-(\bbr^\ell) \ {\rm and}\ x_1\cdots x_k=|y_1\cdots y_\ell|\}.
$$

\Lemma {4.9}  {\sl
Let $E \equiv H+Q^+(\rn)$.  Then $E$ is fibred over $Q^+(\bbr^k)$, 
where the fibre $E_x$ of $E$ at $x\in Q^+(\bbr^k)$ is the dual MA universal subequation:}
$$
\wt P_{x_1\cdots x_k}(\bbr^\ell) \ =\ (\sim Q^-(\bbr^\ell)) \cup \{ y\in Q^-(\bbr^\ell) : |y_1\cdots y_\ell|\leq x_1\cdots x_k   \}. 
$$

Since it is easy to see that $E$ is closed, $E$ is equal  to  the minimal subequation defined above
for this $H$.

\noindent
{\bf Proof of Lemma 4.9.}  First note that $H$ is fibred over $Q^+(\bbr^k)$ with fibre $H_x$ at 
$x\in Q^+(\bbr^k)$ given by 
$$
H_x \ =\ \{y\in Q^-(\bbr^\ell) : |y_1\cdots y_\ell|=x_1\cdots x_k\}.
$$
Second note that this equals
$$
H_x = \partial \wt P_{x_1\cdots x_k}(\bbr^\ell)
$$
the boundary of the dual MA-subequation at level $c=x_1\cdots x_k$.  Third note that,
since $\wt P_{x_1\cdots x_k}(\bbr^\ell)$ is a subequation, 
$$
H_x + Q^+(\bbr^\ell) = 
 \partial \wt P_{x_1\cdots x_k}(\bbr^\ell) + Q^+(\bbr^\ell) = \wt P_{x_1\cdots x_k}(\bbr^\ell)
$$
Defining $E'$ by its fibres $E_{x_1\cdots x_k}' \equiv  \wt P_{x_1\cdots x_k}(\bbr^\ell)$ 
over $x \in Q^+(\bbr^k)$, we have $E'\ss E$, and it remains to show $E\ss E'$.
But $H\ss E'$, so it is emough to show $E'$ is $Q^+(\rn)$-monotone.
As noted above $E'$ is  $Q^+(\bbr^\ell)$-monotone since $\wt P_{x_1\cdots x_k}(\bbr^\ell)$
is a subequation. Now increasing one of the $x$ coordinates with $x\in Q^+(\bbr^k)$
increases $\wt P_{x_1\cdots x_k}(\bbr^\ell)$ proving that $E'$ is $Q^+(\bbr^k)$-monotone.
Finally the orphant $Q^+(\rn)$ equals the product $Q^+(\bbr^k)\times Q^+(\bbr^\ell)$. \qed

\Lemma {4.10}  {\sl
Let $\wt G \equiv -H+Q^+(\rn)$.  Then $\wt G$ is fibred over $Q^+(\bbr^\ell)$, 
where the fibre $\wt G_y$ of $\wt G$ at $y\in Q^+(\bbr^\ell)$ is the dual MA universal subequation:}
$$
\wt G_y \ =\ \wt P_{|y_1\cdots y_\ell|}(\bbr^k).
$$

The proof of Lemma 4.10 is similar to the one for Lemma 4.9, and is skipped.

\Prop{4.11} {\sl
$H= E\cap(-\wt G)$ is a universal version of a generalized equation with minimum subequation
$E$ and maximum subequation $G$.}

\pf
Note that $(x,y)\in E\  \iff\ x\in Q^+(\bbr^k)$ and $y\in \wt P_{x_1\cdots x_k}(\bbr^\ell)$
by Lemma 4.9.  Note also that
$$
(x,y) \in -\wt G \quad \iff \quad x\in -\wt P_{|y_1\cdots y_\ell|}(\bbr^k) \ \ {\rm and}\ \ y\in Q^-(\bbr^\ell).
$$
by Lemma 4.10.

Now assume $(x,y) \in E\cap (-\wt G)$.  Then $x\in Q^+(\bbr^k)\cap(-\wt P_{|y_1\cdots y_\ell|}(\bbr^k))$
or otherwise said,  $x\in Q^+(\bbr^k)$ and $|y_1\cdots y_\ell| \leq x_1\cdots x_k$.
Also,   $y\in Q^-(\bbr^\ell)\cap\wt P_{x_1\cdots x_k}(\bbr^\ell)$
or otherwise said,  $x\in Q^-(\bbr^\ell)$ and $x_1\cdots x_k \leq |y_1\cdots y_\ell|$.

In summary, if $(x,y)\in E\cap (-\wt G)$, then
$$
x\in Q^+(\bbr^k), \ \ y\in Q^-(\bbr^\ell), 
 \ \   {\rm and}\ \ x_1\cdots x_k  =  |y_1\cdots y_\ell|
$$
that is,  $(x,y)\in \H$.  It is easy to see that $H\ss E\cap (-\wt G)$.\qed

\Remark {4.12. (Universal Equations and G\aa rding/Dirichlet Operators)}  
A closed subset $\L\ss\rn$ which is symmetric under permutations of the coordinates
and satisfies $\L +\rn_+ \ss\L$ is called a {\bf universal eigenvalue subequation}.
There is an obvious one-to-one correspondence between subequations $\bbf\ss\Symn$, which depend only
on the eigenvalues of $A\in\bbf$, and universal subequations $\L\ss \rn$.  However, this $\bbf$
is only one of many subequations determined by $\L$ which are constructed by substituting G\aa rding eigenvalues
for regular eigenvalues as follows.

Let $g$ be a homogeneous polynomial of degree $n>0$, on some $\Sym(\bbr^m)$, which satisfies the conditions
of being a {\bf G\aa rding/Dirichlet}, or {\bf GD}, operator (as defined in [\Hyp, \S 5]).  Then for each $A\in \Sym(\bbr^m)$ this operator
gives $n$ eigenvalues $\l_g(A)$, and so $\L$ determines a subequation in $\bbr^m$ by:
$$
\bbf^g_\L \ =\ \{A \in  \Sym(\bbr^m) : \l_g(A) \in \L\}.
$$
For example, the universal subequation $\L\equiv \{\l\in\rn : \l_j\geq0 \ \forall\, j\}$
determines the {\sl G\aa rding Monge Amp\`ere subequation} 
$\bbf^g_\L = \{ A\in  \Sym(\bbr^m) : \l_{g,1}(A), .. ,  \l_{g,n}(A)\  \text{are all $\geq 0$}\}$,
which is just the closed G\aa rding cone for $g$.

This carries over to generalized equations. 
A {\bf generalized universal equation} is any closed $\L'\ss\rn$ which is an intersection involving two
two universal subequations $\L' = \L_1\cap(-\wt\L_2)$.
For example, the universal Laplacian
$\L = \{\l\in\rn : \l_1+\cdots+\l_n \geq0\}$ determines a G\aa rding Laplacian $\bbf^g_\L$ 
for each GD operator of degree $n$.  Moreover, given a pair of GD operators $g_1, g_2$
of degrees $n_1+n_2=n$,  
one has a  constrained  Laplacian generalized equation induced by the universal version of the constrained Laplacian given in Example 2.3.
Namely, we have
$$
\bbh \ \equiv\   
 \{ A\in \Sym(\bbr^m) : \l_{g_1,j}(a) \geq 0, \l_{g_2,k}(b)\leq0, \ {\rm and}\ \sum_j \l_{g_1,j}(a)
+ \sum_k \l_{g_2,k}(b) =0\}.
$$
\noindent
{\bf Example 4.13. (Twisted G\aa rding MA Generalized Equations).}
Similarly (we leave this to the reader) the universal twisted MA-equation (Example 4.8) spawns 
a huge family of generalized equations.  For instance, in addition to the real version in [\ST], [\SW], one has
a complex  version, a quaternionic version, a Lagrangian version, branched versions of these three,
elementary symmetric versions of these three (the so-called ``hessian equation'' versions),  just to name 
a few.

The Examples 4.1 and 4.2 can also be viewed as ``universal subequations'', spawning many more examples of generalized equations as above.

Since we have no reason to rule out $\bbf=\emptyset$ or $\bbf=\Symn$ as a subequation in this paper,
we have that $\bbh$ equal to plus or  minus a subequation is  an example of  a generalized equation of Type III  or II respectively.

\noindent
{\bf Example 4.14. (Subequations as Generalized Equations).}
For any subequation $\bbe \neq \emptyset$, if we choose $\bbg=\emptyset$, i.e., $-\wt\bbg=\Symn$, then 
$\bbh = \bbe \cap \Symn =\bbe$ is a generalized equation.  
Now since $\overline{\Int\, \bbe}=\bbe$ and $\bbe \ne \emptyset$, we have $\Int\, \bbh\ne \emptyset$.
Also the mirror $\bbh^* = \bbg\cap(-\wt\bbe)=\emptyset$. Hence, $\Int\, \bbh^*=\emptyset$.  In summary, 
if $\bbe\neq\emptyset$ is any subequation, then with $\bbg=\emptyset$, we have $\bbh=\bbe$ and
$\bbh^*=\emptyset$, so that $\bbe$ itself (not $\partial\bbe$) is a generalized equation which falls
in the Existence/Non-Uniqueness case for $\bbh$ (Type III).
Similarly, $-\bbf = \bbe\cap (-\wt\bbg)$ with $\bbe \equiv \Symn$ and $\bbg = \wt\bbf$ is Type II.

\vskip.2in

\centerline{ \bf   Some Type IV Examples of  Non-Uniqueness/Non-Existence.}
 
 \Ex{4.15}   With coordinates $z=(x,y) \in \rn = \bbr^k\times \bbr^\ell$,
 define $\bbe$ by $D_x^2 u\geq0$ and $\wt\bbg$ by
 $D_y^2 u\geq0$ (so that  $\bbg$ is the subaffine subequation on $\bbr^\ell$ considered as a subequation on $\rn$).
 Then with $\bbh \equiv \bbe\cap(-\wt\bbg)$ we have that the $\bbh$-harmonics are continuous functions
 $h(x,y)$ that are separately convex in $x$ and concave in $y$.  
 The mirror $\bbh^*$-harmonics are continuous functions
 $h^*(x,y)$ that are separately subaffine in $x$ and superaffine in $y$.

\vskip.3in

\centerline{\bf Elementary Examples from $\cp$ and $\cpt$}

\Ex{4.16} 

(a) If we take $\bbe=\bbg=\cp$, then $\bbh=\partial \cp$ is the determined equation
whose harmonics are solutions to the real Monge-Amp\`ere equation.

(b)  If we take $\bbe=\bbg=\cpt$, then $\bbh=\partial \cpt = -\partial \cp$
is the determined subaffine equation, whose harmonics are the negatives of the
harmonics for real Monge-Amp\`ere equation.

(c)  If we take $\bbe=\cp$ and $\bbg= \cpt$, then $\bbh= \cp\cap(-\cp) = \{0\}$, and the harmonics
are functions which are both convex and concave, i.e.,  the affine functions.  So far nothing is really new.

(d)    If we take $\bbe=\cpt$ and $\bbg= \cp$, then $\bbh= \cpt\cap(-\cpt) 
= \Symn \sim [(\Int \cp) \cup (-\Int \cp)]$.   The solutions are functions $u$ with $D^2u$ non-definite,
(i.e., never $>0$ nor $<0$) in the $C^2$  case.

\vskip.3in

\centerline{\bf The Affine Generalized Equation}

\Ex{4.17} 
Here we are interested in equations where $\bbh_{\bbe, \bbg}=\{0\}$.
Note that  $\bbh_{\bbe, \bbg} = \bbe\cap(-\wt\bbg) = \bbe \sim (\Int \bbg) = \{0\} \iff \bbe\sim\{0\} \ss\Int\bbg$.
Therefore
$$
\bbh_{\bbe, \bbg}\ =\ \{0\} \qquad\iff\qquad \bbe\sim \{0\} \ss\Int\bbg\ \ {\rm and}\ \ 
0\in\partial \bbe, 0\in\partial \bbg.
$$
\Prop{4.18}  {\sl
Suppose $\bbe$ is a convex cone and let $\bbe^0 \equiv \{B : \bra BA\geq0 \, \forall\, A\in\bbe\}$
be its polar cone.    

Suppose  $\exists \,A\in\Int \bbe^0$ (equivalently $\bbe\ss \D_A = \{B: \bra AB\geq 0\}$) \\such that
$\D_A \ss\bbg$  (equivalently $\wt \bbg\ss \wt \D_A  = \D_A$).  Then $h$ is $\bbh_{\bbe,\bbg}$-harmonic
$\Rightarrow$ $h$ is $\D_A$-harmonic $\Rightarrow$ $h$ is affine.
}

\noindent
{\bf Proof.}
$h$ is $\bbh_{\bbe,\bbg}$-harmonic $\iff$ $h$ is $\bbe (\ss \D_A)$ subharmonic and 
$-h$ is $\wt\bbg (\ss\D_A)$ subharmonic.  Therefore, $h$ is  $\bbh_{\bbe,\bbg}$-harmonic
 $\Rightarrow h$ is $\D_A$-harmonic
$\Rightarrow h$ is smooth $\Rightarrow h$ is affine.\qed

\Prop{4.19}  {\sl
Any affine generalized equation $\bbh_{\bbe,\bbg}=\{0\}$ is Type II, i.e., $\Int\bbh = \emptyset$ and 
$\Int\bbh^*\neq \emptyset$.
}

\noindent
{\bf   Proof.}  Of course $\Int \bbh =\Int\{0\}=\emptyset$.  Now
$\Int \bbh^* = (\Int\bbg)\sim \bbe = \emptyset \iff \Int \bbg\ss\bbe  \iff \bbg\ss\bbe$, which proves 
that $\bbe=\bbg$ is Type I. Hence  $\bbh=\partial \bbe$ and so  $\partial \bbe=\{0\}$.  This is impossible
for a subequation $\bbe$,
so $\Int \bbh^*\neq \emptyset$.\qed

Recall for Type II that $\bbh=\partial \bbe \cap \partial \bbg = \{0\}$.

The {\bf canonical pair} for $\bbh=\{0\}$ is given as follows.

\noindent
$\bbem = \{0\}+\cp = \cp, 
\quad 
-\wt{\bbg}_{\rm max} = \{0\}-\cp = -\cp,
\quad
\wt{\bbg}_{\rm max} = \cp,
\quad
{\bbg}_{\rm max} = \cpt$.

We have $\bbh=\{0\}=\cp\cap(-\cp)$.  

In this case an $\bbh_{\cp, \cpt}$ harmonic $h$ is affine since $\pm h$ are both convex.

\vfill\eject

%\vskip.2in
\noindent
{\headfont 5.   Unsettling  Questions.}

In Section 8 of [\IDP] we posed several such questions, starting with the single-valuedness
of operators and the following equivalent  restatement of that question.

\noindent
{\bf (CCQuest) Constant Coefficient Subequation Question:}  Can a pair of subequations $\bbe, \bbg$, with 
disjoint equations, i.e., $\partial \bbe \cap \partial \bbg=\emptyset$, have a simultaneous harmonic $h$?

Of course a simultaneous harmonic  $h$ cannot be $C^2$ since one would have 
$D^2_x h \in \partial\bbe \cap \partial \bbg=\emptyset$.  One reason for the success of
viscosity theory is that the intuition gained from examining classical situations carries over
to the viscosity approach.  In fact, one can show that any  simultaneous harmonic 
must be quite bizarre, but this is short of non-existence.
(A known result which has some kinship to this open question is the  fact that an  arbitrary
upper semi-continuous function has an upper test function at a dense set of points, 
cf.  [\NOTES, Lemma 6.1$'$].)

In this section we first discuss some equivalent versions of the question above. 
We then examine an  extension of this question to a natural one for any generalized equation.

It is easy to see from positivity that the boundary of a  subequation can be expressed as the graph
of a continuous function over the hyperplane $\{\tr A=0\}$.
Consequently, if $\partial \bbe \cap \partial \bbg = \emptyset$, we might as well assume
$\bbe\ss\bbg$.  Now the hypothesis of (CCQuest) can be reformulated as follows.
$$
\partial \bbe \cap \partial \bbg = \emptyset \ \ {\rm and}\ \ \bbe\ss\bbg
\quad\iff\quad 
\bbe \ss \Int\bbg\quad\iff\quad 
\partial\bbe \ss \Int\bbg.
\eqno{(5.1a)}
$$
We leave the proof to the reader.
The condition $\bbe\ss\Int \bbg$ in (5.1a) is obviously equivalent to:
$$
\text{The generalized equation $\bbh\equiv \bbe\cap (-\wt\bbg) = \bbe\cap(\sim \Int \bbg)$ is empty.}
\eqno{(5.1b)}
$$

\Lemma{5.1} {\sl
Under the equivalent assumptions of (5.1)}
$$
\text{$h$ is $\bbh_{\bbe,\bbg}$-harmonic
\ \ $\iff$ \ \ 
$h$ is both $\partial\bbe$- and  $\partial\bbg$-harmonic.}
\eqno{(5.2)}
$$

\noindent
{\bf Proof.}
First suppose $h$ is $\bbh_{\bbe,\bbg}$-harmonic, i.e., $H$ is $\bbe$-subharmonic and
$-h$ is $\wt\bbg$-subharmonic.  By (5.1) $\bbe\ss\bbg$ and hence $\wt\bbg\ss\wt\bbe$ also.
Thus $h$ is both $\bbe$ and $\bbg$ subharmonic and $-h$ is 
 both $\wt\bbg$ and $\wt\bbe$ subharmonic, which proves that $h$ is both 
 $\partial\bbe$ and $\partial\bbg$ harmonic.
 
 For the converse suppose that $h$ is both   $\partial\bbe$ and $\partial\bbg$ harmonic.
 Then $h$ is $\bbe$ subharmonic and $-h$ is $\wt\bbg$-subharmonic, so that 
 $h$ is $\bbh_{\bbe,\bbg}$-harmonic.\qed

This equivalence (5.2) means that we can reformulate the above (CCQuest)
concerning simultaneous harmonics for $\partial\bbe$ and $\partial \bbg$ as follows.

\noindent
{\bf (CCQuest)$'$:}  Does there exist a subequation pair $\bbe, \bbg$ 
defining the generalized empty equation $\bbh=\emptyset$ with the property that 
$\bbh_{\bbe, \bbg}$ has a harmonic?

\noindent
{\bf Summary.} There are lots of subequations pairs  $\bbe, \bbg$ defining this generalized empty equation
$\bbh = \emptyset$.  For some of these pairs we can prove that $\bbe\cap (-\wt\bbg)$ has no 
harmonics.  For example, this holds if $\bbe$-harmonics and $\bbg$-harmonics are always
$C^2$ because of Lemma 5.1. To conclude we note that for $\bbh=\emptyset$ the canonical
defining pair is $(\emptyset, \emptyset)$, so that for
 this pair $\bbe_{\rm min}, \bbg_{\rm max}$ defining $\bbh=\emptyset$ there are also no harmonics.

Now we can broaden our question as follows.

\noindent
{\bf Broadened  Equation Question:}  Given a generalized equation $\bbh$ does there
exists a subequation pair $\bbe, \bbg$ defining $\bbh$ so that $\bbh \equiv \bbe\cap(-\wt \bbg)$
has a harmonic which is not a harmonic for $\bbh\equiv \bbe_{\rm min}\cap ( - \wt  \bbg_{\rm max})$?  

Note that by Proposition 3.1  $\bbe_{\rm min}\ss\bbe$ and $ - \wt  \bbg_{\rm max} \ss -  \wt  \bbg$
so that $\bbh\equiv \bbe_{\rm min}\cap ( - \wt  \bbg_{\rm max})$ harmonics are always 
$\bbh\equiv \bbe\cap ( - \wt  \bbg)$ harmonics.

As for (CCQuest), any   such harmonic $h$ in this Broadened  Equation Question must be weird and pathological,  much worse than $C^2$ for sure.

\vfill\eject
%\vskip.2in
\noindent
{\headfont 6.  The  General Case of the Main Theorem.}

For clarity and simplicity we have been restricting attention to pure second-order constant
coefficient subequations $\bbe$ and $\bbg$ to define a generalized equation $\bbh = \bbe\cap (-\wt\bbg)$
in $\rn$. 
 However, the main Theorem 2.6 holds for completely general 
subequations on  manifolds, as defined in [\DDR], and we give that result in this section.
For  general definitions we refer to [\DDR].  However, there are many
interesting cases  which the reader could keep in mind (without consulting [\DDR]), namely:
constant coefficient subequations $\bbe$ and $\bbg$ (not necessarily pure second-order) in $\rn$,  
variable coefficient subequations (constraint sets for subsolutions) on domains in $\rn$,
subequations on riemannian manifolds given canonically
by O($n$)-invariant equations in $\rn$,
subequations   on hermitian manifolds given canonically
by U($n$)-invariant equations in $\bbc^n$, etc.

Let $J^2(X)$  be the 2-jet bundle on a manifold $X$. (When $X=\rn$ this is just the bundle
$\rn\times(\bbr\oplus \bbr^n\oplus\Symn)$ over $\rn$ of order-2 Taylor expansions.)

\noindent
{\bf THEOREM 6.1} {\sl
Let $\O\ss\ss X$ be a domain in a manifold $X$,
and suppose $\bbe, \bbg\ss J^2(X)$ are two subequations.
Consider the  generalized equation $\bbh \equiv \bbe\cap(-\wt\bbg)$.

(a)   $\Int\, \bbh = \emptyset \quad\Rightarrow\quad$ uniqueness for the $\bbh$-(DP) on
$\O$, assuming that comparison holds for $\bbe$ and $\bbg$ on $\O$.

(b)   $\Int\, \bbh^* = \emptyset \quad\Rightarrow\quad$ existence holds for the $\bbh$-(DP) on
$\O$, assuming that existence for the $\bbe$-(DP) holds on $\O$.

(c)  There exists $h\in C^2(\O)\cap C(\ob)$ with $J_x^2 h \in \Int\,\bbh$ for all $x\in \O$
$\quad\Rightarrow\quad$ non-uniqueness for the $\bbh$-(DP) on $\O$ for the boundary values 
$\vf \equiv h\bigr|_{\bo}$.

(d)  There exists $f\in C^2(\O)\cap C(\ob)$ with $J_x^2 f \in \Int\,\bbh^*$ for all $x\in \O$
$\quad\Rightarrow\quad$ non-existence for the $\bbh$-(DP) on $\O$ for the boundary values
$\vf \equiv f\bigr|_{\bo}$, assuming that comparison holds for $\bbe$ and $\bbg$ on $\O$.
}

{\bf Proof of Assertion (a).}  We begin by noting that assertions (1.2) -- (1.5) hold for general subequations
as defined in [\DDR].  Our definition of $\bbh$ is the same as in Definition 2.2, and the assertion
(2.5) carries over.  As a result, Lemma 2.7 and Corollary 2.8 hold in this general case.
We now look at the proof of  Proposition 2.11, which carries over and says that under the assumption of comparison (C) for  both $\bbe$ and $\bbg$ we have Part (a).
\qed

{\bf Proof of Assertion (c).}  This follows exactly the argument given for Proposition 2.12. \qed

{\bf Proof of Assertion (b).}   This follows exactly the argument given for Proposition 2.14. \qed

{\bf Proof of Assertion (d).}  We are assuming that comparison (C) holds on $\O$ for both
$\bbe$ and $\bbg$.  This means that Proposition 2.16 holds, and therefore also Proposition 2.17
is valid.  This establishes Part (d). \qed

\Ex{6.2. (Generalized Constant Coefficient Equations in $\rn$)} 
Here  a  subequation is, by definition (cf. [\DDR], [\Survey]), a closed subset 
$$\bbf \ss J^2 \equiv \bbr \oplus \rn \oplus \Symn$$ such that $\bbf + (r, 0, P) \ss \bbf$
for $r\leq 0$ and $P\geq 0$ and such that $\bbf = \overline{\Int \, \bbf}$.
The topological condition $\bbf = \overline{\Int \, \bbf}$ was not part of our definition of a 
subequation in Section 1, which was for pure second-order subequations.
This allowed us to use the simpler definition that $\bbf$ is closed, since
  $\bbf = \overline{\Int \, \bbf}$ then follows easily from the positivity condition (1.1).

With regards to Assertion (a) comparison does not  hold for all such equations.  
However it does hold for many interesting classes, for instance, all gradient free ones.
Other such classes  can be found in [\CHLP]. 

On the other hand, {\bf  existence does hold for all these equations} $\bbf\ss J^2$, under the hypothesis that  the domain
$\O\ss\ss\rn$ has a smooth strictly $\bbf$ and $\wt\bbf$ convex boundary.  (See Theorem 12.7 in [\DDR].)
Now in Assertion (b) existence is only required for $\bbe$.  Therefore, Assertion (b) holds for  $\bbe, \bbg\ss J^2$ provided $\bo$ is strictly   $\bbe$ and $\wt\bbe$ convex.

\Ex{6.3. (Generalized Equations on an open set $X\ss \rn$)} 
The general  subequation here is a closed subset of the 2-jet bundle
$$\bbf \ss J^2(X)  \equiv  X \times (\bbr \oplus \rn \oplus \Symn)$$ 
such that
$$
\text{
$\bbf + (x; r, 0, P) \ss \bbf$ for $r\leq 0$ and $P\geq 0$ and for all $x\in X$,
}
$$
 $\bbf = \overline{\Int \, \bbf}$, and  for the fibres $\bbf_x$ we have
$$\bbf_x = \overline{\Int_x \, \bbf_x} \qquad \text{ and \qquad  $\Int_x\, \bbf_x = (\Int\, \bbf)_x$.}
$$
These are  barebones hypotheses needed for the constraint set for subsolutions
of  a nonlinear equation corresponding to $\partial \bbf$.

This is  the general case for domains $\O\ss\ss X\ss\rn$, and so the comparison and existence hypotheses in Theorem
6.1 need to be verified, but, of course, the literature is enormous.

For subequations on manifolds given by ``universal'' equations, much has been done
in [\DDR].  We shall now look at some cases.

\Ex{6.4. (Universal Subequations Defined on any Riemannian Manifold)} Let $\bbf\ss J^2$
be a subequation (as in Example 6.2) which is invariant under the natural action of the orthogonal
group O$(n)$  (or SO$(n)$).  Then $\bbf$ determines an invariant subequation $\bbf_X\ss J^2(X)$
on any riemannian (or oriented riemannian) manifold $X$ as follows.

Every $C^2$-function $u$ on $X$ has a riemannian hessian $\Hess\, u$, which is a section of 
the bundle $\Sym(X)$ of symmetric 2-forms on $X$, given at $x\in X$ by
$$
\left\{ \Hess_x u\right\} (V,W) \ \equiv\ V_x W \, u - (\nabla_V W)_x\, u
$$
where $\nabla$ is the Levi-Civita connection on $TX$.  Note that $\nabla_V W -\nabla_W V = VW-WV= [V,W]$,
so the symmetry and the tensorial properties of $\Hess\, u$ follow.

Now this riemannian Hessian gives a splitting of the 2-jet bundle
$$
J^2(X) \  \cong \ X\times (\bbr\oplus T^*X \oplus \Sym(X)),
$$
and the orthogonally invariant subequation $\bbf$ canonically determines 
a subequation $\bbf_X \ss  J^2(X)$ as follows.  Any  orthonormal frame field $e_1, ... , e_n$
for $TX$ on an open set $U\ss X$ determines an orthonormal framing of 
$J^2(U) \cong U\times (\bbr\oplus \rn\oplus \Symn)$. Via this framing, $\bbf$ determines 
a subequation on $U$.  However, if we use a different frame field $e_1', ... , e_n'$, the two framings
of $J^2(U)$ differ pointwise by O$(n)$-transforms.  By the  O$(n)$-invariance of $\bbf$
the subequation on $U$ are the same.   This also means that on two different open sets
$U,V \ss X$ the two subequations agree on $U\cap V$.  Hence, we have a well-defined global subequation
$\bbf_X \ss J^2(X)$.

For example, if $\bbf = \{(r,p,A) : \tr(A)\geq  0\}$, we get the subequation $\D u = \tr\{\Hess\, u\}\geq 0$
for the riemannian Laplacian.   If $\bbf = \{(r,p,A) : \det(A)\geq  0\}$, 
we get the real Monge-Amp\`ere subequation $\det\{\Hess\, u\}\geq 0$.
If $\bbf = \{(r,p,A) : p^t A p\geq  0\}$, one gets the infinite Laplacian on $X$.

The questions of comparison and of existence of solutions for the Dirichlet problem on manifolds are addressed 
 in [\DDR].  A cone subequation $M$ on $X$ is a cone monotonicity subequation for $\bbf_X$ if
$\bbf_X + M \ss \bbf_X$.  Then  for such equations we have the following from Thm. 13.2 and  Thm. 10.1 in [\DDR].
(See section 14 of [\DDR]  for examples.)

\Theorem{([\DDR])} {\sl Suppose $X$ admits a $C^2$ strictly $M$-subharmonic function.  Then comparison for $\bbf_X$ holds on
any domain $\O\ss\ss X$, and if $\bo$ is smooth and strictly $\bbf_X$ and $\wt\bbf_X$ convex, then existence holds
for the Dirichlet problem for all boundary functions $\vf \in C(\bo)$.
}

This construction has important generalizations.

\Ex{6.5. (Universal Subequations Defined on a Riemannian Manifold with $G$-Structure)} 
We now assume that the riemannian manifold $X$ can be covered by open sets $U$, with
an orthogonal tangent frame field $e^U \equiv (e_1, ... , e_n)$ on $U$, such that on the intersection
$U\cap V$ of two such, the change of frames from $e^U$ to $e^V$ always lies in a given compact subgroup
$G\ss {\rm O}(n)$.

For example, if $X^{2n}$ has an orthogonal almost complex structure $J$, then $X$ has an U$(n)$-structure.

If the euclidean subequation $\bbf$   is $G$-invariant, then the above construction gives a canonical subequation
on any riemannian manifold with $G$-structure.  For example,  for $(X^{2n}, J)$ above, we can define the 
complex Monge-Amp\`ere operator.

The Theorem at the end of Example 6.4 extends to these cases.

\def\GGG{{{\bf G \!\!\!\! l}}\ }

\Ex{6.6. (Geometric Cases)} 
Of particular importance are the {\sl geometric cases} 
given by a closed  subset $\GGG \ss G(p, \rn)$
of the Grassmannian of $p$-planes in $\rn$.  We assume that $\GGG$ is invariant
under a closed subgroup $G\ss O(n)$.  Then we consider the universal  euclidean subequation
$$
\bbf_\GGG \ \equiv \ \left \{(r,p,A) : \tr\left(A\bigr|_L\right)\geq 0 \ {\rm for \ all }\ L\in\GGG\right\}.
$$

This subequation now carries over to any riemannian manifold with $G$-structure.
For instance, suppose $\GGG$ is the set of special Lagrangian
$n$-planes in $\bbc^n$.  Then we get a subequation on any Calabi-Yau manifold $X$.  
If $\GGG$ is the set of associative 3-planes in $\bbr^7$, then we get a subequation
on any  7-manifold $X$ with holonomy $G_2$.

Theorems  in [\DDR] apply to these cases, but there is a better theorem in [\Geo].
We define the {\bf $\GGG$-core} of $X$ to be the set
$$
{\rm Core}_\GGG (X) \ \equiv\ \{x\in X : \text{no smooth strictly $\bbf_\GGG$-subharm. function is strict at $x$}\}
$$

\Theorem {([\Geo, Thm.\ 7.6 and Thm.\ 7.7])} {If  ${\rm Core}_\GGG (X) = \emptyset$,   then comparison for $\bbf_\GGG(X)$ holds on
any domain $\O\ss\ss X$, and if $\bo$ is smooth and strictly $\bbf_\GGG$  and $\wt{\bbf_\GGG}$ convex, then existence holds
for the Dirichlet problem for all boundary functions $\vf \in C(\bo)$.
}

\Remark {6.7}  In Section 2 we made a remark  which does not carry over to general 
subequations.  Finite intersections  of subequations are not always subequations.
There are classes of subequations where this is true (see [\CHLP]).
%that intersections and unions are again subequations in the same class (see [\CHLP]). 
However in general this  means that one  could expand  the definition of a generalized equation to cover many-fold intersections and unions of subequations.  This will be done elsewhere.

%%%%%%%%%%%%%%%%%%%%%%%%%%%%%%%%%%%%%%%%%%%%%%%%%%%%%
%%%%%%%%%%%%%%%%%%%%%%%%%%%%%%%%%%%%%%%%%%%%%%%%%%%%%
%%%%%%%%%%%%%%%%%%%%%%%%%%%%%%%%%%%%%%%%%%%%%%%%%%%%%
%%%%%%%%%%%%%%%%%%%%%%%%%%%%%%%%%%%%%%%%%%%%%%%%%%%%%
%%%%%%%%%%%%%%%%%%%%%%%%%%%%%%%%%%%%%%%%%%%%%%%%%%%%%
%%%%%%%%%%%%%%%%%%%%%%%%%%%%%%%%%%%%%%%%%%%%%%%%%%%%%

\vfill\eject

\centerline{\headfont Appendix  A.  The Quasi-Convexity Characterization of C$^{1,1}$.}
 \medskip 
Interestingly, the condition that a function $u$ be locally $C^{1,1}$ is 
 equivalent to $u$  locally being simultaneously quasi-convex and quasi-concave.
 This was probably first observed by Hiriart-Urruty and Plazanet in [\HP]. An alternate proof appeared in 
 Eberhard  [\Eb] and also in [\HU].
 For the benefit of the reader we include a proof here.

We say that a function is $\l-C^{1,1}$ if it is $C^1$ and
 the first derivative is locally Lipschitz with Lipschitz coefficient $\l$.
 
\Theorem {A.1}
$$
u \ \ {\rm is}  \ \l-C^{1,1} 
\qquad\iff\qquad
{\rm both}\ \ 
\pm u \ \ {\rm are}\ {\rm locally} \  \l-{\rm quasi-convex}
$$

\pf We consider $u$ on a convex set $\O$.
 
($\Rightarrow$)  Suppose that $u$ is $\l$-$C^{1,1}$, i.e., $u\in C^1$ and $|D_xu -D_yu|\leq\l |x-y|$
for all $x,y\in\O$.  Set $f\equiv u+{\l\over 2}|x|^2$. Then 
$$
D_xf -D_yf\ =\ \l(x-y) + D_xu -D_yu,
$$
and hence
$$
\begin{aligned}
\bra{D_xf -D_yf}{x-y} \ &=\ \l|x-y|^2 + \bra{D_xu -D_yu}{x-y}   \cr
&\geq \ \l|x-y|^2  - |D_xu -D_yu||x-y|   \cr
&= \ (\l|x-y|  - |D_xu -D_yu|)  |x-y|   \ \geq\ 0. 
\end{aligned}
$$
This form of monotonicity of $Df$ is one of the standard definitions of $f$ being convex.
The same proof works for $-u$

($\Leftarrow$ ) We state the converse as a proposition.

\Prop{A.2}
{\sl
If $u$ and $-u$ are $\l$-quasi-convex on a convex domain $\O$, then $u\in C^{1}$ and}
$$
|D_xu -D_yu|\ \leq \ \l|x-y|, \  \ {\rm i.e., \ } \ u\  {\rm is\ \ } \l- C^{1,1}, \ \ {\rm on}\ \  \O.
$$
\noindent
{\bf Proof.}
We first show that this is true if $u\in C^\infty$. 
Note that  $\pm u$ are $\l$-quasi-convex 
\ \ $\iff$\ \ $D^2_x u +\l I\geq0$ and  $-D^2_x u +\l I\geq0$  for all $x$
\ \ $\iff$\ \ $-\l I \leq D^2_x u \leq \l I$ for all $x\in\O$.  Fix $x,y\in\O$.  By the Mean Value Inequality in [\RU],
$|D_xu -D_yu| \leq |D^2_\x u| |x-y|$ for some $\x\in [x,y]$, and hence
$|D_xu -D_yu| \leq \l|x-y|$.  (Here $|A| \equiv \sup\{|A(e)| :|e|=1\}$ is the operator norm of the symmetric transformation
$A= D^2_\x u$, which is equal to the max of $|\bra{Ae}e|$ over unit vectors $e$.)

In general, since the graph of $u+{\l\over 2}|x|^2$ has a supporting hyperplane from below and
 the graph of  $u-{\l\over 2}|x|^2$ has a supporting hyperplane from above, at every point, 
 the function $u$ is differentiable everywhere.  By partial continuity of the first derivative
 for quasi-convex functions (see for example Lemma 1.3 in [\NOTES]), we have $u\in C^1$.
 
 Now standard convolution $u^\e \equiv u * \vf_\e$ works just fine to complete the proof
 since $\pm u^\e$ is $\l$-quasi-convex by the next lemma, and the fact that
  $u\in C^1 \ \Rightarrow \ Du^\e \to Du$ locally uniformly. \qed

\Lemma{A.3}
{\sl
$u$ is $\l$-quasi-convex \ \ $\Rightarrow$\ \ $u^\e \equiv u * \vf_\e$ is  $\l$-quasi-convex.
}

\noindent
{\bf Proof.}
Suppose $u$ is  $\l$-quasi-convex, i.e., $f\equiv u+{\l\over 2}|x|^2$ is convex.
Standard convolution of $f$ with an approximate identity $\vf_\e$ based on $\vf$
(i.e., $\vf_\e(x) \equiv {1\over \e^n}\vf({x\over \e})$) yields $f^\e \equiv f * \vf_\e$ smooth and convex.
Note that $(|x|^2* \vf_\e) = |x|^2 +\bra ax+c$ preserves $|x|^2$ modulo an affine function,
since $\int|x+\e y|^2 \vf(y)\,dy = |x|^2+\e\bra ax+C\e^2$
where $\bra ax = 2\int \bra xy \vf(y) \, dy $ and
$C =\int|y|^2\vf(y)\,dy$.
Therefore, $D^2 f^\e = D^2 u^\e +\l I$, proving that each $u^\e$ is $\l$-quasi-convex.\qed

%\vfill\eject
\vskip.3in

\centerline{\bf REFERENCES}

\noindent
[\CHLP] 
 M. Cirant, F. R. Harvey,  H. B. Lawson, Jr. and K. R. Payne, {\sl  Comparison principles by monotonicity and
 duality for constant coefficient nonlinear potential theory and pde's}.

\noindent
[\CIL]   M. G. Crandall, H. Ishii and P. L. Lions {\sl
User's guide to viscosity solutions of second order partial differential equations},  
Bull. Amer. Math. Soc. (N. S.) {\bf 27} (1992), 1-67.

\noindent
[\CRA]  M. G. Crandall,  {\sl  Viscosity solutions: a primer},  
pp. 1-43 in ``Viscosity Solutions and Applications''  Ed.'s Dolcetta and Lions, 
SLNM {\bf 1660}, Springer Press, New York, 1997.

\noindent
[\DHGL]  A. Daniilidis, M. Haddou, E. Le Gruyer, and O. Ley, {\sl  
Explicit formulas for $C^{1,1}$ Glaesner-Whitney extensions of 1-Taylor fields in Hilbert spaces,}
Proc. A. M. S.,  {\bf  146}, no. 10 (2018), 4487-4495.

\noindent
 {[\Eb]}  A. Eberhard, {\sl Prox-regularity and subjets}, in: A. Rubinov (Ed.), ÔOptimization and Related TopicsÕ, Applied Optimization Volumes, Kluwer Academic Publishers, Dordrecht, 2001, pp. 237-313.

\noindent
 {[\PTCG]}  F. R. Harvey and H. B. Lawson, Jr,   {\sl  An introduction to potential theory in calibrated geometry},
    Amer. J. Math.  {\bf 131} no. 4 (2009), 893-944.  ArXiv:math.0710.3920.

\noindent
 {[\PTCGG]} 
 ----------- ,       {\sl  Duality of positive currents and plurisubharmonic functions in calibrated geometry},
   Amer. J. Math.    {\bf 131} no. 5 (2009), 1211-1240. ArXiv:math.0710.3921.

\noindent
[\DD] 
 ----------- ,     {\sl  Dirichlet duality and the non-linear Dirichlet problem},
    Comm. on Pure and Applied Math. {\bf 62} (2009), 396-443.

\noindent
[\DDR] 
----------- ,    {\sl   Dirichlet duality and the non-linear Dirichlet problem on Riemannian manifolds}, J. Diff. Geom.  {\bf 88} No. 3 (2011), 395-482.  ArXiv:0907.1981.

\noindent
[\Survey] 
----------- ,   {\sl  Existence, uniqueness and removable singularities
for nonlinear partial differential equations in geometry},
  pp. 102-156 in ``Surveys in Differential Geometry 2013'', vol. 18,  
H.-D. Cao and S.-T. Yau eds., International Press, Somerville, MA, 2013.
ArXiv:1303.1117.

\noindent
[\Hyp] 
----------- , {\sl  Hyperbolic polynomials and the Dirichlet problem},    ArXiv:0912.5220.

\noindent
[\Geo] 
----------- ,  {\sl  Geometric plurisubharmonicity and convexity - an introduction},  Advances in Math.  {\bf 230} (2012), 2428-2456.   ArXiv:1111.3875.

\noindent
[\NOTES]     \ \----------,   {\sl Notes on the differentiation of quasi-convex functions},     
ArXiv:1309.1772.

\noindent
[\Gar] 
----------- , {\sl  G\aa rding's theory of hyperbolic polynomials},  {\sl Communications in Pure and Applied Mathematics}  {\bf 66} no. 7 (2013), 1102-1128.

\noindent
[\Pot] 
----------- ,  {\sl  Potential theory on almost complex manifolds},    Ann. Inst.  Fourier,  Vol. 65 no. 1 (2015), p. 171-210.  ArXiv:1107.2584.

\noindent
[\Lag] 
----------- ,  {\sl  Lagrangian potential theory and a Lagrangian equation of  Monge-Amp\`ere type}.   ArXiv:1712.03525.

 \noindent
[\SMP] 
----------- ,   {\sl  Characterizing the strong maximum principle for constant coefficient subequations}, Rendiconti di Matematica  {\bf 37} (2016), 63-104.  ArXiv:1303.1738.

 \noindent
[\IDP]    
----------- ,   {\sl The inhomogeneous Dirichlet problem for natural operators on manifolds}.   ArXiv:1805.11121.

 \noindent
[\Edge]    
----------- ,   {\sl  Pluriharmonics in general potential theories}.  ArXiv:1712.03447.

\noindent
 {[\HP]} J.-B. Hiriart-Urruty, Ph. Plazanet, {\sl MoreauÕs Theorem revisited}, Analyse Nonlin\'eaire (Perpignan, 1987), Ann. Inst. H. Poincar\'e {\bf 6} (1989), 325-338.

 \noindent
 [\HU]  J.-B. Hiriart-Urruty, {\sl How to regularize a difference of convex functions}, J. of Math. Anal. and Appls., 
 {\bf 162} (1991), 196-209.

   \noindent
[\Kry]    N. V. Krylov,    {\sl  On the general notion of fully nonlinear second-order elliptic equations},    Trans. Amer. Math. Soc.   {\bf  347}   (1995), 857-895.

\noindent
[\GRU]  E. Le Gruyer,  {\sl Minimal Lipschitz extensions to differentiable functions defined on a
Hilbert space}, Geom. Funct. Anal.  {\bf 19} (2009), no. 4, 1101Ð1118.

 \noindent
[\ST]    
J. Streets and  G. Tian,  {\sl A parabolic flow of pluriclosed metrics}, Int. Math. Res. Not. {\bf 16} (2010), 3101-3133.

 \noindent
[\SW]    
J. Streets and  M. Warren,  {\sl  Evans-Krylov estimates for a non-convex Monge-Amp\`ere equation},
ArXiv:1410.2911.

 \noindent
[\RU]    W, Rudin,  {\sl Principles of Mathematical Analysis}, McGraw-Hill,   New York, NY,    1953.

\end{document}